%% file: modular.tex
\documentclass[centertags,leqno]{article}
\usepackage{def}
\title{Reduction of covers and Hurwitz spaces}
\author{Irene I.\ Bouw and Stefan Wewers\\[1ex]
        University of Pennsylvania}
\date{}

\begin{document}

\maketitle

\begin{abstract}
  In this paper we study the reduction of Galois covers of curves,
  from characteristic zero to positive characteristic. The starting
  point is a recent result of Raynaud, which gives a criterion for
  good reduction for covers of the projective line branched at three
  points. We use the ideas of Raynaud to study the case of covers of
  the projective line branched at four points. Under some condition on
  the Galois group, we generalize the criterion for good reduction of
  Raynaud. As a new ingredient, we use the Hurwitz space of such
  covers. Combining our results on reduction of covers with the
  Hurwitz space approach, we are able to describe the reduction of the
  Hurwitz space modulo $p$ and compute the number of covers with
  good reduction. \\[2ex]
\end{abstract}
 
\vspace{1cm}

\input{intro}

\input{sec1}

%\input{sec2}

\input{sec3}

\input{sec4}

\input{sec5}

\input{examples}

\bibliographystyle{abbrv} \bibliography{../hurwitz}

\vspace{5mm}
\flushright{DRL, University of Pennsylvania\\
            209 South 33rd Street\\
            Philadelphia, PA 19104-6395\\
            bouw@math.upenn.edu\\
            wewers@math.upenn.edu}

\end{document}

%% file: intro.tex
%\documentclass{article}
%\usepackage{def}

%\begin{document}

%\hspace*{\fill}
%\begin{minipage}[b]{55mm}
%  There is nothing either good or bad,\\
%  But thinking makes it so.\\[2ex]
%  \hspace*{\fill} Shakespeare, \sl Hamlet
%\end{minipage}

\section*{Introduction}

Let $R$ be a complete discrete valuation ring whose residue field $k$
is algebraically closed of characteristic $p$ and whose quotient field
$K$ is of characteristic zero.  Let $f_K\!:Y_K\to\PP^1_K$ be a Galois
cover defined over $K$. We ask ourselves whether $f_K$ has good
reduction. In case $p$ divides the order of the Galois group, this
is a hard question. We cannot expect good reduction, in general.
 
A recent paper of Raynaud gives a criterion for good reduction in a
first case. Let $f_K\!:Y_K\to\PP^1_K$ be a Galois cover branched at
$0$, $1$ and $\infty$. Suppose that $p$ strictly divides the order of
the Galois group $G$, but not the ramification indices of $f_K$ and
that the center of $G$ is trivial. Let $P\cong\ZZ/p$ be a $p$-Sylow of
$G$, and denote by $n:=[N_G(P):C_G(P)]$ the index of the centralizer
of $P$ in the normalizer of $P$. Let $e$ be the absolute ramification
index of $p$ in $K$. It is shown that the cover has good reduction,
provided that $en< p-1$. The idea of the proof is to suppose that
$f_K$ has bad reduction and to study the {\em semistable reduction} of
$f_K$. Raynaud proves general structure results on the semistable
reduction. Using these techniques, he proves that bad reduction
implies $en\geq p-1$.

In this paper, we follow the approach of Raynaud. We consider the
reduction of $G$-covers $f_K\!:Y_K\to\PP^1_K$ branched at four points.
We suppose that $\,p||\,|G|$, but that $p$ does not divide the
ramification indices of $f_K$. This is the next case to study after
the result of Raynaud. However, we put a much stronger condition on
the group $G$. In particular, we assume that $n=2$. The criterion for
good reduction given by Raynaud extends to our situation. The stronger
condition on the Galois group allows us to get a stronger result. For
instance, we are able to describe the semistable model of $f_K$.

When passing from $3$ to $4$ branch points, a new aspect arises: the
reduction of $f_K$ might depend on the position of the branch points.
It is therefore natural to study the corresponding {\em Hurwitz
  space}, i.e.\ the moduli space of $G$-covers of a certain type, and
its reduction to positive characteristic. Hurwitz spaces were first
introduced in a purely geometric context, but have since then proved
to be useful for studying arithmetic aspects of covers as well, see
e.g.\ \cite{Fried87}, \cite{FriedVoe91}. There are many variants of
Hurwitz spaces. To fix ideas, let $G$ be a group, $r\geq 3$ and denote
by $H:=H_r^{\rm in}(G)$ the Hurwitz space parameterizing $G$-Galois
covers of $\PP^1$ with $r$ branch points. Then $H$ is a smooth variety
defined over $\QQ$. It is known that $H$ has good reduction to
characteristic $p$ provided $p\nmid |G|$.

Recently, Abramovich and Oort \cite{AbrOort98} suggested a definition
of an arithmetic compactification of Hurwitz spaces. We follow, and
somewhat simplify, this approach. The idea is to take the closure of
$H$ inside a bigger moduli space which parameterizes maps between
stably marked curves. We obtain an algebraic space $\hb$ which is
proper over $\ZZ$ and contains $H$ as a dense open subspace. The
complement $\hb\bad:=\hb-H$ is a closed subspace supported in positive
characteristic and corresponds to $G$-covers with bad reduction. Let
us denote by $\hb\good$ the closure of $H\otimes\FF_p$ inside
$\hb\otimes\FF_p$. In the case of $r=4$ branch points, there is a
finite map $\hb\good\to\PP^1_\lambda$. Let $d\good$ be its degree. For
a ``generic'' choice of $\lambda\in K-\{0,1\}$, there will be exactly
$d\good$ nonisomorphic $G$-covers $f_K:Y_K\to\PP^1_K$ branched in
$0,1,\infty$ and $\lambda$ which have good reduction. More precisely,
there is a finite set $\lambdab_1,\ldots,\lambdab_m\in\FFbp$ of
exceptional values such that $d\good$ is the number of covers as above
with good reduction provided $\lambda\not\equiv\lambdab_i\pmod{p}$ for
all $i$.  The $\lambda_i$'s are the images of the points where
$\hb\good$ intersects $\hb\bad$.  One way to actually compute the
number $d\good$ and the exceptional values $\lambdab_i$ would be to
determine the precise structure of $\hb\bad$.  This seems a difficult
problem, in general. In this paper, we are able to solve it in a
special case, thanks to a surprising connection with modular curves.

\bigskip\noindent Even though our definition of $\hb$ works without
any extra assumption, it seems to be a hard problem to describe the
structure of $\hb$ and $\hb\bad$ in general. To our knowledge, the
only case that has been studied is the case of modular curves, which
admit an interpretation as Hurwitz spaces. Let us sketch this
correspondence in the case of $X_1(p)$. Let $f_K:Y_K\to\PP^1_K$ be a
$G$-cover branched at $4$ points of order $2$, where $G$ is a dihedral
group of order $2p$. We may identify $Y_K$ with an elliptic curve
$E_K$ over $K$; the cover $f_K$ factors through a $p$-cyclic isogeny
$\pi_K:E_K\to E_k'$.  Moreover, the choice of an element $\sigma$ of
$G$ of order $p$ defines a $p$-torsion point $P:=\sigma(0)\in E_K[p]$
generating the kernel of $\pi_K$. The pair $(E_K,P)$ corresponds to a
point on $X_1(p)$. This gives an identification $H\cong X_1(p)$, for a
suitable Hurwitz space $H$. By the results of Katz and Mazur
\cite{KatzMazur}, the reduction of $X_1(p)$ to characteristic $p$ is
well understood.  One can check that the subspace $\hb\bad$ of the
arithmetic compactification of $H$ corresponds to the component of
$X_1(p)\otimes\FF_p$ parameterizing pairs $(E,P)$ such that $P=0$ and
hence the isogeny $\pi:E\to E'=E/\gen{P}$ is inseparable. Even without
using the very precise results of \cite{KatzMazur}, the theory of
elliptic curves gives the following result on good reduction of Galois
covers.  Suppose the elliptic curve $E'_K$ given by the equation
$y^2=x(x-1)(x-\lambda)$ has good ordinary reduction. Then there are
precisely $p-1$ nonisomorphic $G$-covers branched in $0$, $1$,
$\infty$ and $\lambda$, with ramification index $2$, which have good
reduction. On the other hand, there is no such cover with good
reduction if $E'_K$ has supersingular reduction.

\subsection*{Results}

In this paper, we look at the following situation. Let $G$ be a finite
group and $p$ an odd prime which strictly divides $|G|$. We assume
that the normalizer of a $p$-Sylow of $G$ is a dihedral group. Let $K$
be as in the beginning, and let $f_K:Y_K\to\PP^1_K$ be a $G$-Galois
cover branched at $4$ points, of order prime-to-$p$. We prove that the
cover $f_K$ has either good reduction or a very specific type of bad
reduction, which we call {\sl modular reduction}.

We will briefly explain what this means. Assume that $f_K$ has bad
reduction. Following Raynaud \cite{Raynaud98}, \S 3.2, we associate to
the $G$-cover $f_K$ a $\Delta$-cover $g_K:Z_K\to\PP^1_K$, called the
{\em auxiliary cover}. Here $\Delta$ is a subgroup of $G$, and $g_K$
is branched in the same points as $f_K$ and has bad reduction. The
statement that $f_K$ has modular reduction of {\em level} $N$ means
essentially that $\Delta$ is a dihedral group of order $2N$ (where
$p|N$) and that $g_K$ has ramification of order $2$. In particular,
$g_K$ gives rise to a $K$-point on $X_1(N)$. This is the link between
our results and modular curves.

To explain the construction of $g_K$, we assume for simplicity that
the normalizer of a $p$-Sylow of $G$ is of order $2p$ and that the
branch points of $f_K$ do not coalesce modulo $p$. Let $f:Y\to X$ be
the special fiber of the semistable model of $f_K$. The curve $X$
consists of 5 components: the strict transform of the original
component $X_0$, and 4 tails $X_1,\ldots,X_4$ containing the
specializations of the branch points $x_i$. The cover $f$ is
inseparable over $X_0$ and separable over the tails. Let $E$ be a
component of $Y$ above $X_0$. The decomposition group 
$\Delta:=D(E)\subset G$ is dihedral of order $2p$, the inertia group
$I(E)$ cyclic of order $p$ (see Fig.\ \ref{modauxfig} in Section
\ref{figpage}). We obtain $g:Z\to X$ by replacing, for $i=1,\ldots,4$,
the (disconnected) $G$-cover $f^{-1}(X_i)\to X_i$ by a $\Delta$-cover
$Z_i\to X_i$ which is locally, i.e.\ in an \'etale neighborhood of
$E$, isomorphic to $f^{-1}(X_i)\to X_i$ and tamely ramified above
$x_i\in X_i$. This is possible in a unique way, by the Katz--Gabber
Lemma,
\cite{Katz86}. Using formal patching one can show that $g:Z\to X$ is
the reduction of a $\Delta$-cover $g_K:Z_K\to\PP^1_K$ which, in some
sense, contains all the information about the bad reduction of $f_K$.

\begin{figure}[htb]
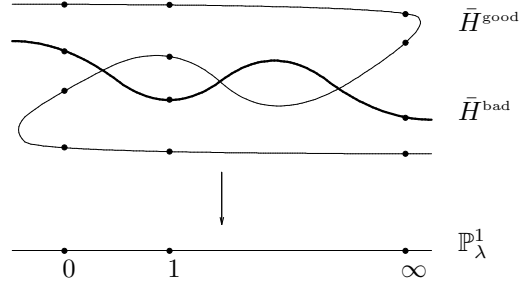

\begin{center}
\include{goodbad}
\end{center}
\caption{\label{goodbadfig} the Reduction Theorem}
\end{figure}

Let $H\subset H_4\inn(G)$ be the subset of the Hurwitz space
corresponding to $G$-covers with $4$ branch points and prime-to-$p$
ramification. Denote by $\hb$ its arithmetic compactification and by
$\hb\bad\subset\hb\otimes\FF_p$ the subspace corresponding to bad
reduction in characteristic $p$. The map $f:Y\to X$ discussed above
corresponds to a $k$-point on $\hb\bad$; the associated map $g:Z\to X$
corresponds essentially to a $k$-point on $X_1(N)\otimes\FF_p$. The
formal patching argument mentioned above can be used to show that the
deformation theory of $f$ and $g$ are equivalent. This gives a strong
connection between the subspace $\hb\bad\subset\hb$ and the component
of $X_1(N)\otimes\FF_p$ corresponding to inseparable isogenies. Using
the results of \cite{KatzMazur} on the reduction of $X_1(N)$, we prove
our Reduction Theorem, which describes $\hb\otimes\FF_p$. It can be
roughly stated as follows (see Fig.\
\ref{goodbadfig}). The subspaces $\hb\good$ and $\hb\bad$ are smooth
curves over $\FF_p$; they intersect transversally in the {\em
  supersingular points}, i.e.\ the points of $\hb\bad$ with
supersingular $\lambda$-value. The scheme $\hb\otimes\FF_p$ is not
reduced, in general; the irreducible components of $\hb\bad$ have
multiplicity $p-1$ or $(p-1)/2$. Moreover, each irreducible component
of $\hb\bad$ is essentially the reduction of a modular curve. 
\vspace{5mm}

We call the $\QQb$-rational points of $\hb$ lying above $0$, $1$ and
$\infty$ the {\em cusps} of $H$. A cusp corresponds to a degenerate
cover $f_{\QQb}:Y_{\QQb}\to X_{\QQb}$, obtained from a smooth
$G$-cover by coalescing of branch points. The curve $X_{\QQb}$ is the
union of two projective lines intersecting in one point. The cover
$f_{\QQb}$ is ramified above the singular point of $X_{\QQb}$, say of
order $n$. Our next result, which we call the Cusp Principle, states
that this cusp has bad reduction (i.e.\ its closure in $\hb$ meets
$\hb\bad$) if and only if $p|n$.  Essentially, this is an application
of Raynaud's result, since the degenerate cover $f_{\QQb}$ is build up
from two $3$ point covers. The Cusp Principle is the key result in our
calculation of the number of $G$-covers with good reduction. The point
here is that, via the Hurwitz classification and the braid action, one
can explicitly compute the set of cusps of $H$ and decide for each of
them whether they have good or bad reduction. In particular, given a
finite group $G$ and an odd prime $p$ verifying all the assumptions
made above, we can compute two numbers, $d$ and $d\bad$, such that
\[    \renewcommand{\arraystretch}{1.5}
     |\mathop{\rm Cov}(G,\lambda)\good| \;=\; \left\{\;\;
        \begin{array}{ll} 
           d- d\bad,\qquad & 
              \text{\rm if $\lambda$ is ordinary,} \\
           d-\frac{p+1}{p}d\bad,\qquad & 
              \text{\rm if $\lambda$ is supersingular.}
        \end{array}
     \right.
\]
Here $\mathop{\rm Cov}(G,\lambda)\good$ is the set of isomorphism
classes of $G$-covers of $\PP^1$ with good reduction, branched in $0$,
$1$, $\infty$ and $\lambda$.

Let $G=PSL_2(\ell)$, where $\ell$ is an odd prime different from $p$
such that $p$ exactly divides the order of $G$, and consider
$G$-covers branched in $4$ points of order $\ell$. We have computed
the cusps of the corresponding Hurwitz spaces, using the computer
program {\em ho} \cite{ho}, for $\ell\leq 31$. From this information,
we can deduce the complete structure of $\bar{H}^{\rm bad}\otimes
\FF_p$ and the number of covers with good reduction.

\bigskip\noindent This paper owes a lot to Raynaud. It started as an
attempt to understand a talk he gave in Oberwolfach, June 1997. In
this talk Raynaud presented Example \ref{Raynaudexa}. In the problem
session of the same conference, he gave a similar problem as an
exercise. In a way, this paper is our solution of this exercise. We
would also like to thank Bas Edixhoven for a helpful conversation and
for sending his manuscript \cite{Edix}, and Andrew Kresch for comments
on an earlier version of Section \ref{complete}.  The second author
gratefully acknowledges financial support from the Deutsche
Forschungsgemeinschaft.

\subsection*{Notation}

In this paper we will understand by a {\sl semistable curve} a flat
projective morphism $X\to S$ of schemes whose geometric fibers are
reduced connected curves having at most ordinary double points as
singularities. We write $X\sm$ for the subset of smooth points of the
morphism $X\to S$. A {\em mark} on $X/S$ is a closed subscheme
$D\subset X\sm$ which is finite and \'etale over $S$. The pair
$(X/S,D)$ is called a marked semistable curve. A {\em stably marked
  curve} is either a pointed stable curve $(X/S;x_i)$ in the sense of
\cite{Knudsen83} or a marked semistable curve $(X/S,D)$ which becomes
a pointed stable curve after an \'etale base change $S'\to S$. By an
{\em algebraic stack} we mean an algebraic stack in the sense of
Deligne--Mumford \cite{DelMum69}.

%\bibliographystyle{abbrv}
%\bibliography{../hurwitz}

%\end{document}

%% file: goodbad.tex
%\documentclass{article}
%\usepackage{def}

%\begin{document}

\setlength{\unitlength}{0.0003in}
{\renewcommand{\dashlinestretch}{30}
\begin{picture}(8172,4395)(0,-10)
\put(915,97){\circle*{90}}
\put(915,1870){\circle*{90}}
\put(915,2842){\circle*{90}}
\put(915,3520){\circle*{90}}
\put(915,4327){\circle*{90}}
\put(2715,97){\circle*{90}}
\put(2715,1800){\circle*{90}}
\put(2715,2690){\circle*{90}}
\put(2715,3427){\circle*{90}}
\put(2715,4320){\circle*{90}}
\put(6765,97){\circle*{90}}
\put(6765,1762){\circle*{90}}
\put(6765,2380){\circle*{90}}
\put(6765,3670){\circle*{90}}
\put(6765,4160){\circle*{90}}
\drawline(3615,1447)(3615,547)
\drawline(3615,1447)(3615,547)
\drawline(3585.000,667.000)(3615.000,547.000)(3645.000,667.000)
\drawline(15,97)(7215,97)
\put(880,-350){$0$}
\put(2680,-350){$1$}
\put(6680,-350){$\infty$}

\thicklines
\drawline(15,3697)
   (71.250,3696.341)
   (127.500,3694.363)
   (183.750,3691.067)
   (240.000,3686.453)
   (296.250,3680.521)
   (352.500,3673.270)
   (408.750,3664.700)
   (465.000,3654.812)
   (521.250,3643.606)
   (577.500,3631.082)
   (633.750,3617.239)
   (690.000,3602.078)
   (746.250,3585.599)
   (802.500,3567.801)
   (858.750,3548.685)
   (915.000,3528.250)
   (971.250,3506.497)
   (1027.500,3483.426)
   (1083.750,3459.036)
   (1140.000,3433.328)
   (1196.250,3406.302)
   (1252.500,3377.957)
   (1308.750,3348.294)
   (1365.000,3317.312)
   (1421.250,3285.013)
   (1477.500,3251.395)
   (1533.750,3216.458)
   (1590.000,3180.203)
   (1646.250,3142.630)
   (1702.500,3103.738)
   (1758.750,3063.528)
   (1815.000,3022.000)
   (1871.250,2981.131)
   (1927.500,2942.898)
   (1983.750,2907.303)
   (2040.000,2874.344)
   (2096.250,2844.021)
   (2152.500,2816.336)
   (2208.750,2791.287)
   (2265.000,2768.875)
   (2321.250,2749.100)
   (2377.500,2731.961)
   (2433.750,2717.459)
   (2490.000,2705.594)
   (2546.250,2696.365)
   (2602.500,2689.773)
   (2658.750,2685.818)
   (2715.000,2684.500)
   (2771.250,2685.818)
   (2827.500,2689.773)
   (2883.750,2696.365)
   (2940.000,2705.594)
   (2996.250,2717.459)
   (3052.500,2731.961)
   (3108.750,2749.100)
   (3165.000,2768.875)
   (3221.250,2791.287)
   (3277.500,2816.336)
   (3333.750,2844.021)
   (3390.000,2874.344)
   (3446.250,2907.303)
   (3502.500,2942.898)
   (3558.750,2981.131)
   (3615.000,3022.000)
   (3671.250,3062.869)
   (3727.500,3101.102)
   (3783.750,3136.697)
   (3840.000,3169.656)
   (3896.250,3199.979)
   (3952.500,3227.664)
   (4008.750,3252.713)
   (4065.000,3275.125)
   (4121.250,3294.900)
   (4177.500,3312.039)
   (4233.750,3326.541)
   (4290.000,3338.406)
   (4346.250,3347.635)
   (4402.500,3354.227)
   (4458.750,3358.182)
   (4515.000,3359.500)
   (4571.250,3358.182)
   (4627.500,3354.227)
   (4683.750,3347.635)
   (4740.000,3338.406)
   (4796.250,3326.541)
   (4852.500,3312.039)
   (4908.750,3294.900)
   (4965.000,3275.125)
   (5021.250,3252.713)
   (5077.500,3227.664)
   (5133.750,3199.979)
   (5190.000,3169.656)
   (5246.250,3136.697)
   (5302.500,3101.102)
   (5358.750,3062.869)
   (5415.000,3022.000)
   (5471.250,2980.472)
   (5527.500,2940.262)
   (5583.750,2901.370)
   (5640.000,2863.797)
   (5696.250,2827.542)
   (5752.500,2792.605)
   (5808.750,2758.987)
   (5865.000,2726.688)
   (5921.250,2695.706)
   (5977.500,2666.043)
   (6033.750,2637.698)
   (6090.000,2610.672)
   (6146.250,2584.964)
   (6202.500,2560.574)
   (6258.750,2537.503)
   (6315.000,2515.750)
   (6371.250,2495.315)
   (6427.500,2476.199)
   (6483.750,2458.401)
   (6540.000,2441.922)
   (6596.250,2426.761)
   (6652.500,2412.918)
   (6708.750,2400.394)
   (6765.000,2389.188)
   (6821.250,2379.300)
   (6877.500,2370.730)
   (6933.750,2363.479)
   (6990.000,2357.547)
   (7046.250,2352.933)
   (7102.500,2349.637)
   (7158.750,2347.659)
   (7215.000,2347.000)

\thinlines
\drawline(15,4327)
   (67.635,4326.997)
   (120.070,4326.986)
   (172.307,4326.969)
   (224.344,4326.945)
   (276.183,4326.914)
   (327.822,4326.876)
   (379.262,4326.832)
   (430.503,4326.780)
   (481.545,4326.722)
   (532.387,4326.657)
   (583.031,4326.585)
   (633.475,4326.506)
   (683.721,4326.420)
   (733.767,4326.327)
   (783.614,4326.228)
   (833.262,4326.121)
   (882.710,4326.008)
   (931.960,4325.888)
   (981.011,4325.761)
   (1029.862,4325.627)
   (1078.514,4325.486)
   (1126.967,4325.338)
   (1175.221,4325.184)
   (1223.276,4325.022)
   (1271.132,4324.854)
   (1318.789,4324.679)
   (1366.246,4324.497)
   (1413.505,4324.308)
   (1460.564,4324.113)
   (1507.424,4323.910)
   (1554.085,4323.701)
   (1600.547,4323.484)
   (1646.810,4323.261)
   (1692.873,4323.031)
   (1738.738,4322.794)
   (1784.403,4322.551)
   (1829.869,4322.300)
   (1875.136,4322.042)
   (1920.204,4321.778)
   (1965.073,4321.507)
   (2009.743,4321.229)
   (2054.214,4320.944)
   (2098.485,4320.652)
   (2142.557,4320.353)
   (2186.431,4320.048)
   (2230.105,4319.735)
   (2273.580,4319.416)
   (2316.855,4319.090)
   (2359.932,4318.757)
   (2402.810,4318.417)
   (2445.488,4318.070)
   (2487.968,4317.717)
   (2530.248,4317.356)
   (2572.329,4316.989)
   (2614.211,4316.614)
   (2655.894,4316.233)
   (2697.377,4315.845)
   (2738.662,4315.451)
   (2779.747,4315.049)
   (2820.634,4314.640)
   (2861.321,4314.225)
   (2901.809,4313.803)
   (2942.098,4313.374)
   (2982.188,4312.938)
   (3022.078,4312.495)
   (3061.770,4312.045)
   (3101.262,4311.588)
   (3140.555,4311.125)
   (3179.650,4310.654)
   (3218.545,4310.177)
   (3257.241,4309.693)
   (3295.737,4309.202)
   (3334.035,4308.704)
   (3372.133,4308.200)
   (3410.033,4307.688)
   (3447.733,4307.170)
   (3485.234,4306.644)
   (3522.536,4306.112)
   (3559.639,4305.573)
   (3596.543,4305.027)
   (3633.248,4304.475)
   (3669.753,4303.915)
   (3706.059,4303.348)
   (3742.167,4302.775)
   (3778.075,4302.195)
   (3813.784,4301.608)
   (3884.604,4300.413)
   (3954.629,4299.191)
   (4023.856,4297.941)
   (4092.287,4296.664)
   (4159.922,4295.359)
   (4226.760,4294.027)
   (4292.802,4292.668)
   (4358.047,4291.281)
   (4422.495,4289.866)
   (4486.147,4288.424)
   (4549.003,4286.955)
   (4611.062,4285.458)
   (4672.324,4283.934)
   (4732.790,4282.382)
   (4792.460,4280.802)
   (4851.333,4279.196)
   (4909.409,4277.562)
   (4966.689,4275.900)
   (5023.173,4274.211)
   (5078.860,4272.494)
   (5133.750,4270.750)
   (5187.844,4268.978)
   (5241.141,4267.179)
   (5293.642,4265.353)
   (5345.347,4263.499)
   (5396.255,4261.618)
   (5446.366,4259.709)
   (5495.681,4257.772)
   (5544.199,4255.809)
   (5591.921,4253.817)
   (5638.846,4251.799)
   (5684.975,4249.752)
   (5730.308,4247.679)
   (5774.843,4245.578)
   (5818.583,4243.449)
   (5861.526,4241.293)
   (5903.672,4239.109)
   (5945.022,4236.898)
   (5985.575,4234.660)
   (6025.332,4232.394)
   (6064.292,4230.101)
   (6102.456,4227.780)
   (6139.823,4225.431)
   (6176.394,4223.056)
   (6212.168,4220.652)
   (6281.327,4215.763)
   (6347.300,4210.765)
   (6410.087,4205.656)
   (6469.688,4200.438)
   (6526.102,4195.109)
   (6579.331,4189.671)
   (6629.374,4184.123)
   (6676.230,4178.465)
   (6719.901,4172.697)
   (6760.386,4166.819)
   (6831.797,4154.734)
   (6890.464,4142.210)
   (6936.387,4129.246)
   (6990.000,4102.000)
   (7021.641,4035.906)
   (7019.004,3995.477)
   (7004.062,3950.125)
   (6976.816,3899.852)
   (6937.266,3844.656)
   (6885.410,3784.539)
   (6821.250,3719.500)
   (6784.556,3685.135)
   (6744.785,3649.539)
   (6701.938,3612.713)
   (6656.016,3574.656)
   (6607.017,3535.369)
   (6554.941,3494.852)
   (6499.790,3453.104)
   (6441.562,3410.125)
   (6380.259,3365.916)
   (6315.879,3320.477)
   (6248.423,3273.807)
   (6177.891,3225.906)
   (6141.471,3201.495)
   (6104.282,3176.775)
   (6066.324,3151.749)
   (6027.598,3126.414)
   (5988.102,3100.772)
   (5947.837,3074.822)
   (5906.803,3048.565)
   (5865.000,3022.000)
   (5822.922,2995.688)
   (5781.064,2970.188)
   (5739.426,2945.502)
   (5698.008,2921.629)
   (5656.809,2898.569)
   (5615.830,2876.321)
   (5575.071,2854.887)
   (5534.531,2834.266)
   (5494.211,2814.457)
   (5454.111,2795.462)
   (5414.231,2777.280)
   (5374.570,2759.910)
   (5335.129,2743.354)
   (5295.908,2727.610)
   (5256.907,2712.680)
   (5218.125,2698.562)
   (5179.563,2685.258)
   (5141.221,2672.767)
   (5103.098,2661.088)
   (5065.195,2650.223)
   (5027.512,2640.170)
   (4990.049,2630.931)
   (4952.805,2622.504)
   (4915.781,2614.891)
   (4878.977,2608.090)
   (4842.393,2602.103)
   (4806.028,2596.928)
   (4769.883,2592.566)
   (4733.958,2589.018)
   (4698.252,2586.282)
   (4627.500,2583.250)
   (4557.627,2583.470)
   (4488.633,2586.941)
   (4420.518,2593.665)
   (4353.281,2603.641)
   (4286.924,2616.868)
   (4221.445,2633.348)
   (4156.846,2653.079)
   (4093.125,2676.062)
   (4030.283,2702.298)
   (3968.320,2731.785)
   (3907.236,2764.524)
   (3847.031,2800.516)
   (3787.705,2839.759)
   (3729.258,2882.254)
   (3671.689,2928.001)
   (3615.000,2977.000)
   (3558.530,3026.263)
   (3501.621,3072.801)
   (3444.272,3116.614)
   (3386.484,3157.703)
   (3328.257,3196.067)
   (3269.590,3231.707)
   (3210.483,3264.622)
   (3150.938,3294.812)
   (3090.952,3322.278)
   (3030.527,3347.020)
   (2969.663,3369.036)
   (2908.359,3388.328)
   (2846.616,3404.896)
   (2784.434,3418.738)
   (2721.812,3429.856)
   (2658.750,3438.250)
   (2595.249,3443.919)
   (2531.309,3446.863)
   (2466.929,3447.083)
   (2402.109,3444.578)
   (2336.851,3439.349)
   (2271.152,3431.395)
   (2205.015,3420.716)
   (2138.438,3407.312)
   (2071.421,3391.185)
   (2003.965,3372.332)
   (1936.069,3350.755)
   (1867.734,3326.453)
   (1798.960,3299.427)
   (1729.746,3269.676)
   (1660.093,3237.200)
   (1590.000,3202.000)
   (1520.566,3165.789)
   (1452.891,3130.281)
   (1386.973,3095.477)
   (1322.812,3061.375)
   (1260.410,3027.977)
   (1199.766,2995.281)
   (1140.879,2963.289)
   (1083.750,2932.000)
   (1028.379,2901.414)
   (974.766,2871.531)
   (922.910,2842.352)
   (872.812,2813.875)
   (824.473,2786.102)
   (777.891,2759.031)
   (733.066,2732.664)
   (690.000,2707.000)
   (648.691,2682.039)
   (609.141,2657.781)
   (571.348,2634.227)
   (535.312,2611.375)
   (468.516,2567.781)
   (408.750,2527.000)
   (356.016,2489.031)
   (310.312,2453.875)
   (240.000,2392.000)
   (190.781,2336.102)
   (155.625,2280.906)
   (134.531,2226.414)
   (127.500,2172.625)
   (134.531,2119.539)
   (155.625,2067.156)
   (190.781,2015.477)
   (240.000,1964.500)
   (293.613,1939.979)
   (339.536,1928.311)
   (398.203,1917.039)
   (469.614,1906.163)
   (510.099,1900.873)
   (553.770,1895.682)
   (600.626,1890.589)
   (650.669,1885.596)
   (703.898,1880.702)
   (760.312,1875.906)
   (819.913,1871.210)
   (882.700,1866.612)
   (948.673,1862.113)
   (1017.832,1857.713)
   (1053.606,1855.550)
   (1090.177,1853.412)
   (1127.544,1851.298)
   (1165.708,1849.209)
   (1204.668,1847.145)
   (1244.425,1845.106)
   (1284.978,1843.091)
   (1326.328,1841.102)
   (1368.474,1839.136)
   (1411.417,1837.196)
   (1455.157,1835.280)
   (1499.692,1833.389)
   (1545.025,1831.523)
   (1591.154,1829.681)
   (1638.079,1827.864)
   (1685.801,1826.072)
   (1734.319,1824.305)
   (1783.634,1822.562)
   (1833.745,1820.844)
   (1884.653,1819.151)
   (1936.358,1817.482)
   (1988.859,1815.839)
   (2042.156,1814.219)
   (2096.250,1812.625)
   (2151.140,1811.055)
   (2206.827,1809.510)
   (2263.311,1807.990)
   (2320.591,1806.495)
   (2378.667,1805.024)
   (2437.540,1803.578)
   (2497.210,1802.156)
   (2557.676,1800.760)
   (2618.938,1799.388)
   (2680.997,1798.041)
   (2743.853,1796.718)
   (2807.505,1795.420)
   (2871.953,1794.147)
   (2937.198,1792.899)
   (3003.240,1791.675)
   (3070.078,1790.477)
   (3137.713,1789.302)
   (3206.144,1788.153)
   (3275.371,1787.028)
   (3345.396,1785.928)
   (3416.216,1784.853)
   (3451.925,1784.325)
   (3487.833,1783.802)
   (3523.941,1783.286)
   (3560.247,1782.777)
   (3596.752,1782.273)
   (3633.457,1781.775)
   (3670.361,1781.284)
   (3707.464,1780.799)
   (3744.766,1780.320)
   (3782.267,1779.847)
   (3819.967,1779.381)
   (3857.867,1778.920)
   (3895.965,1778.466)
   (3934.263,1778.018)
   (3972.759,1777.576)
   (4011.455,1777.141)
   (4050.350,1776.711)
   (4089.445,1776.288)
   (4128.738,1775.871)
   (4168.230,1775.460)
   (4207.922,1775.055)
   (4247.812,1774.656)
   (4287.902,1774.264)
   (4328.191,1773.878)
   (4368.679,1773.498)
   (4409.366,1773.124)
   (4450.253,1772.756)
   (4491.338,1772.394)
   (4532.623,1772.039)
   (4574.106,1771.690)
   (4615.789,1771.347)
   (4657.671,1771.010)
   (4699.752,1770.680)
   (4742.032,1770.355)
   (4784.512,1770.037)
   (4827.190,1769.725)
   (4870.068,1769.419)
   (4913.145,1769.119)
   (4956.420,1768.826)
   (4999.895,1768.538)
   (5043.569,1768.257)
   (5087.443,1767.982)
   (5131.515,1767.713)
   (5175.786,1767.451)
   (5220.257,1767.194)
   (5264.927,1766.944)
   (5309.796,1766.700)
   (5354.864,1766.462)
   (5400.131,1766.230)
   (5445.597,1766.005)
   (5491.262,1765.785)
   (5537.127,1765.572)
   (5583.190,1765.365)
   (5629.453,1765.164)
   (5675.915,1764.969)
   (5722.576,1764.781)
   (5769.436,1764.599)
   (5816.495,1764.422)
   (5863.754,1764.253)
   (5911.211,1764.089)
   (5958.868,1763.931)
   (6006.724,1763.780)
   (6054.779,1763.635)
   (6103.033,1763.496)
   (6151.486,1763.363)
   (6200.138,1763.236)
   (6248.989,1763.115)
   (6298.040,1763.001)
   (6347.290,1762.893)
   (6396.738,1762.791)
   (6446.386,1762.695)
   (6496.233,1762.606)
   (6546.279,1762.522)
   (6596.525,1762.445)
   (6646.969,1762.374)
   (6697.613,1762.309)
   (6748.455,1762.250)
   (6799.497,1762.198)
   (6850.738,1762.151)
   (6902.178,1762.111)
   (6953.817,1762.077)
   (7005.656,1762.049)
   (7057.693,1762.028)
   (7109.930,1762.012)
   (7162.365,1762.003)
   (7215.000,1762.000)

\put(7665,97){$\PP^1_\lambda$}
\put(7665,2347){$\hb\bad$}
\put(7665,3877){$\hb\good$}
\end{picture}
}

%\end{document}

%%% Local Variables: 
%%% mode: latex
%%% TeX-master: t
%%% End: 

%% file: sec1.tex
%\documentclass{article}
%\usepackage{def}
%\begin{document}

\section{Complete Hurwitz spaces} \label{complete}

The goal of this section is to define arithmetic compactifications of
Hurwitz spaces for $G$-covers. For a given Hurwitz space $H$, such a
compactification $\hb$ should be a proper model of $H$ over $\ZZ$
whose points in positive characteristic correspond to the reductions
of the covers which are parameterized by $H$. To make this precise, one
first has to give the definition of the reduction of a $G$-cover. This
is done in Section \ref{reduction}.

Our definition of the complete Hurwitz space for $G$-covers
essentially follows the approach of \cite{AbrOort98}. We let $\Hb$ be
the closure of the moduli stack of $G$-covers inside a bigger moduli
stack parameterizing certain maps between stably marked curves. Then
we define $\hb$ as the coarse moduli space associated to $\Hb$. This
is the content of Section \ref{hurwitz}. Section \ref{quotientmodel}
discusses a technical problem that arises from our definition. The
reader who is not interested in this abstract approach may wish to
skip these two sections.

\subsection{Reduction of $G$-covers}  \label{reduction}

In this section we give some terminology and recall some general facts
concerning the reduction of $G$-covers to positive characteristic. We
closely follow \cite{Raynaud98}, \S 2.  However, our definition of the
model of a $G$-cover is not exactly the same as Raynaud's. Also, since
we allow the bottom curve to degenerate, we have to consider a
slightly more general situation than in \cite{Raynaud98}, \S 2.

\begin{defn} \label{Gcoverdef}
  Let $K$ be a field and $G$ a finite group. A {\em $G$-cover} defined
  over $K$ is a finite separable morphism $f:Y\to X$ of smooth
  projective and geometrically irreducible $K$-curves together with an
  isomorphism $G\cong\Aut(Y/X)$ such that $|G|={\rm deg}\,f$. 
  We say that a $G$-cover $f$ is {\em tame} if it is tamely ramified.
\end{defn}
 
Throughout this section, we assume the following situation. Let $R$ be
a complete discrete valuation ring with quotient field $K$ of
characteristic zero, and residue field $k=\bar{k}$ of characteristic
$p>0$. Let $f_K\!:Y_K\to X_K$ be a $G$-cover defined over $K$ ($f_K$
is automatically tame). Write $x_{1,K},\ldots,x_{r,K}\in X_K(\Kb)$ for
the branch points and $y_{1,K},\ldots,y_{n,K}\in Y_K(\Kb)$ for the
ramification points of $f_K$. We assume that $2g+r\geq 3$, where $g$
is the genus of $X_K$.

We would like to define a model of the $G$-cover $f_K$ over the ring
$R$. After replacing $K$ by a finite extension $K'/K$ and $R$ by its
integral closure in $K'$, we may assume that the ramification points
$y_{i,K}$ of $f_K$ are $K$-rational. After a further extension of $K$,
we may assume that the smooth stably marked curve $(Y_K;y_{i,K})$
extends to a stably marked curve $(Y_R;y_{i,R})$ over $R$,
\cite{Knudsen83}. In particular, $Y_R$ is semistable over $R$ and the
points $y_{i,K}$ specialize to pairwise distinct, smooth points $y_i$
on the special fiber $Y$ of $Y_R$. Since the stably marked model is
unique, the action of $G$ on $Y_K$ extends to $Y_R$.  Let $X_R:=Y_R/G$
be the quotient scheme and $X$ the special fiber of $X_R$. By
\cite{Raynaudfest}, Appendice, $X_R$ is again a semistable curve over
$R$. Since the ramification points $y_{i,K}$ specialize to pairwise
distinct smooth points on $Y$, the branch points $x_{j,K}$ specialize
to pairwise distinct smooth points $x_j\in X$. According to
\cite{Knudsen83}, the stably marked curve $(X_K;x_{j,K})$ over $K$
extends to a stably marked curve $(X_{0,R};x'_{j,R})$ over $R$, and we
have a well defined contraction morphism $X_R\to X_{0,R}$ sending
$x_{j,R}$ to $x'_{j,R}$. Let $f_R:Y_R\to X_R$ and $f_{0,R}:Y_R\to
X_{0,R}$ be the natural maps and $f:Y\to X$, $f_0:Y\to X_0$ the
induced maps on the special fibers.

\begin{defn} \label{modeldef}
  Let $f_R:Y_R\to X_R$ and $f_{0,R}:Y_R\to X_{0,R}$ be as above.  We
  call $f_R:Y_R\to X_R$ the {\em quotient model} and $f_{0,R}:Y_R\to
  X_{0,R}$ the {\em stable model} over $R$ of the $G$-cover $f_K$.
\end{defn}

The quotient model and the stable model of $f_K$ exist after a finite
extension of $K$. It is clear that these models are stable with
respect to any further extension of $K$. In many places it is better
to work with the quotient model, because it is a finite map. However,
the stable model is easier to study the moduli of. Therefore, our
definition of a complete Hurwitz space will be based on the stable
model.

\begin{defn} \label{reductiondef}\ 
  \begin{enumerate} 
  \item 
    The $G$-cover $f_K$ has {\em (potentially)
    good reduction} if (after a finite extension of $K$) the special
    fiber $f:Y\to X$ of the quotient model of $f_K$ is a tame $G$-cover.
  \item 
    The $G$-cover $f_K$ has {\em (potentially) admissible
    reduction} if (after a finite extension of $K$) the special fiber
    $f:Y\to X$ of the quotient model of $f_K$ is a tame admissible
    cover. By this we mean that $f$ is finite, separable, tamely
    ramified over the smooth locus $X\sm$ and has {\em tame admissible
    ramification} over the ordinary double points of $X$, see
    \cite{HarMum82}, \S 4, or \cite{Saidi97}, where such a cover is
    called {\em kumm\'erien}.
  \item
    The $G$-cover $f_K$ has {\em bad reduction} if it does not have
    potentially admissible reduction.  
  \end{enumerate}
\end{defn}

For the rest of this subsection we will omit the word ``potentially''
and assume that $K$ is chosen such that the quotient and stable model
exist over $R$. Note that $f_K$ has good reduction if and only if it
has admissible reduction and $X_R$ is smooth over $R$. If $f_K$ has
admissible reduction then $X_R=X_{0,R}$, i.e.\ quotient and stable
model are the same. 

Let $W$ be a component of $X$. We call $W$ an {\em original component}
if it is the strict transform of a component of $X_0$ (otherwise, the
map $X\to X_0$ contracts $W$). We will say that $f$ is separable over
$W$ if for one (and therefore for all) components $Z$ of $Y$ above $W$
the restriction $f|_Z:Z\to W$ is a separable morphism. Equivalently,
the inertia group of $Z$ (the group of elements of $G$ acting
trivially on $Z$) is trivial. Proposition \ref{admredprop} below
extends \cite{Raynaud98}, Corollaire 2.4.9, to our situation, which
includes admissible reduction. The proof is essentially the same, with
\cite{Saidi97}, Th\'eor\`eme 3.2, as additional ingredient.

\begin{prop} \label{admredprop}
  The $G$-cover $f_K$ has admissible reduction if and only if $f$ is
  separable over the original components of $X$. 
\end{prop}

From Proposition \ref{admredprop} we can deduce the following well
known fact.

\begin{cor} \label{admredcor}
  Assume that the order of $G$ is prime to the characteristic of $k$.
  Then $f_K$ has potentially admissible reduction. If in addition the
  branch points $x_{j,K}$ of $f_K$ specialize to pairwise distinct
  points on the special fiber of a smooth model of $X_K$, then
  $f_K$ has potentially good reduction.
\end{cor} 

The following example plays a central role in this paper.

\begin{exa} \label{ellexa}
  Let $R$ and $K$ be as before. Choose four $R$-rational points
  $x_1,\ldots,x_4$ on $\PP^1$ such that $x_i\not\equiv x_j\pod{\m}$,
  where $\m\lhd R$ is the maximal ideal. Then $(\PP^1_R;x_i)$ is a
  smooth, stably marked curve over $R$. Let $G$ be the dihedral group
  of order $2p$, where $p$ is an odd prime, equal to the residue
  characteristic of $R$. Let $f_K:E_K\to\PP^1_K$ be a $G$-cover
  branched only in the $4$ points $x_i$, with ramification of order
  $2$. We have $4p$ ramification points $y_{i,j}$ on $E_K$, where
  $1\leq i\leq 4$, $1\leq j\leq p$ and $f_K(y_{i,j})=x_i$. The curve
  $E_K$ has genus $1$. After extending $K$ we may assume that the
  $y_{i,j}$ are $K$-rational.  Choosing e.g.\ $y_{1,1}$ as the origin
  gives $E_K$ the structure of an elliptic curve. Moreover, $f_K$ can
  be written as the composition
  \[ 
       f_K: E_K \lpfeil{p} E_K' \lpfeil{2} \PP^1_K 
  \]
  of a $p$-cyclic isogeny of elliptic curves and a cyclic cover of
  degree $2$. The cover $f_K$ extends to a finite flat morphism
  $f_R:E_R\to\PP^1_R$. Moreover, $E_R$ is an elliptic curve and $f_R$
  factors through an isogeny $\pi_R:E_R\to E'_R$.
  
  There are two cases to consider. First, $\pi_R$ might be \'etale.
  In this case, the ramification points $y_{i,j}$ extend to disjoint
  sections $y_{i,j}:\Spec R\to E_R$ and $f:E_R\to\PP^1_R$ is tamely
  ramified along the sections $x_i:\Spec R\to \PP^1_R$.  In other
  words, $f_K$ has good reduction. Now assume that $\pi_R$ is not
  \'etale. Then its restriction $\pi:E\to E'$ to the special fiber is
  purely inseparable, and for fixed $i$ the $p$ points $y_{i,j}$,
  $j=1,\ldots,p$ specialize to the same point of $E$. We see that
  $(E_R;y_{i,j})$ is not a stably marked curve. Let $(Z_R;y_{i,j})$ be
  the extension of $(E_K;y_{i,j})$ to a stably marked curve over $R$
  and $q_R:Z_R\to E_R$ the contraction morphism. We can identify $E$
  with its strict transform in $Z_R$. The special fiber $Z$ of $Z_R$
  has exactly $5$ components $E,\,Z_1,\ldots,Z_4$.  For
  $i=1,\ldots,4$, the curve $Z_i$ is smooth and of genus $0$,
  connected to $E$ in one point and contains the specialization of the
  points $y_{i,j}$ for $j=1,\ldots,p$. Let $X_R:=Z_R/G$ be the
  quotient.  The special fiber $X:=X_R\otimes_R k$ has $5$ components
  $X_0,\ldots,X_4$. In fact, $X_0$ is the original component, and for
  $i=1,\ldots,4$, the component $X_i$ is the image of $Z_i$ and
  contains the specialization of $x_i$. The restriction of $f:Z\to X$
  to $X_i$ is a $G$-cover $Z_i\to X_i$ ramified in two points, with
  ramification of order $2$ and $2p$ (so it is not tame).
\end{exa}

\subsection{Complete Hurwitz stacks} \label{hurwitz}

In this section we define the concept of a {\em complete Hurwitz
stack}, following the idea of \cite{AbrOort98}. Since there are many
different versions of Hurwitz stacks, we do this first in detail for
one specific kind, namely for $\H_{[r]}\inn(G)$, the inner Hurwitz
stack for $G$-covers of genus $0$ curves with unordered branch
points. Then we define several variants of the above. In Section
\ref{coarse} we look at Hurwitz spaces as coarse moduli spaces.

In this paper we only consider Hurwitz spaces for Galois covers.
Moreover, the target curve of the cover will always be of genus $0$
and is considered ``up to isomorphism''. Thus, we only consider
``reduced'' Hurwitz spaces in the terminology used by Fried
\cite{Fried87}. The genus zero assumption is made only to simplify the
notation. It is easy to extend all our definitions to nonreduced
Hurwitz spaces. It seems much less trivial to do the same for Hurwitz
spaces parameterizing non-Galois covers.  For instance, in
\cite{AbrOort98} a completion of the classical Hurwitz stack for
simple covers is constructed, using the moduli space of stable maps as
an ambient space. This construction is more involved than the one we
give here. Another problem, discussed in \cite{AbrOort98}, Section 4.1,
is that ``taking quotients by finite groups does not commute with base
change''.  Therefore, the method proposed in \cite{AbrVis98} of
studying complete moduli of non-Galois covers by going to the Galois
closure probably does not work very well with our definition of
complete Hurwitz stacks. A related problem is discussed in Section
\ref{quotientmodel}.

%We denote by $\Mb_{g,[n]}$ the stack of stably $n$-marked
%curves. Obviously, $\Mb_{g,[n]}$ can be identified with the quotient
%stack $[\Mb_{g,n}/{\mathfrak{S}_n}]$, where $\Mb_{g,n}$ is the stack
%of $n$-pointed stable curves and $\mathfrak{S}_n$ is the symmetric
%group. Therefore, $\Mb_{g,[n]}$ is an algebraic stack, smooth and
%proper over $\ZZ$. 

\subsubsection{The ambient stack} \label{ambient}

Let $(X/S,C)$ and $(Y/S,D)$ be stably marked curves, defined over the
same scheme $S$.  A morphism of stably marked curves is an
$S$-morphism $f:Y\to X$ such that $f(D)\subset C$. To ease notation,
we will usually write $X$ and $Y$ instead of $(X/S,C)$ and
$(Y/S,D)$. We fix an integer $r>0$ and let $\SB\rr$ be the following
category.  Objects of $\SB\rr$ are morphisms $f:Y\to X$ between stably
marked curves such that $X$ has genus $0$ and is stably $r$-marked. A
morphism between an $S$-object $f:(Y,D)\to(X,C)$ and an $S'$-object
$f':(Y',D')\to(X',C')$ of $\SB\rr$ consists of a Cartesian diagram
\begin{equation} \label{Smorph}
  \begin{CD}
     Y'  @>>> Y \\
     @V{f'}VV  @V{f}VV \\
     X'  @>>> X \\
     @VVV      @VVV    \\
     S'   @>>>  S      \\
  \end{CD}
\end{equation}
such that $D'=D\times_S S'$ and $C'=C\times_S S'$. It is clear that
$\SB\rr$ is a stack. 

Let $G$ be a finite group and $\SB\rr\inn(G)$ the following category.
Objects of $\SB\rr\inn(G)$ (over a scheme $S$) are pairs $(f,\sigma)$,
where $f:Y\to X$ is an object of $\SB\rr$ defined over $S$ and
$\sigma:G\to\Aut(Y/X)$ is an action of $G$ on $Y$ commuting with $f$
such that the induced action on every geometric fiber of $f$ is
faithful (equivalently, $\sigma$ induces a closed immersion
$\tilde{\sigma}:G_S\inj\AUT(Y/X)$ of group schemes). Mostly we will
omit the map $\sigma$ and simply write $f:Y\to X$ for an object of
$\SB\rr\inn(G)$. A morphism between two objects $f':Y'\to X'$ and
$f:Y\to X$ of $\SB\rr\inn(G)$ is an $\SB\rr$-morphism \zgl{Smorph} such
that the top arrow $Y'\to Y$ is $G$-equivariant.  Again it is clear
that $\SB\rr\inn(G)$ is a stack.

\begin{prop} \label{stackprop}
  The stacks $\SB\rr$ and $\SB\rr\inn(G)$ are algebraic, separated and
  locally of finite type over $\ZZ$. 
\end{prop}

\proof The stack $\Mb_{g,[n]}$ classifying stably $n$-marked curves of
fixed genus $g$ is algebraic, separated and of finite type over $\ZZ$.
A standard Hilbert scheme argument (see e.g.\ \cite{MumFog},
Chap. 0.5) shows that $\SB\rr$ is algebraic, separated and locally of
finite type over $\ZZ$. Let $f:Y\to X$ be an object of $\SB\rr$
defined over a scheme $S$. Since $\AUT(Y/X)$ is a finite $S$-group
scheme, the functor $\Hom_S(G_S,\AUT(Y/X))$ is represented by a finite
$S$-scheme. It is clear that the natural morphism
$S\times_{\SB\rr}\SB\rr\inn(G)\inj\Hom_S(G_S,\Aut(Y/X))$ is a locally
closed immersion. Therefore, the forgetful morphism
$\SB\rr\inn(G)\to\SB\rr$ is relatively representable, separated and of
finite type.  This completes the proof of the proposition.  \Endproof

\subsubsection{The complete inner Hurwitz stack} \label{inner}

Let $r$ and $G$ be as before. We define the Hurwitz stack $\Hin\rr(G)$
as follows. Objects of $\Hin\rr(G)$ over a scheme $S$ are morphisms
$f:Y\to X$ between smooth $S$-curves, together with an action of $G$
on $Y$, commuting with $f$, such that the following holds. The curve
$X/S$ has genus $0$ and the geometric fibers of $f$ are tame
$G$-covers (see Definition \ref{Gcoverdef}) with exactly $r$ branch
points. Morphisms in $\Hin\rr(G)$ are Cartesian diagrams of the form
\zgl{Smorph} such that the top horizontal arrow is $G$-equivariant. We
call an object $f:Y\to X$ of $\Hin\rr(G)$ a {\em tame $G$-cover},
defined over $S$. It is proved e.g.\ in \cite{diss} that $\Hin\rr(G)$
is an algebraic stack, smooth and of finite type over $\ZZ$.

Let $f:Y\to X$ be an $S$-object of $\Hin\rr(G)$. Then $f$ is finite
and tamely ramified along a divisor $C\subset X$ which is finite
\'etale of degree $r$ over $S$. Hence $(X/S,C)$ is a (smooth) stably
$r$-marked curve. Moreover, $(Y/S,D)$ is a (smooth) stably marked
curve, where $D:=f^{-1}(C)\subset Y$ is the (reduced) inverse image of
$C$. Therefore, we obtain a natural monomorphism
\begin{equation} \label{hurwinjeq}
    \Hin\rr(G) \;\inj\; \SB\rr\inn(G),
\end{equation}
identifying $\Hin\rr(G)$ with a full subcategory of $\SB\rr\inn(G)$.
We will show in Proposition \ref{properprop} below that
\zgl{hurwinjeq} is a locally closed immersion. Note however that we do
not need this fact to make the following definition.

\begin{defn} \label{hurwdef1}
  The {\rm complete Hurwitz stack} $\Hinb\rr(G)$ is the closure of
  $\Hin\rr(G)$ inside $\SB\rr\inn(G)$ (i.e.\ the smallest closed
  substack of $\SB\rr\inn(G)$ containing $\Hin\rr(G)$ as a full
  subcategory).
\end{defn}

Proposition \ref{stackprop} shows that $\Hinb\rr(G)$ is an algebraic
stack, separated and locally of finite type over $\ZZ$.  Let $k$ be an
algebraically closed field and $f:Y\to X$ an object of $\Hinb\rr(G)$
defined over $k$. Choose an \'etale neighborhood $U\to\Hinb\rr(G)$ of
the point $s:\Spec k\to\Hinb\rr(G)$ corresponding to $f$ and let
$\eta:\Spec K\to U$ be a generic point of some irreducible component
of $U$. By \cite{Hartshorne}, Exercise II.4.11, $\eta$ extends to a
morphism $\eta_R:\Spec R\to U$, where $R$ is a discrete valuation ring
of $K$ with residue field $k$, and the restriction of $\eta_R$ to the
special point is equal to $s$. The morphism $\Spec R\to\Hinb\rr(G)$
corresponds to an $R$-object $f_R:Y_R\to X_R$ of $\Hinb\rr(G)$ with
special fiber $f:Y\to X$. Since $\Hin\rr(G)$ is dense in
$\Hinb\rr(G)$, the generic fiber $f_K:Y_K\to X_K$ of $f_R$ is actually
an object of $\Hin\rr(G)$, i.e.\ a tame $G$-cover.  Moreover, $K$ has
characteristic $0$. We are essentially (modulo taking the completion
of $R$) in the situation of Section \ref{reduction}. It is clear that
$f_R:Y_R\to X_R$ is the stable model of the $G$-cover $f_K$. In
particular, $f_K$ has good reduction if and only if $f$ is an object
of $\Hin\rr(G)$. We say that $f$ is a {\em bad cover} if $f_K$ has bad
reduction (Definition \ref{reductiondef}).  By Proposition
\ref{admredprop}, $f$ is a bad cover if and only if some irreducible
component of $Y$ has a nontrivial inertia group (with respect to the
action of $G$).

\begin{lem} \label{badlem}
  There is a unique closed reduced substack
  $\Hinb\rr(G)\bad\subset\Hinb\rr(G)$ characterized by the following
  property. For an algebraically closed field $k$, a $k$-object
  $f:Y\to X$ of $\Hinb\rr(G)$ is a bad cover if and only if it is an
  object of the substack $\Hinb\rr(G)\bad$.
\end{lem}

\proof By \cite{Hartshorne}, Lemma II.4.5, it suffices to show that
the subset of bad covers is stable under specialization. More
precisely, let $R$ be a discrete valuation ring and $f_R:Y_R\to X_R$
an object of $\Hinb\rr(G)$ defined over $R$. Assume that the generic
fiber $f_K:Y_K\to X_K$ of $f_R$ is a bad cover. We have to show that
the special fiber $f:Y\to X$ is a bad cover, too. As remarked above,
$f_K$ (resp.\ $f$) is a bad cover if and only if some irreducible
component of $Y_K$ (resp.\ $Y$) has a nontrivial inertia group.
Clearly, this property is stable under specialization.  \Endproof

By the lemma, $\Hinb\rr(G)\adm:=\Hinb\rr(G)-\Hinb\rr(G)\bad$ is a
dense open substack of $\Hinb\rr(G)$. The discussion before Lemma
\ref{badlem} shows that the geometric points of $\Hinb\rr(G)\adm$
correspond to tame admissible covers $f:Y\to X$ over $k$ which arise
as the reduction of tame $G$-covers. It follows that an object $f:Y\to
X$ of $\Hinb\rr(G)\adm$ (defined over an arbitrary scheme $S$) lies in
the full subcategory $\Hin\rr(G)$ if and only if $X/S$ is smooth.
Since smoothness of $X/S$ is an open condition on $S$, $\Hin\rr(G)$ is
an open substack of $\Hinb\rr(G)\adm$, and hence an open substack of
$\Hinb\rr(G)$.

\begin{prop} \label{properprop} 
  The algebraic stack $\Hinb\rr(G)$ is reduced, proper and of finite
  type over $\ZZ$ and contains $\Hin\rr(G)$ and $\Hinb\rr(G)\adm$ as
  dense open substacks.
\end{prop}

\proof It only remains to show that $\Hinb\rr(G)$ is reduced and
proper. We know that the dense open substack $\Hin\rr(G)$ is reduced,
therefore $\Hinb\rr(G)$ is reduced as well. To prove properness, we
apply \cite{EGA2}, Corollaire 7.3.10. We are immediately reduced to
the following situation. Let $R$ be a complete discrete valuation ring
with residue field $k=\bar{k}$ of characteristic $p$ and quotient
field $K$ of characteristic $0$. Let $f_K:Y_K\to X_K$ be an object of
$\Hin\rr(G)$, i.e.\ a $G$-cover. This is the situation of Section
\ref{reduction}.  After a finite extension of $K$, $f_K$ has a stable
model $f_R:Y_R\to X_R$ over $R$ (Definition \ref{modeldef}).
Obviously, $f_R$ is an object of $\SB\rr\inn(G)$.  Therefore, $f_R$ is
an object of $\Hinb\rr(G)$, by Definition \ref{hurwdef1}. This proves
that $\Hinb\rr(G)$ is proper.  \Endproof

The stack $\Hin\rr(G)$ is called the {\em inner Hurwitz stack} for
$G$-covers with $r$ (unordered) branch points. We will say that the
stack $\SB\rr\inn(G)$ is the {\em ambient stack} of $\Hin\rr(G)$. We
will call the stack $\Hinb\rr(G)$ the {\em completion} of
$\Hin\rr(G)$.

\begin{rem} \label{admrem}
  It is easy to check that the stack $\Hinb\rr(G)\adm$ can be
  identified with the compactification of $\Hin\rr(G)$ constructed
  in \cite{diss} or \cite{AbrVis98}. It follows from loc.cit.\ that
  $\Hinb\rr(G)\adm$ is smooth over $\ZZ$ and proper over $\ZZ[1/|G|]$.
\end{rem}

\begin{var} \label{variant1}
  Let $\H\rr\ab(G)$ be the stack whose objects are $S$-morphisms
  $f:Y\to X$ which are locally on $S$ tame $G$-covers. More precisely,
  after an \'etale localization $S'\to S$, there exists an action of
  $G$ on $Y$ such that $f:Y\to X$ becomes an object of
  $\H\rr\inn(G)$. We call $\H\rr\ab(G)$ the {\em absolute Hurwitz
  stack} for $G$-covers with $r$ (unordered) branch points, see also
  \cite{FriedVoe91}. We embed $\H\rr\ab(G)$ into an ambient stack
  $\SB\ab\rr(G)$. Objects of $\SB\ab\rr(G)$ are pairs $(f,\G)$, where
  $f:Y\to X$ is an object of $\SB\rr$ defined over a scheme $S$ and
  $\G\subset\AUT(Y/X)$ is an \'etale subgroup scheme which becomes
  isomorphic to the constant $S$-group scheme $G_S$ after an \'etale
  localization of $S$. We define the completion $\Hb\rr\ab(G)$ as the
  closure of $\H\rr\ab(G)$ inside $\SB\rr\ab(G)$.
\end{var}

\begin{var} \label{variant2}
  Let $\H_r\inn(G)$ and $\H_r\ab(G)$ be the inner resp.\ absolute
  Hurwitz stack for $G$-covers with $r$ {\em ordered} branch
  points. We can embed $\H_r\inn(G)$ into an ambient stack
  $\SB\inn_r(G)$. Objects of $\SB\inn_r(G)$ are objects $f:Y\to X$ of
  $\SB\rr\inn(G)$ together with $r$ sections $x_1,\ldots,x_r:S\to X$
  such that $C=\cup_i x_i(S)$ is the mark of the stably marked curve
  $X$. We define the completion $\Hb_r\inn(G)$ as the closure of
  $\H_r\inn(G)$ inside $\SB_r\inn(G)$. Similar for $\Hb_r\ab(G)$.
\end{var}
We obtain a diagram
\begin{equation} \label{basicdiag}
\begin{CD}
      \Hinb_r(G) @>>> \Hinb\rr(G)  \\
      @VVV                   @VVV             \\
      \Habb_r(G) @>>> \Habb\rr(G).  \\
\end{CD}
\end{equation}
The vertical arrows in \zgl{basicdiag} are principal
$\Aut(G)$-bundles, the horizontal arrows are principal
$\mathfrak{S}_r$-bundles. All stacks in \zgl{basicdiag} are algebraic,
reduced and proper and of finite type over $\ZZ$.

\subsubsection{Complete Hurwitz stacks for a given type} 
\label{morevariants}

With $G$ and $r$ as before, let $\Cl=(C_1,\ldots,C_r)$ be an $r$-tuple
of conjugacy classes of $G$. We denote by
$\QQ(\Cl)\subset\QQ(\zeta_n)$ the field generated by the values
$\chi(\tau)$, where $\chi$ runs over all irreducible characters of $G$
and $\tau\in C_i$, for $i=1,\ldots,r$.  In other words, $\QQ(\Cl)$ is
the minimal number field over which every class $C_i$ becomes
rational. Let us make a ``choice of $n$th root of unity over
$\QQ(\Cl)$'', i.e.\ we choose an orbit under
$\Gal(\bar{\QQ}/\QQ(\Cl))$ of primitive $n$th roots of unity, see
\cite{SerreTopics}, Section 8.2.1. Let $\Lambda\subset \QQ(\Cl)$ be a
Dedekind domain with fraction field $\QQ(\Cl)$. We define an open
substack
\begin{equation}
   \H_r\inn(\Cl)_\Lambda \subset \H_r\inn(G)\otimes_{\ZZ}\Lambda,
\end{equation}
corresponding to $G$-covers with {\em inertia type} $\Cl$. To be more
precise, let $f:Y\to(X;x_i)$ be an object of
$\Hin_r(G)\otimes\Lambda$, defined over an algebraically closed field
$k$. We say that $f$ has inertia type $\Cl$ if $C_i$ is the conjugacy
class associated to the branch point $x_i$ (with respect to a
canonical choice of $m_i$th root of unity, induced by the natural map
$\Lambda\to k$), compare \cite{Voelklein}, Section 2.2.1. An object of
$\Hin_r(G)\otimes\Lambda$, defined over an arbitrary scheme $S$, is
said to have inertia type $\Cl$ if all its geometric fibers have
inertia type $\Cl$. We define the completion $\Hb_r\inn(\Cl)_\Lambda$
as the closure of $\H_r\inn(\Cl)_\Lambda$ inside
$\SB_r\inn(G)\otimes_{\ZZ}\Lambda$. Clearly, $\Hb_r\inn(\C)_\Lambda$
is an algebraic stack, proper and of finite type over $\Lambda$. If
the choice of $\Lambda$ is understood, we will omit it from the
notation, and we say that $\Hb_r\inn(\Cl)$ is defined over $\Lambda$,
or that $\Lambda$ is the {\em domain of definition} of
$\Hb_r\inn(\Cl)$.

In a similar way, we can define further variants of complete Hurwitz
stacks: $\Hb_r\ab(\Cl)$, $\Hb_{[r]}\inn(\Cl)$ and $\Hb_{[r]}\ab(\Cl)$.
The domains of definition of these stacks are Dedekind domains whose
fraction fields are suitable subfields of $\QQ(\Cl)$. For instance,
$\Hb_{[r]}\inn(\Cl)$ can be defined over the ring of integers of the
smallest field over which $\Cl$ as a tuple is rational, compare
\cite{Voelklein}, Definition 3.15.

\subsubsection{Complete Hurwitz spaces as coarse moduli spaces}
\label{coarse}

Let $\Hb:=\Hb_r\inn(\Cl)$ be the complete Hurwitz stack over
$\Lambda$, as defined in Section \ref{morevariants}. We denote by
$\hb:=\hb_r\inn(\Cl)$ the associated coarse moduli space, and call it
the {\em complete Hurwitz space} over $\Lambda$ for $G$-covers with
inertia type $\Cl$. By construction \cite{KeMo97}, $\hb$ is an
algebraic space, proper and of finite type over $\Lambda$, and
contains the usual Hurwitz space $H=H_r\inn(\Cl)$ as a dense open
subspace.  Actually, $H$ is a scheme and is smooth over $\Lambda$, see
\cite{diss}. We define a closed subspace $\hb\bad$ as the image of the
natural morphism $\Hb\bad\to\hb$. Thus, $\hb\adm:=\hb-\hb\bad$ is the
coarse moduli space associated to $\Hb\adm=\Hb-\Hb\bad$ and contains
$H$ as a dense open subset. According to \cite{diss}, $\hb\adm$ is a
normal scheme.

Let $\Hb\to\Mb_{0,r}$ be the natural forgetful morphism. It is known
that $\Mb_{0,r}$ is represented by a smooth projective scheme over
$\ZZ$, see \cite{GHP}. By the universal property of the coarse moduli
space, we obtain a morphism $\hb\to\Mb_{0,r}$.

\begin{prop} \label{coarseprop}
  Assume that the complete Hurwitz stack $\Hb$ is normal and that the
  forgetful morphism $\Hb\to\Mb_{0,r}\otimes\Lambda$ is relatively
  representable and finite. Then the complete Hurwitz space $\hb$ is a
  normal scheme, finite over $\Mb_{0,r}\otimes\Lambda$.
\end{prop}

\proof If $\Hb$ is normal and $\Hb\to\Mb_{0,r}\otimes\Lambda$ finite,
then the coarse moduli space $\hb$ is the normalization of
$\Mb_{0,r}\otimes\Lambda$ in $H$, see \cite{DelRap}, Proposition
IV.3.10.  This proves the proposition.  \Endproof

%\begin{rem} \label{coarserem}
%  Let $\tilde{H}\bad$ be the coarse moduli space associated to
%  $\Hb\bad$. The natural map 
%  \begin{equation} \label{coarseeq1}
%         \tilde{H}\bad \To \hb\bad 
%  \end{equation} 
%  induces a bijection on geometric points. However, it might not be an
%  isomorphism, in general.
%\end{rem}

\begin{var} \label{variant3}
  In the same manner, we define complete Hurwitz spaces
  $\hb\rr\inn(G)$, $\hb_r\ab(G)$, etc. Proposition \ref{coarseprop}
  applies as well. There are natural maps between all these variants.
  For instance, diagram \zgl{basicdiag} induces an analogous diagram
  of finite morphisms between algebraic spaces. But unlike the maps in
  \zgl{basicdiag}, these maps are in general not \'etale.
\end{var}

\subsection{Quotient model versus stable model} \label{quotientmodel}

Let us fix a finite group $G$ and an integer $r\geq 0$. Let
$\Hb:=\Hinb_r(G)$ be the complete inner Hurwitz stack defined in
Section \ref{hurwitz}. In Section \ref{reduction} we have defined two
different models of a $G$-cover $f_K:Y_K\to X_K$, where $K$ is a
complete discrete valued field: the quotient model $f_R:Y_R\to X_R$
and the stable model $f_{0,R}:Y_R\to X_{0,R}$. The definition of $\Hb$
was made such that the stable model $f_{0,R}:Y_R\to X_{0,R}$ of the
$G$-cover $f_K$ is an object of $\Hb$.  One drawback of the stable
model is that it is in general not a finite map. Over the discrete
valuation ring $R$, one can recover the quotient model from the stable
model by taking the quotient scheme $X_R:=Y_R/G$. The problem is that
taking quotients does not commute with arbitrary base change. So it is
not clear how to define a quotient model of an object $f_0:Y\to X_0$
of $\Hb$ over an arbitrary scheme $S$. In this section we propose a
definition that works well when $S$ is either the spectrum of an
algebraically closed field or a normal scheme $S$ with function field
of characteristic $0$. This will be enough for our purposes.

\begin{defn}
  Let $S$ be a scheme and $f_0:Y\to X_0$ an object of $\Hb$ defined
  over $S$. Let $f:Y\to X$ be a finite morphism between marked
  semistable curves commuting with the action of $G$ on $Y$. We say
  that $f$ is a {\em quotient model} of $f_0$ if $f_0$ is the
  composition of $f$ with an $S$-morphism $X\to X_0$ and if for every
  geometric fiber $f_s:Y_s\to X_s$ of $f$ the following holds.
  \begin{enumerate} \item The natural morphism $Y_s/G\to X_s$ induces
  a bijection on geometric points.  \item Let $Y_i$ be a component of
  $Y_s$ and $X_i$ the component of $X_s$ under $Y_i$. Then the degree
  of $Y_i$ over $X_i$ is equal to the order of the stabilizer
  $D(Y_i)\subset G$ of $Y_i$.  \end{enumerate}
\end{defn}

\begin{prop} \label{quotprop}
  Let $S$ be either the spectrum of an algebraically closed field $k$
  or a normal scheme, generically of characteristic $0$. Let $f_0:Y\to
  X_0$ be an object of $\Hb$ over $S$. Then there exists a unique
  quotient model $f:Y\to X$ of $f_0$. 
\end{prop} 

\proof Let us first assume that $S$ is a normal scheme, generically of
characteristic $0$. We may assume that $S=\Spec R$ is local, with
algebraically closed residue field $k$. The generic fiber $f_K:Y_K\to
X_K$ is an admissible cover. In particular, $f_K$ is a quotient model
of itself.  Let $X:=Y/G$. It follows from \cite{deJong97}, Proposition
4.2, and \cite{KatzMazur}, Corollary A.7.2.2, that $f:Y\to X$ is a
quotient model of $f_0$.  To show that it is unique, suppose we have
another quotient model $f':Y\to X'$ of $f_0$. Since $f'$ commutes with
the action of $G$ on $Y$, there exists a morphism $\kappa:X\to X'$ of
$S$-schemes such that $\kappa\circ f=f'$.  We claim that $\kappa$ is
an isomorphism of $X_0$-schemes.  Since $X$ and $X'$ are flat over $S$
and $\kappa$ is the identity on the generic fiber, it suffices to
prove that the restriction of $\kappa$ to the special fiber is an
isomorphism of $X_0$-schemes.  Hence we have reduced the proposition
to the case $S=\Spec k$. Let us now prove this case. We have seen
before that every $k$-object $f_0:Y\to X_0$ of $\Hb$ is the reduction
of a $G$-cover $f_K:Y_K\to X_K$, where $K$ is the quotient field of a
discrete valuation ring $R$ with residue field $k$.  Therefore, the
special fiber $f:Y\to X$ of the quotient model of $f_K$ is a quotient
model of $f_0$. It remains to prove its uniqueness. The quotient
$X':=Y/G$ is a semistable curve over $k$ and the natural map $X'\to X$
is a bijection on geometric points, by assumption. Let $Y_i$ be a
component of $Y$ and $X_i$ resp.\ $X_i'$ the component of $X$ resp.\ 
of $X'$ under $Y_i$. It follows from \cite{Hartshorne}, Proposition
IV.2.5, that $X_i$ is the $n$th Frobenius twist $(X_i')^{F^n}$ of
$X_i'$, where $p^n$ is the order of the inertia group $I(Y_i)\subset
G$ of $Y_i$. Clearly, this characterizes the map $X'\to X$ and hence
the quotient model $f:Y\to X$ uniquely.  \Endproof

%\bibliographystyle{abbrv} \bibliography{../hurwitz}

%\end{document}

%%% Local Variables: 
%%% mode: latex
%%% TeX-master: t
%%% End: 

%% file: sec3.tex
%\documentclass{article}
%\usepackage{def}

%\begin{document}

\section{Semistable reduction} \label{semistable}

In this section the notations and conventions are as in Section
\ref{reduction}. Let us recall them briefly. Let $R$ be a complete
discrete valuation ring with quotient field $K$ of characteristic zero
and residue field $k=\bar{k}$ of characteristic $p$. Let $f_K\!:Y_K\to
X_K$ be a $G$-cover of smooth curves defined over $K$. Let $g=g(X_K)$
and $r$ the number of branch points of $f_K$. Let $f_R\!:Y_R\to X_R$
be the quotient model of $f_K$ defined in Definition \ref{modeldef}.
Write $f\!:Y\to X$ for its special fiber. The branch points of $f_K$
specialize to distinct points $x_1,\ldots,x_r$ of the smooth locus of
$X$. Let $m_i$ be the ramification index of $x_i$ in $f_K$. Recall
that we also defined a different model for $f_K$, called the stable
model $f_{0,R}\!:Y_R\to X_{0,R}$.  It is obtained from $f_R$ by
composing with a contraction map $X_R\to X_{0,R}$.

 We will determine the structure of the reduction $f\!:Y\to X$
in case the cover $f_K$ has bad reduction, under suitable assumptions
(Condition \ref{cond1} below). All results of this section rely on the
results of \cite{Raynaud94} and \cite{Raynaud98}. Results from these
papers are recalled very briefly and the reader is referred to these
papers for more details. 

Condition \ref{cond1} plays an essential role in the rest of the
paper.  It will enable us to compute the structure of $f\!:Y\to X$.
Most importantly, we will assume that the normalizer of a $p$-Sylow
group $P$ of $G$ is a dihedral group. We will define an {\sl auxiliary
  cover} $g\!:Z\to X$. The condition on the normalizer $N_G(P)$ will
imply that this auxiliary cover is equivariant under a (dihedral)
subgroup of $N_G(P)$. In Section \ref{reductionthm} we will relate the
deformation theory of $f$ to the deformation theory of $g$.

In Section \ref{inertia} we will recall some results on inertia and
decomposition groups of points and components of $Y$. In Section
\ref{uconn} a technical lemma is proved. This enables us to extend
some constructions and notations from \cite{Raynaud98} to the case
that $X_0$ is not smooth in Section \ref{auxsec}. In Section
\ref{modularred} we show that $f\!:Y\to X$ has a fairly simple
structure, which we call {\sl modular reduction}. In Section
\ref{boundarypts} we study the reduction of the degenerate covers in
characteristic zero, i.e.\ the covers corresponding to cusps of the
Hurwitz space.

\subsection{Inertia and decomposition groups}\label{inertia}

\begin{lem} 
  A singular point of $Y$ maps to a singular point of $X$, i.e.\ $G$
  acts on $Y$ {\sl without inversions}.
\end{lem}

\proof 
The ramification points $y_{i,K}$ specialize to distinct smooth points
of $Y$, by construction. It follows from \cite{Raynaud98}, Proposition
2.3.2.b, that there are no inversions.
\Endproof

Note in particular that the above lemma implies that the irreducible
components of $Y$ do not intersect themselves. For an irreducible
component $Z$ of $Y$, we denote by $I(Z)$ its inertia group, and
$D(Z)$ its decomposition group. Analogously, for a closed point $y$ of
$Y$, we will write $I(y)$ for its inertia group $I(y)$. It is equal to
the decomposition group of $y$, since $k$ is algebraically closed.

\begin{lem}\label{inertia1lem} 
\begin{itemize}
\item[(a)] Let $Z$ be an irreducible component of $Y$. Then $I(Z)$ is a
  $p$-group and a normal subgroup of the decomposition group $D(Z)$.
\item[(b)] Let $y_R\!:\Spec(R)\to Y_R$ be a section whose image is
  contained in the smooth locus of $Y_R$. Let $Z$ be the irreducible
  component of $Y$ on which ${y}:=y_R\otimes k$ lies. Write
  $m=p^\alpha n$, with $\gcd(n,p)=1$, for the ramification index of
  $y_K$ in $f$. Then $I(Z)$ is normal in $I({{y}})$ and
  $I({{y}})/I(Z)$ is cyclic of order $n$.
\item[(c)] Any branch point $x_{i,K}$ of $f_K$ whose ramification
  index $m_i$ is prime-to-$p$, specializes to a component of $X$ over
  which $f$ is separable.
\end{itemize}
\end{lem}

\proof Part (a) and (b) follow from \cite{Raynaud94}, Lemme 6.3.3. Let
$x_{i,R}$ be a branch point such that the ramification index $m_i$ of
$x_{i,K}$ is prime-to-$p$. Suppose it specializes to a component $W$
of $X$ over which $f$ is inseparable. Then the inertia group of any
component $Z$ of $Y$ mapping to $W$ is nontrivial. Hence
$x_i:=x_{i,R}\otimes k$ will be branched of order $p^am_i$ with $a>0$,
by Part (b). This is in contradiction with the assumption that the
ramification points specialize to distinct points on $Y$. This proves
(c). \Endproof

\begin{lem}\label{inertia2lem}
  Let $y$ be a singular point of $Y$. Let $Z_1,Z_2$ be the two
  irreducible components of $Y$ passing through ${y}$.
\begin{itemize}
\item[(a)] The inertia group $I(y)$ is an extension of a cyclic group
  of order prime-to-$p$ by a $p$-group.
\item[(b)] The groups $I({Z_1})$ and $I({Z_2})$ are normal subgroups of
  $I(y)$ and $\mathopen<I({Z_1}),I({Z_2})\mathclose>$ is the $p$-Sylow group of $I(y)$.
\end{itemize}
\end{lem}

\proof Part (a) follows from \cite{Raynaud98}, Proposition 2.3.2.a,
since $G$ acts without inversions. Part (b) is proved analogous to
\cite{Raynaud94}, Lemme 6.3.6.iii. The assumptions in that lemma
differ from the assumptions in the present case, but the proof goes
through.  \Endproof

\subsection{First properties}\label{uconn}

We will suppose now that $f_K\!: Y_K\to X_K$ is a $G$-cover branched
at four points $x_1,\ldots, x_4$, where $X_K$ has genus zero. Let
$m_i$ be the ramification index of a point above $x_i$. We suppose
that $f_K$ has bad reduction, and denote by $f\!:Y\to X$ the special
fiber of the quotient model of $f_K$.

\begin{notation}
  Let $U$ be the union of the components $W$ of $X$ such that $f$ is
  inseparable over $W$. Let $P$ be a $p$-Sylow group of $G$. We denote
  by $N_G(P)$ the normalizer of $P$ in $G$ and by $C_G(P)$ the
  centralizer of $P$ in $G $.
\end{notation}

\begin{cond}\label{cond1} 
  In the rest of the paper we will assume the following conditions to
  hold:
\begin{tabbing}
(c) \quad $m_1,\ldots,m_4$ 
prime-to-$p$\qquad,\= \kill
 (a)\quad  $p\neq 2$, \> (c) \quad  $p||\, |G|$,\\
(b)\quad  $m_1,\ldots,m_4$ 
prime-to-$p$,\>
(d)\quad $N_G(P)$ is a dihedral group.
\end{tabbing}
\end{cond}

\begin{exa} 
  Here are some examples of groups for which Condition \ref{cond1}.(d)
  is satisfied.  Let $G=PSL_2(\ell^\alpha)$, where $\ell>2$ is a
  prime.  Suppose that $p\neq \ell$ is a prime exactly dividing
  $|G|=(\ell^\alpha-~1)\ell^\alpha(\ell^\alpha+1)/2.$ Then $N_G(P)$ is
  a dihedral group of order $\ell^\alpha-1$ or order $\ell^\alpha+1$,
  \cite{Huppert}, Abschnitt II.8. Special cases are $G=PSL_2(5)=A_5$
  and $G=PSL_2(9)=A_6$.
\end{exa}

Recall that there is a map $X\to X_0$ which contracts some components;
the strict transforms of the components of $X_0$ in $X$ are called the
{\sl original} components.  Since we assumed that $r=4$, there are two
possibilities for $X_0$: it can be smooth or not. In case $X_0$ is
singular, it will consist of two genus zero components which meet in a
unique point. We will call these components $W_1$ and $W_2$ and will
also write $W_1,W_2$ for their strict transform in $X$. Similarly, we
will write $X_0$ for the strict transform of $X_0$ in $X$, in case
$X_0$ is smooth.

Suppose $X_0$ is not smooth. Since $X_0$ is stably marked, there will
be exactly two branch points specializing to each of the components
$W_i$. In $X$ the two components $W_1$ and $W_2$ are connected by a
chain of $\PP^1$'s, since the dual graph of $X$ is a tree. Let
$\Lambda$ the union of $W_1$, $W_2$ and the components connecting the
two. We will say that a component $W$ of $X$ is a {\sl tail} if it is
not contained in $\Lambda$ and meets the rest of $X$ in a unique
point.

In case $X_0$ is smooth we just put $\Lambda=\{X_0\}$. The definition
of {\sl tail} then becomes the usual definition of tail, i.e.\ one
views the dual graph of $X$ to be oriented from $X_0$.

\begin{lem}\label{uconnlem}
  Suppose that $f_K$ has bad reduction.  Let $W$ be an irreducible
  component of $X$. Then $f|_W$ is separable if and only if $W$ is a
  tail.
\end{lem} 

In case $X_0$ is smooth, the statement of this lemma is the same as
\cite{Raynaud98}, Lemme 3.1.2, except that the model we are looking at
is slightly different from the one considered in that paper. The
definition of tail is made in such a way that $f$ will be separable
exactly over the tails of $X$. The proof of Lemma \ref{uconnlem}
relies on Condition \ref{cond1}. If one does not assume the condition,
it will not be true in general that $f$ will be exactly separable over
the tails.  For counter-examples see for example \cite{AbrOort98} for
$p=2$ or \cite{Saidicr} for general $p$.

\bigskip
\proof
We split the proof up in two parts. First we consider the statement of
the lemma for components $W$ which are not contained in $\Lambda$. Then
we show that $f$ is inseparable for all components contained in
$\Lambda$.

Let $W$ be a component of $X$ such that $f|_W$ is separable and $W$ is
not contained in $\Lambda$. The proof of \cite{Raynaud98}, Proposition
2.4.8, carries over to this situation and shows that $W$ is connected
to the rest of $X$ in a single point, i.e.\ $W$ is a tail.

Now suppose that $W$ is a tail of $X$. Let $Z$ be a component of $Y$
which maps to $W$. In case there is a branch point specializing to
$W$, the cover is separable over $W$ by Lemma \ref{inertia1lem}.(c)
and Condition \ref{cond1}.(b). Suppose there are no branch points
specializing to $W$ and $Z\to W$ is inseparable.  Then $I(Z)$ is a
nontrivial $p$-group which is a normal subgroup of $D(Z)$. Let
$Z'=Z/I(Z)$. Since $p$ exactly divides the order of $G$, it follows
that $D(Z)/I(Z)$ is of order prime-to-$p$.  Then $Z'\to W$ is Galois
of prime-to-$p$ order and branched at at most one point, hence
trivial. Since $Z\to Z'=W$ is purely inseparable, $Z$ is a component
of $Y$ of genus zero which meets the rest of $Y$ in a single point. By
assumption, there is no ramification point specializing to $Z$.  This
contradicts the minimal character of $Y$. Hence $f|_W$ is separable.

We are now going to prove the lemma for components contained in
$\Lambda$, i.e.\ we have to show that $f$ is inseparable over the
components contained in $\Lambda$. 

Suppose that $X_0$ is smooth. The cover $f$ is inseparable over
$X_0$, by Proposition \ref{admredprop}. This finishes the proof of the
lemma for $X_0$ smooth.

Suppose that $X_0$ is singular and suppose that there exists a
component of $\Lambda$ over which $f$ is separable.  Let $U$ be the
union of the components of $X$ over which the cover is inseparable. If
$f$ is separable over both $W_1$ and $W_2$, then the reduction is
admissible by Proposition \ref{admredprop}. It is no restriction to
suppose that $f$ is inseparable over $W_1$. Let $U'$ be the connected
component of $U$ containing $W_1$. The assumption on $\Lambda$ implies
that $U'$ does not contain $W_2$.  The number of branch points
specializing to $U'$ is at most two. Let $\bar{U}'$ be the union of
$U'$ and the components of $X$ which are adjacent to $U'$.

We orient the dual graph of $\bar{U}'$ starting from $W_1$.  Let
$\BB'$ be the set of tails of $\bar{U'}$. Let $b_{0}$ be the unique
component of $\bar{U}'$ on the geodesic connecting $W_1$ and $W_2$
over which $f$ is separable. By assumption, such a component exist.
For $b\in \BB'-\{b_{0}\}$, we say that the corresponding tail $X_b$ is
{\sl primitive} if one of the branch points specializes to this
component. Otherwise, we will call the tail {\sl new}.  Denote the set
of primitive (resp.\ new) tails by $\BB'_{{\rm prim}}$ (resp.\ 
$\BB'_{{\rm new}}$). The tail $X_b$ of $\bar{U}'$ meets the rest of
$\bar{u}'$ in a unique point $x_b$. Let $y_b$ be a point of $Y$ mapping to
$x_b$. By Lemma \ref{inertia1lem}, $I(y_b)$ is an extension of a
cyclic group of order $n_b$ prime-to-$p$, by a cyclic group of order
$p$. Denote the conductor by $h_b$, and put $\sigma_b=h_b/n_b$.
Recall that the ramification at $y_b$ is wild. This means that
$I(y_b)$ is a subgroup of $N_G(P)$. We have assumed (Condition
\ref{cond1}) that $N_G(P)$ is a dihedral group. It follows from
\cite{Raynaud98}, Lemme 1.1.2, that $\sigma_b\equiv 1/2\pmod{\ZZ}$. In
particular, for $b\in \BB'$, we have $\sigma_b\geq 1/2$.

Analogous to \cite{Raynaud98}, Section 3.4, one proves the following
vanishing cycle formula:
\[\sum_{b\in
  \BB'_{{\rm new}}} (\sigma_b-1)=-2 +\sum_{b\in \BB'_{{\rm prim}}}
(1-\sigma_b) +(1-\sigma_{b_0}).
\] 
(One proves this formula by constructing an ``auxiliary cover'' $Z'\to
\bar{U}'$ and showing that it can be lifted to characteristic zero.
The vanishing cycle formula follows then from the Riemann--Hurwitz
formula applied to the generic fiber of the lift.)  Furthermore, for
$b\in \BB'_{{\rm new}}$ we have that $\sigma_b-1\geq 1/2$. This
follows from a genus consideration, \cite{Raynaud98}, Proposition
3.3.5. The inequalities for $\sigma_b$ together with the fact that
$|\BB'_{\rm prim}|\leq 2$, implies that
\[
   \frac{|\BB'_{\rm new}|}{2}\leq-2+\frac{|\BB'_{\rm prim}|+1}{2}\leq
   -\frac{1}{2}.
\]
Which is impossible. This concludes the proof.
\Endproof

\begin{rem}
The notation and results used in the above proof will be introduced
and explained in more detail in the next section. The reason for
introducing the notation twice is that now that we proved Lemma
\ref{uconnlem}, the notation can be simplified considerably.
\end{rem}

\subsection{The auxiliary cover}\label{auxsec}

In this subsection we will suppose that Condition \ref{cond1} is
satisfied.  Furthermore, we suppose that the cover $f_K$ has bad
reduction.  We start by introducing some more notation.  Similar
notation is used in the proof of Lemma \ref{uconnlem}.

Let $\BB$ be the set of tails of $X$.  Every tail $X_b$ of $X$
contains a unique singular point $x_b$ of $X$. Choose a singular point
$y_b$ of $Y$ mapping to $x_b$ and let $Y_b$ be the component of $Y$
through $y_b$ which is mapping to $X_b$. Denote by $h_b$ the conductor
of $Y_b\to X_b$ at $y_b$, and by $n_b$ the order of the prime-to-$p$
ramification. Put
$\sigma_b=h_b/n_b$; this is the jump in the higher ramification groups
of $D({y_b})$, in the upper numbering.

Let $P$ be a $p$-Sylow group of $G$.  Let $U$ be the union of all
components of $X$ over which $f$ is inseparable. Let $V$ be any
connected component of $f^{-1}(U)$. Let $y$ be a singular point of $V$
and $Z_1$ and $Z_2$ be the components passing through $y$. By Lemma
\ref{uconnlem}, the inertia groups $I({Z_1}),I({Z_2})$ are cyclic of
order $p$, since $Z_1$ and $Z_2$ do not map to tails of $X$.
Therefore, by Lemma \ref{inertia2lem}, we have that
$I({Z_1})=I({Z_2})$ and both are equal to the $p$-part of $I(y)$. Here
we use that $p$ exactly divides the order of $G$. It follows that all
irreducible components of $V$ have the same $p$-Sylow subgroup of $G$
as inertia group. Therefore, there is a connected component $V$ of
$f^{-1}(U)$ such that $I(V)=P$. We will always assume $V$ to be chosen
like this.
  
As a consequence of Lemma \ref{uconnlem}, the construction of {\sl
auxiliary covers} as in \cite{Raynaud98}, Section 3.2, goes through in
our slightly different context. Since our notation differs from the
notation of \cite{Raynaud98}, we will recall the result.

\begin{prop} \label{auxprop}
  There exists a cover $g_K\!:Z_K\to X_K$ which is Galois with group
  $D(V)$ and has a quotient model $g_R\!:Z_R\to X_R$. It is uniquely
  characterized by the following properties:
\begin{itemize}
\item[(a)] There exists a suitable open $\Omega$ of $X_R$, which
  contains $X_b-\{x_b\}$ for $b\in \BB$, such that $Z_R\to X_R$ is
  tamely ramified over $\Omega$ and unramified outside the sections
  $x_i$ ($i=1\ldots 4$).
\item[(b)] There exists an \'etale neighborhood $X'_R\to X_R$ of
  $U\subset X\subset X_R$ such that
  $Y'_R\to X'_R\simeq \Ind_{D(V)}^G(Z'_R\to X'_R)$, where
  $Z'_R=Z_R\times_{X_R} X'_R$ and $Y'_R=Y_R\times_{X_R}X'_R$.
\end{itemize}
  We call $g_K$ the auxiliary cover associated to $f_K$.
\end{prop}

\proof
\cite{Raynaud98}, Proposition 3.2.6.
\Endproof

The (special fiber of the) auxiliary cover looks as follows. As above,
we let $V$ be a connected component of $f^{-1}(U)$ with inertia group
$P$. Restricted to $U$, the auxiliary cover is just $V\to U$. Let
$X_b$ be a tail of $X$ and $x_b$ the unique point of $X_b$ which is
singular in $X$; we may suppose $x_b=\infty$. Let $\Delta_b$ be the inertia
group of a point $y_b$ above $x_b$ which lies on $V$. Then, by the
Katz--Gabber Lemma \cite{Katz86}, there exist a cover $Z_b\to X_b$
unbranched outside $0, \infty$ and at most tamely ramified at 0 which
locally around $\infty$, if we induce it up to a $G$-cover, agrees
with $Y\to X$. Now $Z\to X|_{X_b}=\Ind_{\Delta_b}^{D(V)} (Z_b\to X_b)$.  Note
that by construction, the $\sigma_b$ for $b\in\BB$ for the cover
$g\!:Z\to X$, are the same as for the original cover.

\begin{lem}\label{sigmalem}

For $b\in \Bn$, we have $\sigma_b-1\geq  1/2.$
\end{lem}

\proof This follows from \cite{Raynaud98}, Proposition 3.3.5 and Lemme
1.1.2, and the assumption that $N_G(P)$ is a dihedral group. \Endproof

The following formula reflects the condition that the genus of $Y$ has
to be equal to the genus of the generic fiber $Y_K$. It is proved
in \cite{Raynaud98}, Section 3, using the auxiliary cover $g_R\!:Z_R\to
X_R$. 

\begin{prop}[Vanishing  cycle formula]\label{vcprop}
\[
   \sum_{b\in\Bn}(\sigma_b-1)=-2+\sum_{b\in\Bp}(1-\sigma_b).
\]
\end{prop}

\subsection{Modular reduction} \label{modularred}

Let $f_K\!:Y_K\to X_K$ be a $G$-cover for which Condition \ref{cond1}
holds.  In this subsection we will show that $f_K$ has either good
reduction or {\sl modular reduction}. Essentially, the property of
having modular reduction means that the auxiliary cover introduced in
the previous section is a cover $Z_K\to\PP^1_K$ with Galois group a
dihedral group and ramification of order 2 as discussed in Example
\ref{ellexa}.

\begin{defn}\label{modulardef} 
  Suppose $f_K$ has bad reduction. We will say that $f_K$ has {\em
    modular reduction} if the following conditions are satisfied.
\begin{itemize}
\item[(a)] The curve $X$ has four primitive tails and no new tails.
  Every irreducible component of $X$ is either an original component
  or a tail.
\item[(b)] Let $E$ be a connected component of $f^{-1}(X_0)$. Then
  $D(E)$ is a dihedral group of order $2N$ for some $N$ divisible by
  $p$.
\item[(c)] Let $X_b$ be a tail of $X$. Let $x_b$ be the unique
  singular point of $X_b$ in $X$ and let $y_b$ be a point of $Y$ mapping
  to $x_b$. Then the inertia group $I(y_b)$ is a dihedral group of order
  $2p$ and $\sigma_b=1/2$.
\end{itemize}
If $f_K\!:Y_K\to X_K$ has modular reduction, then we will say that the
special fiber $f\!:Y\to X$ of the quotient model of $f_K$ is of {\sl
  modular type}.
\end{defn}

The integer $N$ defined above will be called the {\sl level} of $f$.
This terminology reflects the relation between $f$ of modular type and
modular curves.

\begin{rem} 
  Suppose that $f_K\!:Y_K\to X_K$ has modular reduction of level $N$.
  The auxiliary cover $g_K\!:Z_K\to X_K$ corresponding to $f_K$ is
  a Galois cover with Galois group the dihedral group $\Delta$ of
  order $2N$, branched at $x_1,\ldots, x_4$ of order 2. The reduction
  $g\!:Z\to X$ is inseparable over $X_0$; it is separable over the
  tails. Note that $g^{-1}(X_0)$ may be identified with $E$ in the
  definition; it is, after choice of a base point, a generalized
  elliptic curve.
\end{rem}

\begin{figure}[htb]
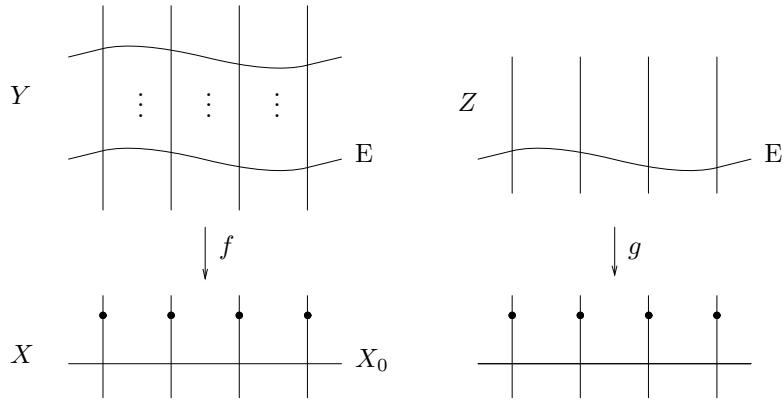
 \label{figpage}
\begin{center}
\include{modaux}
\end{center}
\caption{\label{modauxfig} modular reduction of level $N=p$}
\end{figure}

\begin{prop}\label{modularprop} 
  Let $f_K\!:Y_K\to X_K$ be as before, in particular we suppose that
  Condition \ref{cond1} is satisfied and that $f_K$ has bad reduction.
  Then $f_K$ has modular reduction.
\end{prop}

\proof Suppose that $f_K$ has bad reduction.
Proposition  \ref{vcprop} implies that 
\[\sum_{b\in \Bn}\sigma_b -1=-2+\sum_{b\in\Bp}1-\sigma_b.\]
Recall that $\sigma_b\geq 3/2$ for $b\in \Bn$ and $\sigma\geq 1/2$
for $b\in\Bp$, Lemma \ref{sigmalem}. Hence $|\Bn|/2\leq 0$. This
implies that $|\Bn|=0$ and $\sigma_b=1/2$ for $b\in \Bp$. Moreover,
all $x_i$ specialize to components over which $f$ is separable, so
$|\Bp|=4$.

 The decomposition group of a singular point $y_b$ of $Y$
contains a dihedral group of order $2p$, since 
$\sigma_b=1/2$.  (It cannot be Abelian, since then $\sigma_b$ would
be an integer by the Hasse--Arf Theorem.)  Since $D(y_b)=I(y_b)$ is a
cyclic-by-$p$ group, $D(y_b)$ is isomorphic to a dihedral group of
order $2p$. This proves Part (c) of Definition \ref{modulardef}.

The decomposition group $D(V)$ is a subgroup of $N_G(P)$ which is a
dihedral group, by assumption. Note that $D(V)$ is not cyclic, since it
contains (a conjugate of) the inertia group of some point $y_b$,
which is dihedral. This proves Part (b).

The only thing left to show is the second part of Part (a). Let $W$ be
a component of $X$ which is neither a tail nor an original component
and let $Z$ be a component of $Y$ which maps to $W$. Above we have
shown that there are exactly four tails. This implies that $W$ meets
the rest of $X$ in two points. There is no branch point specializing
to $W$, so the maximal separable subcover $Z'\to W$ of $Z\to W$ is
branched at at most two points.  Moreover, $I_Z$ is a cyclic group of
order $p$, since $W$ is not a tail. This implies that the degree of
$Z'\to W$ is prime-to-$p$, so it is a cyclic cover of $\PP^1$ branched
at two points.  But then $Z$ has genus zero and meets the rest of $Y$
in exactly two points. This contradicts the minimality of $Y$.

This shows that the cover has modular reduction.
\Endproof

\subsection{Reduction of degenerate covers} \label{boundarypts}

In this section we will study the reduction behavior of the degenerate
covers corresponding to the cusps of the Hurwitz space. Note that the
theory of semistable reduction we developed so far does not apply
here, since we always assumed that the generic fibers of our curves
were smooth.

The situation in this section is somewhat different from that in the
previous sections. Let $R$ be a complete discrete valuation ring whose
quotient field $K$ is of characteristic zero and whose residue field
$k=\bar{k}$ has characteristic $p$. Let $(X_K;x_1,\ldots,x_4)$ be a
stably marked $K$-curve of genus zero, which we suppose to be
singular. Let $X_{1,K}$ and $X_{2,K}$ be the two irreducible
components of $X_K$.  Let $\tau$ be the singular point of $X_K$. Let
$f_K\!:Y_K\to X_K$ be an admissible $G$-Galois cover ramified at
$x_1,\ldots ,x_4$.  Let $\rho$ be a point of $Y_K$ mapping to $\tau$
and let $Y_{1,K}$ and $Y_{2,K}$ be the components of $Y_K$ passing
through $\rho$, where we suppose that $Y_{i,K}$ maps to $X_{i,K}$. Let
$G_i$ be the decomposition group $Y_{i,K}$.  Let $X_{i,R}$ and
$Y_{i,R}$ be the closure of $X_{i,K}$ and $Y_{i,K}$ in $X_R$. Since
$\bar{\H}_4^{\rm in}(G)$ is proper, $f_K$ extends uniquely to a map
$f_{0,R}\!:Y_R\to X_{0,R}$ between stably marked curves over $\Spec
R$.  Let $f_i\!:Y_{i,R}\to X_{i,R}$ be the corresponding morphism. We
will denote the special fibers of $Y_R,X_R,Y_{i,R},X_{i,R}$ by
$Y,X,Y_i,X_i$, respectively. The mark $C$ on $X_R$ can be written
as a union $C=C'_1\cup C'_2$, where $C'_i$ is a mark on $X_{i,R}$.
Let $\tau$ be the unique singular point of $X_K$, we denote its
(unique) extension to $X_R$ also by $\tau$. Let
$C_i=C'_i\cup\{\tau\}$.  Let $D'_i$ be the restriction of the mark
of $Y_R$ to $Y_{i,R}$ and write $D_i$ for the union of $D'_i$ with the
points of $Y_{i,R}$ which are singular in $Y_{R}$. The next lemma
follows immediately.

\begin{lem}\label{halff}
  The morphism $f_{i,R}\!:(Y_{i,R},D_i)\to (X_{i,R},C_i)$ constructed
  above is a morphism of stably marked curves.
\end{lem}

The covers $f_{i,K}$ constructed above are covers of a projective line
branched at three points. So for these covers we can apply the
criterion for good reduction proved by Raynaud \cite{Raynaud98}.
Actually, since here we put a stronger condition on $G$ than in
\cite{Raynaud98}, we can show that the covers $f_{i,K}$ have good
reduction iff $p$ does not divide the ramification indices, Lemma
\ref{3pt}. In other words, the condition on the field of definition in
the Theorem of Raynaud is not needed in this case.

Proposition \ref{redbpts} is the key result in this section. It
describes the reduction of the degenerate covers to characteristic
$p$. The idea is that we can understand the reduction of such a
degenerate cover, because we understand the reduction of the two three
point covers of which it is made. Proposition \ref{redbpts} will be
used in Section 5 to explicitly describe the bad part of the Hurwitz
space.

\begin{lem}\label{3pt}
  Let $f_K\!:Y_K\to X_K\cong\PP^1_K$ be a $G$-cover branched at three
  points $0,1,\infty$. Suppose that $p$ exactly divides the order of
  $G$ and the normalizer of a $p$-Sylow group of $G$ is dihedral.
  Denote the ramification indices of $f_K$ by $m_1,m_2,m_3$. Let
  $f_R\!:Y_R\to X_R$ be the quotient model corresponding to
  $X_{0,R}=\PP^1_R$. Then $f_R$ has good reduction if and only if
  $p\nmid m_i$ for $i=1,2,3$.
\end{lem}

\proof
This follows immediately from the vanishing cycle formula
\cite{Raynaud98}, Section 3.4, combined with the estimates for
$\sigma_b$ from Lemma \ref{sigmalem}.
\Endproof

\begin{prop} \label{redbpts}
  Let $f_K\!:Y_K\to X_K$ be as above and let $f\!:Y\to X$ be
  its reduction.  Let $n$ be the order of the ramification of $f_K$
  above the singular point $\tau$. We denote by $\tau_k$ the image
  of $\tau$ on the special fiber. Let $\rho_k$ be a point on $Y$
  above $\tau_k$.
\begin{itemize}
\item[(a)] The cover $f_K$ has admissible reduction iff $p\nmid n$.
\item[(b)] The inertia group  $I(\rho_k)$ has order $n$
\item[(c)] If $f_K$ has bad reduction, then $f$ is of modular type.
\item[(d)] Suppose $f_K$ has bad reduction, then $n$ divides the level
  $N$ of $f$.  Let $Z_1$ and $Z_2$ be the irreducible components of
  $Y$ passing through $\rho_k$. Then $D({Z_1})$ and $D({Z_2})$ are
  dihedral groups of order $2n$.
\end{itemize}
\end{prop}

\proof 
The cover $f\!:Y\to X$ is the reduction of a degenerate cover, by
assumption. However, there will be covers representing points in the
interior of the Hurwitz stack in characteristic zero which specialize
to $f$. Therefore, in case $f_K$ has bad reduction, $f$ will be of
modular type  by Proposition \ref{modularprop}. This proves Part (c).

This implies that in case $f_K$ has bad reduction, $f_0\!:Y\to X_0$
does not contract any components to $\tau$. This is clearly also the
case if $f_K$ does not have bad reduction.  Since $f_i\!:Y_{i,R}\to
X_{i,R}$ for $i=1,2$ are flat, we have $g(Y_i)=g(Y_{i,K})$. Moreover,
$g(Y)=g(Y_K)$. Now
\[
     g(Y)=(|G|/|G_1|)g(Y_1)+(|G|/|G_2|)g(Y_2)+1-
           |G|/|G_1|-|G|/|G_2|+|G|/|I(\rho_k)|
\]
and
\[
    g(Y_K)=(|G|/|G_1|)g(Y_{1,K})+(|G|/|G_2|)g(Y_{2,K})+
           1-|G|/|G_1|-|G|/|G_2|+|G|/n.
\]
This shows that $|I(\rho_k)|=n$. 

Suppose $p|n$. Then Part (b) of Lemma \ref{inertia1lem} implies that
the covers $f_{i,K}$ have bad reduction for $i=1,2$. It follows that
$f_K$ has bad reduction. Conversely, suppose that $f_K$ has bad
reduction. Part (c) implies that $f_K$ has modular reduction. In
particular it follows that $I({Z_i})$ is $p$-cyclic. Part (b) of Lemma
\ref{inertia1lem} implies that $p|n$.

Suppose that $f_K$ has bad reduction.  The decomposition groups
$D(Z_i)$ are subgroups of $N_G(P)$, since the corresponding inertia
groups are nontrivial. Moreover, they contain a dihedral group of
order $2p$, by the definition of modular reduction (Definition
\ref{modulardef}). It follows that the $D(Z_i)$ are dihedral
groups. The maximal separable subcover of $Z_i\to X_i$ is
$f'_i\!:Z'_i:=Z_i/I({Z_i})\to X_i$. Note that $f'_i$ is branched at
three points of order $2,2,n/p$ hence $g(Z'_i)=0$. It follows that the
degree of $f'_i$ is $2n/p$, hence $D(Z_i)$ is a dihedral group of
order $2n$. By definition of the level $N$ of $f$, the subgroup of
$N_G(P)$ generated by $D(Z_1)$ and $D(Z_2)$ is a dihedral group of
order $2N$; this implies that $n|N$.  \Endproof

%\bibliographystyle{abbrv}
%\bibliography{hurwitz}

%\end{document}

%% file: modaux.tex
\setlength{\unitlength}{0.0004in}
{\renewcommand{\dashlinestretch}{30}
\begin{picture}(9486,5499)(0,-10)
\put(924,1377){\blacken\ellipse{90}{90}}
\put(924,1377){\ellipse{90}{90}}
\put(1824,1377){\blacken\ellipse{90}{90}}
\put(1824,1377){\ellipse{90}{90}}
\put(2724,1377){\blacken\ellipse{90}{90}}
\put(2724,1377){\ellipse{90}{90}}
\put(3624,1377){\blacken\ellipse{90}{90}}
\put(3624,1377){\ellipse{90}{90}}
\path(924,297)(924,1647)
\path(1824,297)(1824,1647)
\path(2724,297)(2724,1647)
\path(3624,297)(3624,1647)
\put(6324,1377){\blacken\ellipse{90}{90}}
\put(6324,1377){\ellipse{90}{90}}
\put(7224,1377){\blacken\ellipse{90}{90}}
\put(7224,1377){\ellipse{90}{90}}
\put(8124,1377){\blacken\ellipse{90}{90}}
\put(8124,1377){\ellipse{90}{90}}
\put(9024,1377){\blacken\ellipse{90}{90}}
\put(9024,1377){\ellipse{90}{90}}
\path(474,747)(4074,747)(4074,747)
\path(924,5472)(924,2772)
\path(1824,5472)(1824,2772)
\path(2724,5472)(2724,2772)
\path(3624,5472)(3624,2772)
\path(6324,1647)(6324,297)
\path(7224,1647)(7224,297)
\path(8124,1647)(8124,297)
\path(9024,1647)(9024,297)
\path(5874,747)(9474,747)
\path(5874,747)(9474,747)
\path(6324,4797)(6324,2997)
\path(7224,4797)(7224,2997)
\path(8124,4797)(8124,2997)
\path(9024,4797)(9024,2997)
\path(2274,2547)(2274,1872)
\path(2244.000,1992.000)(2274.000,1872.000)(2304.000,1992.000)
\path(7674,2547)(7674,1917)
\path(7644.000,2037.000)(7674.000,1917.000)(7704.000,2037.000)
\spline(474,3447)
(1374,3672)(3174,3222)(4074,3447)
\put(4250,3400){E}
\spline(474,4797)
(1374,5022)(3174,4572)(4074,4797)
\spline(5874,3447)
(6774,3672)(8574,3222)(9474,3447)
\put(9650,3400){E}
\put(1374,4000){$\vdots$}
\put(2274,4000){$\vdots$}
\put(3174,4000){$\vdots$}
\put(4250,700){$X_0$}
%\put(750,-100){$X_1$}
%\put(1650,-100){$X_2$}
%\put(2550,-100){$X_3$}
%\put(3450,-100){$X_4$}
\put(-300,747){$X$}
\put(-300,4167){$Y$}
\put(5604,4077){$Z$}
\put(2454,2187){$f$}
\put(7854,2187){$g$}
\end{picture}
}

%%% Local Variables: 
%%% mode: latex
%%% TeX-master: t
%%% End: 

%% file: sec4.tex
%%% Local Variables: 
%%% mode: latex
%%% TeX-master: t
%%% End: 

%\documentclass{article}
%\usepackage{def}
%\begin{document}

\section{A Reduction Theorem} \label{reductionthm}

In this section we prove a theorem about the structure of
$\hb\otimes\FF_p$, where $\hb$ is a complete Hurwitz space for
$G$-covers satisfying Condition \ref{cond1}. As explained in the
introduction, this theorem relies on, and in some sense extends, the
results of Katz and Mazur on the reduction of the modular curve
$X_1(p)$. This is somewhat surprising. Like modular curves, the
Hurwitz spaces we look at are curves, defined over small number
fields, and arise as quotients of the upper half plane by discrete
subgroups of ${\rm GL}_2(\ZZ)$, see \cite{Fried87}. However, these
groups are non-congruence subgroups, in general.

\subsection{Statement of the main results} \label{statement}

\subsubsection{} \label{statement1}

Let $G$ be a finite group, $\Cl = (C_1,C_2,C_3,C_4)$ a class vector in
$G$ of length $4$ and $p$ an odd prime. We denote by $m_i$ the order
of the elements of $C_i$. We assume that Condition \ref{cond1} holds,
with respect to $G$, $p$ and $m_i$. 

Let $K:=K(\Cl)\subset\QQ(\zeta_m)$ be the minimal field over which the
classes $C_i$ are rational, see Section \ref{morevariants}. We choose
a prime ideal $\p$ of $K$ dividing $p$, and denote by
$\Lambda:=\O_{K,\p}$ its local ring. We let $\hb:=\hb_4\inn(\Cl)$ be
the complete Hurwitz space over $\Lambda$, as defined in Section
\ref{coarse}. By construction, $\hb$ is an algebraic space, proper and
of finite type over $\Lambda$.  It contains the Hurwitz space
$H:=H_4\inn(\Cl)$ as a dense open subscheme. The scheme $H$ is smooth
over $\Lambda$; its generic fiber $H\otimes K$ is the reduced Hurwitz
curve studied e.g.\ in \cite{Fried87}.
 
Let $\hb\bad\subset\hb$ be the closed subspace corresponding to bad
covers. Since bad covers occur only in positive characteristic,
$\hb\bad$ is a closed subspace of $\hb\otimes\FF_q$, where $\FF_q$ is
the residue field of $\Lambda$.  The complement $\hb\adm:=\hb-\hb\bad$
is a dense open subscheme and corresponds to admissible covers. Let
$\hb\good$ be the closure of $\hb\adm\otimes\FF_q$ in
$\hb\otimes\FF_q$. Note that this is an abuse of notation, since
$\hb\good$ has nontrivial intersection with $\hb\bad$, in general.
There is a natural map $\hb\to\PP^1_\lambda\otimes\Lambda$.  We say
that an $\bar{\FF}_p$-rational point $s$ of $\hb\bad$ is {\em
supersingular} if the corresponding value $\lambda(s)\in\bar{\FF}_p$
is supersingular.

\begin{thm} \label{redthm}\ 
  \begin{enumerate}
  \item 
    The complete Hurwitz space $\hb$ is a normal scheme of
    dimension $2$. The natural map $\hb\to\PP^1_\lambda\otimes\Lambda$
    is finite and flat.
  \item 
    The subspaces $\hb\bad$ and $\hb\good$ are smooth projective
    curves over $\FF_q$. They intersect transversally in the
    supersingular points.
  \end{enumerate}
\end{thm}

In the rest of Section \ref{statement} we state a number of results on
the local structure of $\hb$ and explain how Theorem \ref{redthm} can
be deduced from them. Let us give a brief outline. On the open subset
$\hb\adm\subset\hb$, Theorem \ref{redthm} is known to hold, see
\cite{diss}. Therefore, it suffices to look at $\hb$ in a neighborhood of a
point corresponding to a bad cover. Theorem \ref{precisethm} below
describes the universal deformation ring of such a bad cover.
Essentially, this theorem implies that the complete Hurwitz stack
$\Hb$ associated to $\hb$ is regular and that a version of Theorem
\ref{redthm} holds for $\Hb$. To finish the proof of Theorem
\ref{redthm}, we have to study the {\em monodromy action}, i.e.\ the
action of the group of automorphisms of a bad cover on its universal
deformation ring. Propositions \ref{monoprop} describes this action,
and Theorem \ref{redthm} follows. Under some extra hypotheses
(Condition \ref{cond2}), we can improve our results on the monodromy
action, and we can actually show that $\hb$ is regular. The relevant
statement is made in Proposition \ref{innermonoprop}.

The proofs of Theorem \ref{precisethm}, Proposition \ref{monoprop} and
Proposition \ref{innermonoprop} are postponed to Section
\ref{proofs}.

\subsubsection{The universal deformation ring} \label{precise1}

Let $\Hb:=\Hb_4\inn(\Cl)$ be the complete Hurwitz stack over
$\Lambda$, associated to $\hb$, and let $\Hb\bad,\Hb\good\subset\Hb$
be the closed substacks corresponding to the closed subspaces
$\hb\bad,\hb\good\subset\hb$. We look at the following situation. Let
$k$ be an algebraically closed field of characteristic $p$ and
\[
      f_0:Y\To(X_0;x_i)
\]
an object of $\Hb\bad$, defined over $k$. In the rest of
this section we will mostly write ``$Y$'' instead of ``$f_0:Y\to
(X_0;x_i)$'' for this object; we understand that the curve $Y$ is
equipped with an action of the group $G$, a mark $D\subset Y$ and a
map $f_0$ to the stably marked curve $(X_0;x_i)$. 

Let $R_Y$ be the strict complete local ring of $\Hb$ at the $k$-point
corresponding to $Y$. Let $Y\univ$ be the object of $\Hb_4\inn(\Cl)$
corresponding to the tautological morphism $\Spec R_Y\to\Hb$. By a
general property of algebraic stacks, $Y\univ$ is the {\em universal
  deformation} of $Y$ as object of $\Hb$.  This means the following.
Let $W(k)$ denote the ring of Witt vectors over $k$ and $\Cd_k$ the
category of complete local Noetherian $W(k)$-algebras with residue
field $k$. We let $\Def(Y)$ be the functor which assigns to
$R\in\Cd_k$ the set of isomorphism classes of deformations of $Y$ over
$R$. Then $Y\univ$ defines an equivalence
\[
      \Hom_{\Cd_k}(R_Y,\;\cdot\;) \liso \Def(Y).
\]
The closed substack $\Hb\bad\subset\Hb$ (resp.\ $\Hb\good\subset\Hb$)
corresponds to a subfunctor $\Def(Y)\bad\subset\Def(Y)$ (resp.\ 
$\Def(Y)\good\subset\Def(Y)$), which is represented by a quotient ring
$R_Y\bad$ (resp.\ $R_Y\good$) of $R_Y$.

We denote by $\Def(X_0;x_i)$ the deformation functor for $(X_0;x_i)$ as
object of $\Mb_{0,4}$ and by $R_{X_0}$ the universal deformation ring.
The morphism $\Hb\to\Mb_{0,4}$ induces a transformation
$\Def(Y)\to\Def(X;x_i)$, hence a morphism $R_{X_0}\to R_Y$. The ring
$R_{X_0}$ is of the form $R_{X_0}=W(k)[[w]]$. For instance, if $X_0$
is smooth we may assume that $X_0=\PP^1$, $x_1=0$, $x_2=1$,
$x_3=\infty$, $x_4=\lambda_0\in k-\{0,1\}$ and $w=\lambda-\lambda_0$.
We say that $Y$ is {\em supersingular} (resp.\ {\em ordinary}) if
$\lambda_0$ is supersingular (resp.\ ordinary). If $X_0$ is singular,
we will also say that $Y$ is ordinary.

\begin{thm} \label{precisethm}
  The ring $R_Y$ is regular of dimension $2$ and a finite flat
  extension of $R_{X_0}$. Moreover, there exists a regular sequence
  $(t,\pi)$ for $R_Y$ such that
  \begin{enumerate}
  \item
    $R_Y\bad=R_Y/(\pi)\cong k[[t]]$,
  \item 
    The induced (finite) morphism $R_{X_0}\otimes_{W(k)}k\to
    R_Y\bad$ is inseparable of degree $p$. Its separable part
    $R_{X_0}\otimes_{W(k)}k\to (R_Y\bad)^p\cong k[[t^p]]$ is tamely
    ramified (of degree $d$); if $X_0$ is smooth then $d=1$.
  \item
    If $Y$ is ordinary, then $p=u\,\pi^{p-1}$ for a unit $u\in
    R_Y^\times$ and $R_Y\good=0$.
  \item If $Y$ is supersingular, then $p=u\,\pi^{p-1}$, where $u\in R_Y$
    is a local equation for $\Hb_4\inn(\Cl)\good$. More precisely,
    $R_Y\good=R_Y/(u)\cong k[[\bar{\pi}]]$, and $(u,\pi)$ is a regular
    sequence for $R_Y$.
  \end{enumerate}
\end{thm} 

In Section \ref{dihedral} we will prove this theorem in the case that
the group $G$ is dihedral and the conjugacy classes $C_i$ represent
reflections. In this case, the $G$-cover $f:Y\to X_0$ corresponds
essentially to a point on the modular curve $X_1(N)$, and Theorem
\ref{precisethm} follows from the results of \cite{KatzMazur}. We will
prove the general case in Section \ref{deformation} by reduction to
the dihedral case.

\begin{cor} \label{precisecor1}
  The complete Hurwitz stack $\Hb$ is regular of dimension $2$. The
  natural morphism $\Hb\to\Mb_{0,4}$ is finite and flat. The closed
  substacks $\Hb\bad$ and $\Hb\good$ are smooth over $\FF_q$, and
  intersect transversally in the supersingular point. 
\end{cor}

\subsubsection{The monodromy action} \label{precise2}

Corollary \ref{precisecor1} together with Proposition \ref{coarseprop}
implies Part (i) of Theorem \ref{redthm}. To prove Part (ii)
of Theorem \ref{redthm}, we use the general fact that the coarse
moduli scheme is locally the quotient of (an \'etale cover of) the
corresponding algebraic stack by a finite group action. Actually, it
suffices to look at the strict complete local rings. Hence we are led
to study the monodromy action of the automorphisms of $Y$ on the
universal deformation ring $R_Y$. 

Let $Y$ be as in Section \ref{precise1}. We write $\Aut_k(Y)$ for
the group of $k$-linear automorphisms of $Y$, considered as object of
$\Hb$. Since $(X_0;x_i)$ has no nontrivial automorphism, an element
$\sigma\in\Aut_k(Y)$ is a $k$-automorphism of $Y$ such that
$f_0\circ\sigma=f_0$ and $g\circ\sigma=\sigma\circ g$ for all $g\in G$.
By the universal property of $Y\univ$, for each $\sigma$ there exists a
unique automorphism $\gamma:R_Y\iso R_Y$ such that $\sigma$ lifts to a
unique $\gamma$-semilinear automorphism $\sigma\univ:Y\univ\iso
Y\univ$. We call $\gamma\in\Aut_{\Cd_k}(R_Y)$ the {\em monodromy
  action} of $\sigma$. We define the {\em monodromy group} $\Gamma$
of $Y$ as the image of the homomorphism
$\Aut_k(Y)\to\Aut_{\Cd_k}(R_Y)$.  Let $Y_{\bar{\eta}}$ be the
geometric generic fiber of $Y\univ$ (note that $R_Y$ is a domain, by
Theorem \ref{precisethm}). Identifying the group of
$\bar{\eta}$-automorphisms of $Y_{\bar{\eta}}$ with $C_G$, the center
of $G$, we obtain a natural exact sequence
\begin{equation} \label{monseq1}
   1 \To C_G \To \Aut_k(Y) \To \Gamma \To 1.
\end{equation}
Let $s:\Spec k\to\hb$ be the geometric point of the Hurwitz space
corresponding to $Y$. The strict complete local ring of $\hb$ is the
ring of $\Gamma$-invariants of $R_Y$:
\begin{equation} \label{inveq}
       \Od_{\hb,s} \;=\; R_Y^{\textstyle \;\Gamma}.
\end{equation}

\subsubsection{The absolute monodromy action} \label{precise3}

To study the action of $\Gamma$ on $R_Y$, it will turn out to be
useful to enlarge $\Gamma$ and look at the {\em absolute monodromy
  group} $\Gamma\ab$.

Let $\Hb_4\ab(\Cl)$ be the absolute version of the complete Hurwitz
stack $\Hb$. Let us for the moment consider $Y$ and $Y\univ$ as
objects of $\Hb_4\ab(\Cl)$. In other words, we forget the embedding
$G\inj\Aut(Y/X)$, and only retain its image, see Variant
\ref{variant1}. It is still true that $Y\univ$ is the universal
deformation of $Y$. But we obtain a bigger automorphism group, in
general. We denote by $\Aut_k\ab(Y)$ the group of $k$-automorphisms of
$Y$ as object of $\Hb_4\ab(\Cl)$. An element $\sigma\in\Aut_k\ab(Y)$
is a $k$-automorphism of $Y$ such that $f_0\circ\sigma=f_0$ and
$\sigma\circ g\circ\sigma^{-1}\in G$ for all $g\in G$. We find that
$G$ is a normal subgroup of $\Aut_k\ab(Y)$, and $\Aut_k(Y)$ is the
centralizer of $G$ in $\Aut_k\ab(Y)$. The group $\Aut_k\ab(Y)$ acts on
$R_Y$; this action extends the action of $\Aut_k(Y)$. We write
$\Gamma\ab$ for the image of $\Aut_k\ab(Y)$ in $\Aut_{\Cd_k}(R_Y)$.
Clearly, $\Gamma\subset\Gamma\ab$, and the exact sequence
\zgl{monseq1} becomes
\begin{equation} \label{monseq2}
    1 \To G \To \Aut_k\ab(Y) \To \Gamma\ab \To 1. 
\end{equation}

\begin{prop} \label{monoprop}
  There exist two characters $\chi\bad:\Gamma\ab\To\FF_p^\times$ and
  $\chi\adm:\Gamma\ab\To \mu_d$ (here $\mu_d\subset R_Y$ denotes the
  set of $d$th roots of unity, for $d$ as in Theorem \ref{precisethm}
  (ii)) with the following properties.  The parameters $t,\pi\in R_Y$
  in Theorem \ref{precisethm} can be chosen such that for all
  $\gamma\in\Gamma\ab$
  \[
        \gamma(t) = \chi\adm(\gamma)\cdot t,\qquad
        \gamma(\pi) \equiv \chi\bad(\gamma)\cdot\pi \pmod{\pi^2}.
  \]
  Moreover, the homomorphism 
  $(\chi\bad,\chi\adm):\Gamma\ab\To\FF_p^\times\times\mu_d$
  is injective.   
\end{prop}

We will prove this proposition in Section \ref{monoproof}. In the rest
of this subsection we will show that Theorem \ref{precisethm} and
Proposition \ref{monoprop} together imply Theorem \ref{redthm}. Since
$\Gamma\subset\Gamma\ab$, Proposition \ref{monoprop} holds also for
$\Gamma$.

The closed subscheme $\hb\bad\subset\hb$ is the image of the natural
morphism $\Hb\bad\to\hb$, by definition. It follows that the strict
complete local ring $\Od_{\hb\bad,s}$ is the image of the natural
morphism $R_Y^{\Gamma}\to R_Y\bad$. By Proposition \ref{monoprop}, the
order of $\Gamma$ is prime-to-$p$. Therefore, the natural map
$R_Y^{\Gamma}\to(R_Y\bad)^{\Gamma}$ is surjective, hence
$\Od_{\hb\bad,s}=(R_Y\bad)^{\Gamma}$, see \cite{KatzMazur}, A 7. Now
Theorem \ref{precisethm} (i) and Proposition
\ref{monoprop} imply
\begin{equation} \label{badinveq}
   \Od_{\hb\bad,s} \;=\;      (R_Y\bad)^{\textstyle\Gamma}
                   \;\cong\;  k[[t^{d'}]],
\end{equation}
where $d'|d$ is the order of $\chi\adm(\Gamma)$. We have shown that
$\hb\bad$ is a smooth curve. The same argument shows that
\begin{equation} \label{badinv3}
    \Od_{\hb\good,s} \;=\; (R_Y\good)^{\textstyle\Gamma}
    \;\cong\; \left\{\begin{array}{ll}
         0                     & \text{\rm if $Y$ is ordinary},\\
         k[[\bar{\pi}^\mu]]\quad & \text{\rm if $Y$ is supersingular},\\
              \end{array}\right.
\end{equation}
where $\mu|(p-1)$ is the order of $\chi\bad(\Gamma)$. We conclude that
$\hb\good$ is a smooth curve and that $\hb\good\to\PP^1_\lambda$ is
ramified of order $(p-1)/\mu$ in the supersingular points. Let us
assume that $Y$ is supersingular and denote by $\hb_p\red$ (resp.\ by
$\Hb_p\red$) the closed subscheme $(\hb\otimes\FF_p)\red\subset\hb$
(resp.\ the closed substack $(\Hb\otimes\FF_p)\red\subset\Hb$). By
Theorem \ref{precisethm} (iv), $\Hb_p\red\times_{\Hb}\Spec
R_Y=R_Y/(\pi u)\cong k[[\bar{\pi},\bar{u}\mid\bar{\pi}\bar{u}=0]]$.
Taking invariants and arguing as before, we get
\begin{equation} \label{badinv4}
  \Od_{\hb_p\red,s} \;=\; (R_Y/(\pi u))^{\textstyle\Gamma}
   \;\cong\; k[[\bar{\pi}^\mu,\bar{u}\mid\bar{\pi}^\mu\bar{u}=0]]
\end{equation}
(note that $\chi\adm=1$, since $X_0$ is smooth).  This completes the
proof of Theorem \ref{redthm}, modulo the proofs of Theorem
\ref{precisethm} and Proposition \ref{monoprop}. 

Translating the statements of Theorem \ref{precisethm} into geometric
properties of the map $\hb\to\PP^1_\lambda$, we obtain the following
corollary.

\begin{cor} \label{precisecor3}
  Let $s$ be the point on $\hb\bad$ corresponding to $Y$.
\begin{enumerate}
\item 
  The natural map $\hb\bad\to\PP^1_\lambda\otimes\FF_q$ is finite,
  with inseparability degree $p$. Its separable part
  $(\hb\bad)^{(p)}\to\PP^1_\lambda\otimes\FF_q$ is tamely ramified at
  $\lambda=0,1,\infty$ and \'etale everywhere else. More precisely,
  its ramification index in $s^{(p)}$ is equal to
  $[\Im(\chi\adm):\mu_d]$.
\item 
  The natural map $\hb\good\to\PP^1_\lambda\otimes\FF_q$ is tamely
  ramified at $\lambda=0,1,\infty$ and the supersingular values of
  $\lambda$, and is \'etale everywhere else. Its ramification index in
  $s\in \hb\good\cap\hb\bad$ is $[\chi\bad(\Gamma):\FF_p^\times]$
  (this happens if and only if $Y$ is supersingular).
\item 
  Let $m_s$ be the multiplicity of $\hb\bad$ in a neighborhood of $s$
  in $\hb\otimes\FF_q$, i.e.\ the length of the Artinian local ring
  $\Od_{\hb,\eta}/p\Od_{\hb,\eta}$, where $\eta:\Spec
  k((t^{d'}))\to\hb$ comes from equation \zgl{badinveq}. Then
  $m_s=[\chi\bad(\Gamma):\FF_p^\times]$.
\end{enumerate}
\end{cor}

\subsubsection{Regularity} \label{precise4}

We have a better control of the monodromy group $\Gamma$, if we assume
that --- in addition to Condition \ref{cond1} --- the following holds.
\begin{cond} \label{cond2}\ 
\begin{itemize}
\item[(e)]
  The center $C_G$ of $G$ is trivial.
\item[(f)]
  Let $G'$ be a subgroup of $G$ which contains an element of order $p$
  and an element of one of the classes $C_i$, $i=1,\ldots,4$. Then
  $G=G'$.
\end{itemize}
\end{cond}

\begin{prop} \label{innermonoprop}
  Assume that Condition \ref{cond2} holds. Let $N$ be the level of
  $Y$. Then $\chi\adm|_{\Gamma}=1$ and
  \[ 
    \Gamma \;\cong\;
      \left\{\begin{array}{cll} 
           \ZZ/2 \quad & \text{if}\;\; & N=p, \\ 
           1     \quad & \text{if}\;\; & N>p. \\ 
      \end{array}\right.  
  \]
\end{prop}

We will prove Proposition \ref{innermonoprop} in Section
\ref{innermonoproof}.

\begin{cor} 
  If Condition \ref{cond2} holds, then $\hb$ is regular.  
\end{cor}

\proof It follows from Theorem \ref{precisethm} and Propositions
\ref{monoprop} and \ref{innermonoprop} that $R_Y^{\Gamma}$ has a
regular sequence $(t,\pi')$, where $\pi'\sim\pi^\mu$ and $\mu=1,2$.
\Endproof

\begin{rem} \label{preciserem}
  We do not know of any example in which either $\hb$ or
  $\hb_4\ab(\Cl)$ is not regular. On the other hand, the Hurwitz
  spaces $\hb_{[4]}\inn(\Cl)$ and $\hb_{[4]}\ab(\Cl)$ behave like
  modular curves. We have to replace the $\lambda$-line with the
  $j$-line. For the special values $j=0$ and $j=1728$, the monodromy
  action may become more complicated then in Proposition
  \ref{monoprop}. Therefore, $\hb_{[4]}\inn(\Cl)$ and
  $\hb_{[4]}\ab(\Cl)$ are very often not regular. However, Theorem
  \ref{redthm} remains true.
\end{rem}

\subsection{Dihedral covers and generalized elliptic curves}
\label{dihedral}
  
In this section we show that Theorem \ref{precisethm} is true in the
case the group $G$ is a dihedral group of order $2N$, where
$p\,||\,N$, and $\Cl$ consists of $4$ times the conjugacy class of
reflections. In order to prove this, we relate the deformation theory
of a bad cover, which arises as the reduction of a $G$-cover of type
$\Cl$, to the deformation theory of a generalized elliptic curve
endowed with a certain level structure. Once this is achieved,
Theorem \ref{precisethm} follows from the results of \cite{KatzMazur}
on the reduction of the modular curve $X_1(p)$.

\subsubsection{Generalized elliptic curves} \label{review}

We start by recalling some definitions from \cite{DelRap},
\cite{KatzMazur} and \cite{Edix}. A {\em generalized elliptic curve}
over a scheme $S$ is a semistable curve $E/S$ of genus $1$ together
with a section $0:S\to E^{\rm sm}$ and an $S$-morphism $+:E^{\rm
  sm}\times_S E\to E$, verifying the following properties. The
geometric fibers of $E/S$ are either smooth or ``$n$-gons''. The
restriction of $+$ to $E^{\rm sm}$ gives $E^{\rm sm}$ the structure of
a commutative group scheme with identity $0$.  Moreover, $E^{\rm sm}$
acts on $E$ by ``rotation'', see \cite{DelRap}, Definition II.1.12. If
$E/S$ is smooth of genus $1$ and $0:S\to E$ a section, then there
exists one and only one such group law $+:E\times E\to E$, and we call
$(E,0)$ an {\em elliptic curve} over $S$.

Let $A$ be a finite Abelian group and $E/S$ a generalized elliptic
curve. A {\em weak $A$-structure} on $E$ is a group homomorphism
$\phi:A\to E^{\rm sm}(S)$ such that the Cartier divisor
$\phi(A):=\sum_{a\in A} \phi(a)$ is a subgroup-scheme of $E^{\rm sm}$.
A weak $A$-structure $\phi$ is called an {\em $A$-structure} if
$\phi(A)$ meets every irreducible component of every geometric fiber
of $E/S$.  The following two examples are classical.  For $A=\ZZ/n$,
an $A$-structure is called a {\em $\Gamma_1(n)$-structure}. For
$A=\ZZ/n\times\ZZ/n$, an $A$-structure is called a {\em
  $\Gamma(n)$-structure}.

\subsubsection{$\Gamma_2(N)$-structures} \label{G2N}

We fix an integer $N>0$ and an odd prime number $p$ such that
$p\,||\,N$. We define the Abelian group
\[
     A := \ZZ/2N\times\ZZ/2.
\]
Similarly, we let $A':=\ZZ/2N'\times\ZZ/2$, where $N=pN'$, and let
$\tau:A\to A'$ be the natural projection.

\begin{defn} \label{G2Ndefn}
  A {\em $\Gamma_2(N)$-structure} on a generalized elliptic curve
  $E/S$ is an $A$-structure $\phi:A\to E(S)$, with $A$ as above.
  We say that $\phi$ is {\em \'etale} if the subscheme $\phi(A)\subset
  E$ is \'etale over $S$ (equivalently, the induced map $\phi_s:A\to
  E_s(K)$ is injective for all geometric points $s:\Spec k\to S$). We
  say that $\phi$ is {\em $p$-local} if $\Ker\phi\subset A$ has order
  $p$. 
\end{defn}

The following is obvious from \cite{KatzMazur}:

\begin{rem} \label{G2Nrem}\ 
  \begin{enumerate}
  \item
    If $2N$ is invertible on $S$ then $\phi$ is \'etale.
  \item
    If $\phi$ is $p$-local then $S$ is an $\FF_p$-scheme. Moreover,
    $\phi=\phi'\circ\tau$, where $\phi':A'\to E(S)$ is an \'etale
    $\Gamma_2(N')$-structure.
  \end{enumerate}
\end{rem}

\begin{prop} \label{G2Nprop1}
  Let $G$ be a dihedral group of order $2N$ and $f:E\to\PP^1_K$ a
  $G$-cover, branched in $4$ points with ramification index $2$.
  After a finite extension of $K$, there exists a map $\phi:A\to
  E(K)$, whose image $\phi(A)\subset Y$ is the set of ramification
  points of $f$, such that $(E,\phi)$ is an elliptic curve with
  an \'etale $\Gamma_2(N)$-structure. 
\end{prop}

\proof The cover $f$ factors through an \'etale $N$-cyclic cover
$\pi:E\to E'$, and $E$ and $E'$ are smooth projective curves of genus
$1$. We may assume that all ramification points of $f$ are
$K$-rational. Let us choose one ramification point $0\in E(K)$, and
set $0':=\pi(0)$. Now we can regard $\pi$ is an $N$-cyclic isogeny
between elliptic curves. Moreover, the branch points of $f$ are
precisely the points of $E$ lying above the $2$-torsion points of
$E'$. Therefore, the set of branch points is a subgroup of $E[2N](K)$,
of order $4N$, and contains a point of order $2N$. It is clear that
this subgroup is isomorphic to $A$.  \Endproof

Consider the exact sequence
\begin{equation} \label{Aseq}
    0 \To \ZZ/N \To A \To \ZZ/2\times\ZZ/2 \To 0
\end{equation}
of Abelian groups, where $a\in\ZZ/N$ is send to $(2a,0)\in A$ and
$(a,b)\in A$ is send to $(\bar{a},\bar{b})$. Let $(E,\phi)$ be a
generalized elliptic curve with $\Gamma_2(N)$-structure. The image
$\phi(\ZZ/N)\subset E$ is a subgroup scheme of $E\sm$, finite and flat
over $S$. Since $E\sm$ acts on $E$, we can form the quotient scheme
$E':=E/\phi(\ZZ/N)$. One checks fiber by fiber that $E'/S$ is again a
generalized elliptic curve. Via the exact sequence \zgl{Aseq}, $E'$ is
endowed with a $\Gamma(2)$-structure $\phib:\ZZ/2\times\ZZ/2\to
E'(S)$. Let $[-1]:E'\iso E'$ be the canonical involution (see
\cite{DelRap}, Chapitre II) and $X_0:=E'/\gen{[-1]}$ the quotient. We
write $f_\phi:E\to X_0$ for the natural map. We choose a bijection
$\alpha:\{1,2,3,4\}\cong\ZZ/2\times\ZZ/2$ and let $x_j\in X_0(S)$ be
the image of $\phib(\alpha(j))$, for $j=1,\ldots,4$.

\begin{prop} \label{G2Nprop2}
  Suppose $2$ is invertible on $S$. 
  \begin{enumerate}
  \item
    The curve $(X_0;x_j)$ is stably marked, of genus $0$.
  \item 
    If $\phi$ is \'etale, then $f_\phi:E\to X_0$ is an admissible
    cover, ramified of order $2$ along the sections $x_j$. There is a
    natural $G$-action on $E$, where $G$ is dihedral of order $2N$,
    such that $X_0=E/G$. If in addition $E/S$ is smooth then $f_\phi$
    is a tame $G$-cover.
  \end{enumerate}
\end{prop}

\proof It suffices to prove the proposition in the case $S=\Spec k$,
where $k$ is an algebraically closed field.  Carrying a
$\Gamma(2)$-structure, $E'$ is either smooth or a $2$-gon.  Since
$[-1]$ is the identity on $E'[2]\cong\ZZ/2\times\ZZ/2$, the points
$x_1,\ldots,x_4\in X_0$ are distinct and smooth. If $E$ is smooth then
$X_0\cong\PP^1_k$. Otherwise, $[-1]$ restricts to an involution on
each component of $E'$ and interchanges the two singular point.  In
this case, $X_0$ is the union of two projective lines meeting
transversally in one point, and each component of $X_0$ contains two
of the points $x_1,\ldots,x_4$.  This proves (i).

We have seen that $E'\to X_0$ is ramified at $x_1,\ldots,x_4$ of order
$2$ and \'etale everywhere else.  Assume that $E/S$ is smooth and
$\phi(A)\subset E$ is \'etale. Then $\pi:E\to E'$ is an \'etale
$N$-cyclic isogeny. It follows that $f_\phi$ is a dihedral Galois
cover, ramified at $x_1,\ldots,x_4$ of order $2$. In case $E$ is
singular, the map $f_\phi$ may be ramified in the singular points.
However, using the description of the group law on a N\'eron polygon
given in \cite{DelRap}, Chapitre II, one checks that $f_\phi$ is
admissible.  \Endproof

\subsubsection{Deformation} \label{G2Ndef}

Let $k$ be an algebraically closed field of characteristic $p$. We fix
a generalized elliptic curve $E$ over $k$ and a
$\Gamma_2(N)$-structure $\phi:A\to E(k)$. We assume that $\phi$ is
$p$-local. Let $\Def(E,\phi)$ denote the deformation functor
classifying isomorphism classes of deformations $(E_R,\phi_R)$ of
$(E,\phi)$ over complete local $W(k)$-algebras $R$ with residue field
$k$. Let $\Def(E,\phi)\loc$ be the subfunctor corresponding to
deformations $(E_R,\phi_R)$ where $\phi_R$ is $p$-local.

\begin{prop} \label{G2Ndefprop1}
  The functor $\Def(E,\phi)$ has a universal deformation ring
  $R_\phi$. The ring $R_\phi$ is regular of dimension $2$; there
  exists a regular sequence $(t,\pi)$ with the following properties.
  \begin{enumerate}
  \item
    The ring $R_\phi\loc:=R_\phi/(\pi)\cong k[[t]]$ is the universal
    deformation ring for $\Def(E,\phi)\loc$.
  \item
    If $E$ is ordinary, then $p=u\pi^{p-1}$, where $u\in
    R_\phi^\times$.
  \item
    If $E$ is supersingular, then $p=u\pi^{p-1}$, and $(u,\pi)$ is
    another regular sequence for $R_\phi$. 
  \end{enumerate}
\end{prop}

\proof Let us first assume that $E$ is smooth. We write
$\Def(E,\phi|_{\ZZ/p})$ for the functor classifying deformations of
$E$ together with the $\Gamma_1(p)$-structure $\phi|_{\ZZ/p}$. Since
$A\cong A'\times\ZZ/p$ and the order of $A'$ is prime-to-$p$, the
morphism $\Def(E,\phi)\to\Def(E,\phi|_{\ZZ/p})$ that sends
$(E_R,\phi_R)$ to $(E_R,\phi_R|_{\ZZ/p})$ is an equivalence.
Therefore, if $E$ is smooth, the proposition is a direct consequence
of \cite{KatzMazur}, Section 13.5. Namely, the universal deformation
ring $R_\phi$ can be identified with the strict complete local ring of
the moduli stack $\M(\Gamma_1(p))$ at the point corresponding to
$(E,\phi|_{\ZZ/p})$.

It is clear from \cite{Edix} how to extend this to the general case.
Actually, since $p\,||\,N$, the situation here is somewhat easier than
in \cite{Edix}. Since $\phi$ is $p$-local, the number of components of
$E$ is prime-to-$p$ (more precisely, if $E$ is singular, the number of
components is $2n$, where $n|N'$). By \cite{DelRap}, Th\'eor\`eme
III.1.2, the generalized elliptic curve $E$ admits a universal
deformation $E_{R_0}$ over $R_0=W(k)[[t]]$. As in the smooth case, the
morphism $\Def(E,\phi)\to\Def(E,\phi|_{\ZZ/p})$ is an equivalence,
if we regard $\phi|_{\ZZ_p}$ as a weak $\ZZ/p$-structure. It
follows from \cite{KatzMazur} that $\Def(E,\phi)$ admits a universal
deformation $(E\univ,\phi\univ)$ over a ring $R_\phi$.  In fact,
$\Spec R_\phi$ is a closed subscheme of $E_{R_0}[p]^\times$ and
$E\univ=E_{R_0}\otimes_{R_0}R_\phi$. We regard $t$ as an element of
$R_\phi$ via the natural morphism $R_0=W(k)[[t]]\to R_\phi$. Let us
choose a formal parameter $T$ of $E\univ$ along the $0$-section (see
\cite{KatzMazur}, Section 2.2.3). The point $P\univ:=\phi\univ(2N',0)$
is a {\em point of exact order $p$} on $E\univ$. Since $\phi$ is
$p$-local, $P\univ|_E=\phi(2N',0)=0$. Hence we may regard
$\pi:=T(P\univ)$ as an element of $R_\phi$. We claim that $(t,\pi)$ is
a regular sequence for $R_\phi$ such that (i), (ii) and (iii) hold. We
have already mentioned that, if $E$ is smooth, this is proved in
\cite{KatzMazur}. If $E$ is singular, the situation is essentially the
same as for $E$ smooth and ordinary, see \cite{Edix}. Namely, we may
choose $T$ such that $E_{R_0}[p]^\times=\Spec R_0[\,T\mid
\Phi_p(1+T)=0\,]$, where $\Phi_p(X)=(X^p-1)/(X-1)$.  Therefore,
$R_\phi=W(k)[\zeta_p][[t]]$ and $\pi=\zeta_p-1$. This completes the
proof of the proposition.  \Endproof

Let $(E_R,\phi_R)$ be a deformation of $(E,\phi)$. Following Section
\ref{G2N} we associate to every deformation $(E_R,\phi_R)$ of
$(E,\phi)$ a finite map $f_{\phi_R}:E_R\to X_{0,R}$ and sections
$x_{1,R},\ldots,x_{4,R}:\Spec R\to X_{0,R}$ such that
$(X_{0,R};x_{j,R})$ is a stably marked curve of genus $0$. This gives
rise to a morphism
\begin{equation} \label{G2Ndefeq1}
    \Def(E,\phi) \To \Def(X_0;x_j)
\end{equation}
of deformation functors and hence to a $W(k)$-algebra morphism
$R_{X_0}\to R_\phi$. 

\begin{prop} \label{G2Ndefprop2}
  The ring $R_\phi$ is a finite and flat extension of $R_{X_0}\cong
  W(k)[[w]]$. Modulo $\pi$, we obtain a finite extension 
  \[
      R_{X_0}\otimes_{W(k)}k\cong k[[w]] \To R_\phi\loc\cong k[[t]]
  \]
  with inseparability degree $p$. Its separable part $k[[w]]\to
  (R_\phi\loc)^p=k[[t^p]]$ is tamely ramified of degree $d$. Here
  \[
       d \;=\; \left\{
         \begin{array}{cl}
           1    & \qquad \text{\rm if $E$ is smooth,} \\
           N'/n & \qquad \text{\rm if $E$ is a $2n$-gon.} \\
         \end{array}\right.
  \]
\end{prop} 

\proof Let $E'':=E/\phi(\ZZ/p)$ be the quotient by the subgroup scheme
$\phi(\ZZ/p)\subset E\sm$ and $\phi'':A'\to E''(k)$ the induced
$\Gamma_2(N')$-structure. Clearly, the map $f_\phi:E\to X_0$ induced
by $\phi$ factors through the projection $E\to E''$. The resulting map
$f_{\phi''}:E''\to X_0$ is the map induced by $\phi''$. We see that
the morphism \zgl{G2Ndefeq1} can be written as the composition of two
morphisms, as follows:
\begin{equation} \label{G2Ndefeq2}
   \Def(E,\phi) \To \Def(E'',\phi'') \To \Def(X_0;x_j).
\end{equation}
Let $R_{\phi''}$ be the universal deformation ring of
$\Def(E'',\phi'')$.  By Proposition \ref{G2Nprop2} (ii), the map
$f_{\phi''}:E''\to X_0$ is admissible, and is ``Galois'' with dihedral
Galois group of order $2N'$. If $E$ is singular, then $f_{\phi''}$ is
ramified of order $d$ over the unique singular point of $X_0$, where
$d$ is as in the statement of the proposition. It is not hard to see
that any deformation of the admissible cover $f_{\phi''}$ together
with the group action corresponds to a unique deformation of
$(E'',\phi'')$.  Therefore, it follows from \cite{diss} that
$R_{\phi''}=R_{X_0}[\,z\mid z^d=w\,]\cong W(k)[[z]]$. Here we identify
$R_{X_0}$ with $W(k)[[w]]$, such that, if $X_0$ is singular, $w$ is
the deformation parameter of the singular point. We are
reduced to showing that $R_{\phi''}\to R_\phi$ is finite and flat, and
purely inseparable of degree $p$ modulo $\pi$.

Let $B:=R_{\phi''}\otimes_{W(k)}k\cong k[[z]]$ and $\Bt:=B^{1/p}\cong
k[[z^{1/p}]]$. Let $(E''_{\Bt},\phi''_{\Bt})$ be the deformation of
$(E'',\phi'')$ corresponding to the natural morphism
$R_{\phi''}\to\Bt$. We can define a generalized elliptic curve
$E_{\Bt}$ over $\Bt$ such that $E''_{\Bt}=E_{\Bt}^{(p)}$ is the $p$th
power Frobenius twist of $E_{\Bt}$. The relative Frobenius
$F:E_{\Bt}\to E''_{\Bt}$ is a $p$-cyclic, purely inseparable
``isogeny'' whose kernel is generated by the $0$-section of $E_{\Bt}$
(which is a point of exact order $p$). Moreover, there exists a unique
$p$-local $\Gamma_2(N)$-structure $\phi_{\Bt}:A\to E_{\Bt}(\Bt)$ such
that $\phi''_{\Bt}=F\circ\phi_{\Bt}$. Conversely, let $(E_R,\phi_R)$
be any $p$-local deformation of $(E,\phi)$. By \cite{KatzMazur}, we
can canonically identify $E_R'':=E_R/\phi_R(\ZZ/p)$ with $E_R^{(p)}$
and the quotient map $E_R\to E_R''$ with the Frobenius $F:E_R\to
E_R^{(p)}$.  It follows that $\Bt=R_\phi\loc$. Hence $R_{\phi''}\to
R_{\phi}$ is purely inseparable of degree $p$ modulo $\pi$. Nakayama's
Lemma shows that $R_\phi$ is finite over $R_{\phi''}$. Since a finite
morphism between two regular local rings of the same dimension is
automatically flat (see \cite{AltKlei}, V.3.8), the proposition is
proved.  \Endproof

\subsubsection{Stabilization} \label{stab}

Let $(E_R,\phi_R)$ be a deformation of $(E,\phi)$. A {\em
  stabilization} of $(E_R,\phi_R)$ is a morphism $q_R:Z_R\to E_R$
between semistable $R$-curves together with a map $\psi_R:A\to Z_R(R)$
such that (i) $(Z_R,\psi_R(A))$ is a stably marked curve, (ii)
$\phi_R=q_R\circ\psi_R$, and (iii) $q_R$ is the contraction of the
marked semistable curve $(Z_R,\psi_R(A'))$ (see \cite{Knudsen83}).

\begin{lem} \label{stablem}
  For every deformation $(E_R,\phi_R)$, there exists a stabilization
  $(Z_R,\psi_R)$.  Assume that there exists a dense open subset
  $U\subset\Spec R$ such that $\phi_U$ is \'etale. Then $(Z_R,\psi_R)$
  is unique up to unique isomorphism.
\end{lem}

\proof 
If $\phi_U$ is \'etale, then $(E_U,\phi_U(A))$ is stably marked, so
necessarily $Z_U=E_U$. Therefore, the uniqueness follows from the
fact that the moduli stack of stably marked curves is separated, see
\cite{Knudsen83}. 

The notion of stabilization is certainly compatible with base change
$R\to R'$. Hence it suffices to prove the existence of stabilization
for the universal deformation $(E\univ,\phi\univ)$. Recall that
$\pi=T(\phi\univ(2N',0))\in R_\phi$, where $T$ is a formal parameter
of $E\univ$ along the $0$-section. Relative to $T$, the formal group
of $E\univ\sm$ is given by a power series
$\Phi(T_1,T_2)=T_1+T_2+\ldots\in R[[T_1,T_2]]$. Since $\phi\univ$ is a
group homomorphism, we have
\begin{equation} \label{grouplaw}
    T(\phi\univ(2N'm,0)) \;\equiv\; m\pi \pmod{\pi^2},\qquad m\in\ZZ/p.
\end{equation}
Let $q\univ:Z\univ\to E\univ$ be the blowup of $E\univ$ along the
closed subscheme $\phi\univ(A')\cap(\pi)\subset E\univ$. Since
$\phi\univ(A')$ consists of $|A'|=4N'$ pairwise disjoint sections
$\Spec R_\phi\to E\univ\sm$, the blowup $Z\univ$ is a semistable curve
over $R_\phi$.  Denote its special fiber by $Z$. Then
\[
       Z = E \cup \bigcup_{b\in A'} Z_b,
\]
where $Z_b\cong\PP^1_k$ is connected to $E$ in the point
$z_b:=\phi(b)$. For $b\in A'$, the translate $T_b:=T-\phi\univ(b)$ is
a formal parameter of $E$ along the section $\phi\univ(b)$. By the
definition of $Z\univ$, $\tilde{T}_b:=T_b/\pi$ is a regular function
in a neighborhood of $Z_b-\{z_b\}\subset Z\univ$ and defines an
isomorphism $Z_b\cong\PP^1_k$ mapping $z_b$ to $\infty$. For $a\in A$,
let $\psi\univ(a)$ be the closure of $\phi\univ(a)|_K\in
E\univ(K)=Z\univ(K)$ in $Z\univ$. This defines a map $\psi\univ:A\to
Z\univ(R_\phi)$ such that $\phi\univ=q\univ\circ\psi\univ$. Write
$a=b+c$, with $b\in A'$ and $c\in\ZZ/p$.  Using \zgl{grouplaw} one
finds that $\psi\univ(a)$ meets the special fiber in the point $c\in
Z_b(\FF_p)\cong\FF_p\cup\{\infty\}$.  By construction,
$(Z\univ,\psi\univ)$ is a stabilization of $(E\univ,\phi\univ)$. This
proves the lemma.  \Endproof

As in the proof of the lemma, let $Z\univ$ be the stabilization of the
universal deformation $(E\univ,\phi\univ)$. Let $f_{0,\univ}:Z\univ\to
X_{0,\suniv}$ be the composition of $q\univ$ with the map
$f_{\phi\univ}:E\univ\to X_{0,\suniv}$ induced by $\phi\univ$.  Let
$K$ be the fraction field of $R_\phi$. Since $\phi_K$ is \'etale,
$Z_K=E_K$. By Proposition \ref{G2Nprop1}, $f_{0,K}:Z_K\to X_{0,K}$ is
a $G$-cover, branched at $4$ points of order $2$, where $G$ is
dihedral of order $2N$. In other words, $Z_K$ is a $K$-object of the
Hurwitz stack $\H_4\inn(G)$. By the uniqueness of stabilization, the
action of $G$ extends to $Z\univ$. Therefore, $Z\univ$, together with
the mark $\psi\univ(A)$, the action of $G$ and the map $f_{0,\suniv}$,
is an $R_\phi$-object of the complete Hurwitz stack $\Hb_4\inn(G)$.
In particular, the special fiber $Z$ of $Z\univ$ is a $k$-object of
$\Hb_4\inn(G)\bad$.  Via pullback, we obtain a morphism
\begin{equation} \label{stabeq1}
  \Def(E,\phi) \To \Def(Z)
\end{equation}
of deformation functors, compatible with the morphisms
$\Def(E,\phi)\to\Def(X_0;x_j)$ and $\Def(Z)\to\Def(X_0;x_j)$.

\begin{prop} \label{stabprop}
  The morphism \zgl{stabeq1} is an isomorphism; it induces an
  isomorphism between $\Def(E,\phi)\loc$ and $\Def(Z)\bad$. 
\end{prop}

\proof The first statement is equivalent to the assertion that
$Z\univ$ is the universal deformation of $Z$. By construction, the
maximal subset $U\subset\Spec R_\phi$ such that $\phi_U$ is \'etale is
precisely the maximal subset such that $Z_U\to X_{0,U}$ is an
admissible cover. Therefore, the first statement of the proposition
implies the second.

Let $Z_R$ be a deformation of $Z$. We denote by $\psi:A\to Z(k)$ the
restriction of $\psi\univ$ to $Z$. Clearly, $\psi$ lifts uniquely to a
map $\psi_R:A\to Z_R(R)$ such that $\psi_R(A)$ is the mark of the
stably marked curve $Z_R$.  Let $q_R:Z_R\to E_R$ be the contraction of
the marked semistable curve $(Z_R,\psi_R(A'))$, and let
$\phi_R:=q_R\circ\psi_R$. We claim that the assignment $Z_R\mapsto
(E_R,\phi_R)$ defines a morphism $\Def(Z)\to\Def(E,\phi)$, which is
the inverse of \zgl{stabeq1}. This is clear from the construction and
the following lemma.

\begin{lem} \label{stablem2}
  There exists a unique morphism $+_R:E_R\sm\times E_R\to E_R$ such that
  $(E_R,+_R,\phi_R(0))$ is a generalized elliptic curve and $\phi_R$
  is a $\Gamma_2(N)$-structure. 
\end{lem}

\proof By construction, $E_R$ is a semistable $R$-curve of genus $1$.
We claim that every geometric fiber $E_s$ of $E_R$ is either smooth or
a N\'eron polygon, and that $\phi_R(A)$ meets every irreducible
component of $E_s$. In fact, this is an open condition on $\Spec R$,
and it is true for the special fiber $E$. 

Let us first prove the uniqueness of $+_R$. Suppose a morphism $+_R$
satisfying the conditions of the lemma exists. Since $\phi_R:A\to
E_R\sm(R)$ is a group homomorphism (with respect to the group law
induced by $+_R$), it induces an action of $A$ on $E_R$ via
translation. Since $\phi_R$ is a $\Gamma_2(N)$-structure, $A$ acts
transitively on the set of irreducible components of every geometric
fiber. By definition of this action, $\phi_R$ is equivariant (here $A$
acts on itself by translation). But $(E_R,\psi_R(A'))$ is stably
marked, so there can exists at most one $A$-action on $E_R$ with this
property. Therefore, we can apply \cite{DelRap}, Th\'eor\`eme II.3.2,
to conclude that there exists at most one morphism $+_R$ with the
claimed properties.

The assignment $Z_R\mapsto (E_R,\phi_R)$ is clearly compatible with
base change $R\to R'$. Hence it suffices the prove the existence of
$+_R$ in the case $R=R_Z$, i.e.\ when $Z_R$ is the universal
deformation of $Z$.  Let $U\subset\Spec R$ be the maximal subset such
that $Z_U$ is a $G$-cover. By the construction of the complete Hurwitz
stack $\Hb_4\inn(G)$, $U$ is open and dense. Moreover, $Z_U=E_U$ is a
smooth curve of genus $1$. By \cite{DelRap}, Proposition II.2.7, there
exists a unique structure of elliptic curve on $(E_U,\phi_U(0))$.  As
in the proof of Proposition \ref{G2Nprop1} one shows that $\phi_U$ is
a $\Gamma_2(N)$-structure. In particular, $A$ acts on $E_U$ such that
$\phi_U$ is equivariant. Since $(E_R,\phi_R(A'))$ is stably marked and
$U\subset\Spec R$ is dense, this action extends uniquely to $E_R$, and
$\phi_R$ is equivariant. By
\cite{DelRap}, Proposition II.2.7, the induced action of $A$ on
${\rm Pic}^0E_R/R$ is trivial. Therefore, we can apply
\cite{DelRap}, Th\'eor\`eme II.3.2, to show that there exists a
structure of generalized elliptic curve on $E_R$ such that the action
of $a\in A$ on $E_R$ is given by translation with $\phi_R(a)$. It
remains to show that $\phi_R$ is a weak $A$-structure. But this is a
closed condition on $\Spec R$ (see \cite{KatzMazur}) and it is true on
$U\subset \Spec R$. Hence it is true on $\Spec R$. This concludes the
proof of Lemma \ref{stablem2} and Proposition \ref{stabprop}.
\Endproof

By Proposition \ref{stabprop}, we can identify the universal
deformation rings $R_\phi$ and $R_Z$, and regard $Z\univ$ as the
universal deformation of $Z$. In view of Proposition \ref{G2Ndefprop1}
and Proposition \ref{G2Ndefprop2}, we obtain:

\begin{cor} \label{stabcor}
  Theorem \ref{precisethm} is true for $Y=Z$.
\end{cor}

\begin{rem} \label{stabrem}
  Let $X_0(N)_{\ZZ_{(p)}}$ and $X_1(N)_{\ZZ_{(p)}}$ be the arithmetic
  models over $\ZZ_{(p)}$ of the modular curves $X_0(N)$ and $X_1(N)$,
  as defined in \cite{DelRap} and \cite{KatzMazur}. The results of
  this section imply that $X_0(N)_{\ZZ_{(p)}}\cong\hb_{[4]}\ab(\Cl)$
  and $X_1(N)_{\ZZ_{(p)}}\cong\hb_{[4]}\inn(\Cl)$, where
  \begin{equation} \label{stabremeq} 
      \Cl \;:=\; \left\{
        \begin{array}{cl} 
           (2A,2A,2A,2A) & \quad\text{\rm if $N$ is odd,} \\
           (2A,2A,2B,2B) & \quad\text{\rm if $N$ is even,}\\
        \end{array}\right.
  \end{equation}
  and $2A$, $2B$ denote the conjugacy classes of the ``reflections'' in
  a dihedral group of order $2N$. Let $X_2(N)$ be the coarse moduli
  space for generalized elliptic curves with
  $\Gamma_2(N)$-structure. One can show that
  $X_2(N)_{\ZZ_{(p)}}\cong\hb_4\inn(\tilde{\Cl})$, where
  $\tilde{\Cl}=(2\tilde{A},2\tilde{A},2\tilde{B},2\tilde{B})$ is the
  tuple of conjugacy classes in a dihedral group of
  order $4N$, as in \zgl{stabremeq}. 
\end{rem}

\subsection{Proof of the Reduction Theorem} \label{proofs}

We are now ready to complete the proof of the Reduction Theorem. The
main argument is given in Section \ref{deformation}, where we compare
the deformation theory of $Y$ to the deformation theory of the special
fiber $Z$ of the associated auxiliary cover. Since we have modular
reduction, Theorem \ref{precisethm} follows from the results of
Section \ref{dihedral}. In Section \ref{monoproof} and Section
\ref{innermonoproof} we prove Proposition \ref{monoprop} and
Proposition \ref{innermonoprop}, using the same method: we first
reduce the statements to the dihedral case and then use the results of
Section \ref{dihedral}.

\subsubsection{The auxiliary cover} \label{aux}

Let $f_0:Y\To X_0$ be as in Section \ref{precise1}. As we have seen in
the paragraph following Definition \ref{hurwdef1}, there exists a
complete discrete valuation ring $R$ with residue field $k$ and
quotient field $K$ of characteristic $0$ such that $f_0:Y\to X_0$ is
the reduction of a $G$-cover $f_K:Y_K\to\PP^1_K$. More precisely, the
$G$-cover $f_K$ has a stable model $f_{0,R}:Y_R\to X_{0,R}$ over $R$,
with special fiber $f_0$. We denote by $f_R:Y_R\to X_R$ the quotient
model of $f_K$ and by $f:Y\to X$ its special fiber (see Section
\ref{reduction}). Note that $f_0$ factors through $f$ and that $f$
does not depend on the choice of the lift $f_K$, by Proposition
\ref{quotprop}.

We are assuming that Condition \ref{cond1} holds. Therefore, it
follows from Proposition \ref{modularprop} that the $G$-cover $f_K$
has {\em modular reduction} of level $N$, where $N$ is some integer
diving $|G|$ and divisible by $p$ (see Definition \ref{modulardef}).
Let $g_K:Z_K\to\PP^1_K$ be the auxiliary cover associated to
$f_K$. The cover $g_K:Z_K\to\PP^1_K$ is a $\Delta$-cover, where
$\Delta\subset G$ is a dihedral group of order $2N$; it
has the same branch locus as $f_K$, but with ramification of order
$2$. Let $\Cl\aux$ be the inertia type of $g_K$, and let
$\Hb\aux:=\Hb_4\inn(\Cl\aux)_{\ZZ_{(p)}}$ be the complete Hurwitz
stack over $\ZZ_{(p)}$ for $\Delta$-covers with inertia type
$\Cl\aux$.  By definition, $g_K$ is a $K$-object of $\Hb\aux$. Let
$g_R:Z_R\to X_R$ and $g_{0,R}:Z_R\to X_{0,R}$ be the quotient and the
stable model of $g_K$. By definition, $g_{0,R}$ is an $R$-object of
$\Hb\aux$ extending $g_K$; its special fiber $g_0:Z\to X_0$ is a
$k$-object of $\Hb\aux\bad$. In the sequel, we will denote this object
simply by $Z$. 

According to Section \ref{G2N}, there exists a map $\phi_K:A\to
Z_K(K)$ such that $(Z_K,\phi_K)$ is an elliptic curve with
$\Gamma_2(N)$-structure, and the cover $g_K:Z_K\to X_K$ is induced by
$\phi_K$. We extend $\phi_K$ to a map $\psi_R:A\to Z_R(R)$; let
$q_R:Z_R\to E_R$ be the contraction of $(Z_R,\psi_R(A'))$ and let
$\phi_R:=q_R\circ\psi_R$. By Lemma \ref{stablem2}, $(E_R,\phi_R)$ is a
generalized elliptic curve with $\Gamma_2(N)$-structure. Let
$(E,\phi)$ be the special fiber of $(E_R,\phi_R)$. Clearly, we are in
the situation of Section \ref{stab}. In particular, Theorem
\ref{precisethm} is true for $Y=Z$, by Corollary \ref{stabcor}.

\subsubsection{Proof of Theorem \ref{precisethm}} 
\label{deformation}

By the construction of the auxiliary cover, there exists an \'etale
map $U^{(1)}\to X$, covering $X_0$, and a $G$-equivariant isomorphism
\begin{equation} \label{Indeq}
      Y\times_X U^{(1)} \;\cong\; \Ind_\Delta^G (Z\times_X U^{(1)})
\end{equation}
of $X$-schemes. Over the open subset $U^{(2)}:=X-X_0$, the
map $g$ is tamely ramified along the mark $C\subset X$.  Let
$X\univ:=Z\univ/\Delta$. Since $R_Z$ is regular, the natural morphism
$g\univ:Z\univ\to X\univ$ is a quotient model of $Z\univ$, by
Proposition \ref{quotprop}. In particular, $X\univ$ is a semistable
curve and carries a natural mark $C\univ\subset X\univ$.

\begin{constr} \label{defconstr}
  Let $R\in\C_k$ be Artinian, and let $Z_R$ be a deformation of $Z$
  over $R$.  There exist a unique morphism $R_Z\to R$ such that
  $Z_R=Z\univ\otimes_{R_Z}\!R$. Let $X_R:=X\univ\otimes_{R_Z}\!R$.
  Then $Z_R\to X_R$ is a quotient model of $Z_R$. Note that the pair
  $(U^{(1)},U^{(2)})$, where $U^{(1)}\to X$ and $U^{(2)}\to X$ are
  as above, is an \'etale covering of $X$. For $i=1,2$, let
  $Y^{(i)}:=Y\times_X U^{(i)}$ and $Z^{(i)}:=Z\times_X U^{(i)}$.  For
  $i,j=1,2$, let $U^{(i,j)}:=U^{(i)}\times_XU^{(j)}$,
  $Y^{(i,j)}:=Y\times_XU^{(i,j)}$ and $Z^{(i,j)}:=Z\times_XU^{(i,j)}$.
  Since $R$ is Artinian, the \'etale covering $(U^{(1)},U^{(2)})$ of
  $X$ extends uniquely to an \'etale covering $(U_R^{(1)},U_R^{(2)})$
  of $X_R$. Actually, $U_R^{(2)}=X_R-X_0$. Over $U^{(2)}$, the finite
  map $f:Y\to X$ is tamely ramified along $C$. By a theorem of
  Grothendieck and Murre, there exists a unique extension of
  $Y^{(2)}\to U^{(2)}$ to a finite morphism $Y_R^{(2)}\to U_R^{(2)}$
  which is tamely ramified along $C_R\subset X_R$. Define
\begin{equation} \label{Indeq2}
    Y_R^{(1)} := \Ind_\Delta^G (Z_R\times_{X_R} U_R^{(1)}).
\end{equation}
By \zgl{Indeq}, $Y_R^{(1)}\otimes_R k=Y^{(1)}$. We claim that there
exist $G$-equivariant isomorphisms
\begin{equation} \label{patchingeq}
    \alpha_R^{(i,j)}:Y_R^{(i)}\times_{U_R^{(i)}}U_R^{(i,j)}
    \liso Y_R^{(j)}\times_{U_R^{(j)}}U_R^{(i,j)},\quad i,j=1,2
\end{equation}
of $U_R^{(i,j)}$-schemes which extend the identity on $Y^{(i,j)}$.
For $(i,j)\not=(1,1)$, this follows again from Grothendieck--Murre.
For $(i,j)=(1,1)$, we obtain \zgl{patchingeq} via the canonical
identification
\begin{equation}
    Y_R^{(1)}\times_{U_R^{(1)}}U_R^{(1,1)}\;\cong\;
    \Ind_{H}^G (Z_R \times_{X_R} U_R^{(1,1)}) \;\cong\;
    Y_R^{(1)}\times_{U_R^{(1)}}U_R^{(1,1)}
\end{equation}
(on the left hand side we use the first projection $U_R^{(1,1)}\to
U_R^{(1)}$, on the right hand side we use the second projection). It
is clear that the isomorphisms \zgl{patchingeq} verify the obvious
cocycle condition. Therefore, there exist a finite morphism
$f_R:Y_R\to X_R$ such that $Y_R^{(i)}=Y_R\times_{X_R}U_R^{(i)}$, for
$i=1,2$. 
\end{constr}

Construction \ref{defconstr} associates to any deformation
$Z_R\in\Def(Z)(R)$ over an Artinian ring $R\in\C_k$ a finite map
$f_R:Y_R\to X_R$ which extends $f:Y\to X$. Moreover, the $G$-action on
$Y$ and the mark $D\subset Y$ extend to $Y_R$. Since $X_R$ is
projective over $R$, we can apply Grothendieck's Existence Theorem and
extend this construction to the case of an arbitrary ring $R\in\Cd_k$.
We claim that the constructed curve $Y_R$, together with the
$G$-action, the natural morphism $f_{0,R}:Y_R\to X_{0,R}$ and the mark
$D_R\subset Y_R$, is an object of $\Hb$. It suffices to prove this in
the case $R=R_Z$, $Z_R=Z\univ$. So let $Y_{R_Z}\to X\univ$ be the map
associated to $Z\univ$ by Construction \ref{defconstr}. Since the
generic fiber of $g\univ:Z\univ\to X\univ$ is a $\Delta$-cover, the
generic fiber of the map $Y_{R_Z}\to X\univ$ is a $G$-cover, i.e.\ an
object of $\H$. By construction, $Y_{R_Z}\to X\univ$ is an object of
$\SB\rr\inn(G)$.  Therefore, $Y_{R_Z}\to X\univ$ is an object of
$\Hb$. We have shown that Construction \ref{defconstr} defines a
morphism of functors
\begin{equation} \label{defmorph}
     \Def(Z) \To \Def(Y).
\end{equation}
To complete the proof of Theorem \ref{precisethm}, we have to show
that \zgl{defmorph} is an isomorphism. We need two lemmas.

\begin{lem} \label{deflem2}
  Let $R\in\C_k$ be an Artinian $k$-algebra. Let $Z_R$ be a
  deformation of $Z$ over $R$ and $Y_R$ its image under
  \zgl{defmorph}. If $Y_R\cong Y\otimes_k R$ is a trivial deformation,
  then $Z_R\cong Z\otimes_k R$ is trivial, too.
\end{lem}

\proof It suffices to show that the deformation $(E_R,\phi_R)$ of
$(E,\phi)$ corresponding to $Z_R$ is trivial. Since $Y_R\cong
Y\otimes_k R$, the closed embedding $E\inj Y$ extends to a closed
embedding $E\otimes_k R\inj Y_R$. By Construction \ref{defconstr} and
descent, we obtain a closed embedding $E\otimes_k R\inj Z_R$ extending
the closed embedding $E\inj Z$. Composition with the contraction
morphism $q_R:Z_R\to E_R$ yields a morphism $E\otimes_k R\to E_R$
which restricts to the identity on the special fiber. Since both
$E\otimes_k R$ and $E_R$ are flat and of finite type over $R$, it is
an isomorphism.  By construction, this isomorphism identifies $\phi_R$
with $\phi\otimes_k R$. This proves the lemma.  \Endproof

\begin{lem} \label{deflem1}
  Let $R\in\Cd_k$ be a normal domain, with fraction field of
  characteristic $0$. Let $Y_R$ be a deformation of $Y$ over $R$.
  Then there exists a deformation $Z_R$ of $Z$ over $R$ whose image
  under \zgl{defmorph} is isomorphic to $Y_R$.
\end{lem}

\proof 
Let $X_R':=Y_R/G$. By Proposition \ref{quotprop}, the natural
map $f_R:Y_R\to X_R'$ is a quotient model of $Y_R$. In particular,
$X_R'$ is a semistable curve with special fiber $X$. The image of
$D_R\subset Y_R$ is a mark $C_R'\subset X_R'$. There exists a maximal
open subset $U_R'\subset X_R'$, containing $C_R'$, over which $f_R$ is
tamely ramified along $C_R'$. Clearly, $U_R'$ contains the
generic fiber and $U_R'\otimes_R k=X-X_0$.

We claim that there exists a finite morphism $g_R:Z_R\to X_R'$
extending $g:Z\to X$, with $\Delta$ acting on $Z_R$, characterized by
the following two properties: (i) over $U_R'$, the map $g_R$ is tamely
ramified along $C_R'$, (ii) over an \'etale neighborhood of
$X_R'-U_R'$, $Y_R$ is isomorphic to $\Ind_\Delta^G(Z_R)$. In fact, one
can construct $g_R$ using the same method as in Construction
\ref{defconstr}. It follows that the generic fiber of $g_R$ is a
$\Delta$-cover, hence an object of $\H\aux$. Therefore,
$Z_R\in\Def(Z)(R)$. By construction, $g_R:Z_R\to X_R'$ is a quotient
model of $Z_R$. The uniqueness of the quotient model implies
$X_R'=X\univ\otimes_{R_Z}\!R$. A formal verification shows that $Y_R$
is the image of $Z_R$ under \zgl{defmorph}.  \Endproof

We are now going to complete the proof of Theorem \ref{precisethm}.
The morphism \zgl{defmorph} induces a homomorphism $R_Y\to R_Z$ of
local rings. We have to show that it is an isomorphism.  Let $\m_Y\lhd
R_Y$, $\m_Z\lhd R_Z$ be the maximal ideals.  Lemma \ref{deflem2}
implies that for every $N>0$ the $k$-module $R_Z/\m_Y R_Z$ is
generated by the images of $1$ and $(\m_Z)^N$. It follows that
$R_Z/\m_Y R_Z=k$. So $R_Y\to R_Z$ is surjective, by Nakayama's Lemma.

Let $\p\in\Spec R_Y$ be a generic point. The quotient $A:=R_Y/\p$ is a
complete local domain with residue field $k$ and fraction field $K$ of
characteristic $0$. The integral closure $\tilde{A}$ of $A$ in $K$ is
again a complete local domain with residue field $k$. So by Lemma
\ref{deflem1}, the morphism $R_Y\to\tilde{A}$ factors via $R_Y\to
R_Z$. Therefore, $I:=\Ker(R_Y\to R_Z)\subset\p$. It follows that $I$
is contained in the nilradical of $R_Y$. But $R_Y$ is reduced, so
$I=0$ and \zgl{defmorph} is an isomorphism. This completes the proof
of Theorem \ref{precisethm}.  \Endproof

\subsubsection{Proof of Proposition \ref{monoprop}} \label{monoproof}

We fix an element $\gamma\in\Gamma\ab$ and choose an automorphism
$\sigma:Y\iso Y\in\Aut_k\ab(Y)$ which induces $\gamma$. Such an
automorphism $\sigma$ is unique up to composition with an element of
$G$ (later in the proof we will give a ``canonical'' choice). Recall
that the generalized elliptic curve $E$ is a closed subscheme of $Y$.
Let $E_0$ be the identity component of $E$.  Since $E_0$ is an
irreducible component of $Y$ above $X_0$, we may assume
$\sigma|_{E_0}=\Id_{E_0}$, after composing $\sigma$ with an
appropriate element of $G$. Since $E$ is either irreducible or a
N\'eron polygon, it follows that $\sigma$ induces an automorphism
$\sigma|_E:E\liso E$ of the $k$-curve $E$ which normalizes the action
of $\Delta=D(E)$ and commutes with $f_0|_E:E\to X_0$. By the
construction of $Z$, $\sigma|_E$ extends uniquely to an automorphism
$\sigma_Z:Z\iso Z$ such that $\sigma_Z$ and $\sigma$ agree in an
\'etale neighborhood of $E$. Note that $\sigma_Z$ normalizes the
action of $\Delta$ and commutes with $g_0:Z\to X_0$. In other words,
$\sigma_Z\in\Aut_k\ab(Z)$. Therefore, $\sigma_Z$ lifts to a
$\gamma'$-semilinear automorphism $\sigma_{Z,\suniv}:Z\univ\iso
Z\univ$, for some $W(k)$-algebra automorphism $\gamma'$ of $R_Y=R_Z$.
Let $\hat{E}$ be the formal completion of $Z\univ$ along $E\subset
Z\univ$. It follows from Construction \ref{defconstr} that we can
identify $\hat{E}$ with the formal completion of $Y\univ$ along
$E\subset Y\univ$, and that
$\sigma_{Z,\suniv}|_{\hat{E}}=\sigma\univ|_{\hat{E}}$. In particular,
$\gamma'=\gamma$. Therefore, to prove Proposition \ref{monoprop}, we
may assume that $Y=Z$ and $G=\Delta$.

Recall that there exists a map $\psi\univ:A\to Z\univ(R_Z)$ such that
the following holds: (i) $\psi\univ(A)$ is the mark of the stably
marked curve $Z\univ$, (ii) $(Z\univ,\psi\univ)$ is the stabilization
of $(E\univ,\phi\univ)$ and (iii) $(E\univ,\phi\univ)$ is the
universal deformation of its special fiber $(E,\phi)$. Let $\psi:A\to
Z(k)$ be the restriction of $\psi\univ$ to $Z$ and let $Z_0$ be the
component of $Z$ which meets $E$ in $0\in E$. By the proof of Lemma
\ref{stablem}, we may identify $Z_0$ with $\PP^1_k$ such that
$\infty\in Z_0$ is the point where $Z_0$ meets $E$ and such that
$\psi(2N'a,0)=a\in\FF_p\subset Z_0(k)$. It is clear that $\sigma:Z\iso
Z$ restricts to an automorphism of $Z_0$ which fixes $\infty$ and
permutes the points $\psi(2N'a,0)$. Therefore, $\sigma$ acts on $Z_0$
as $z\mapsto cz+b$, with $b,c\in\FF_p$, $c\not=0$. Composing $\sigma$
by an element of the decomposition group $D(Z_0)\subset G$, we may
assume that $b=0$. This determines $\sigma$ uniquely. Moreover, 
\[
    \chi\bad(\gamma) := c^{-1}
\]
defines a homomorphism $\chi\bad:\Gamma\ab\to\FF_p^\times$.  Since
$\sigma\univ$ is an automorphism of $Z\univ$ as stably marked curve,
it descents to an automorphism $\sigmat\univ:E\univ\iso E\univ$ such
that $q\univ\circ\sigma\univ=\sigmat\univ\circ q\univ$, where
$q\univ:Z\univ\to E\univ$ is the projection morphism.  Let
$P\univ:=\phi\univ(2N',0)\subset E\univ$. Since $\phi\univ$ is a group
homomorphism, we have $a\cdot P\univ=\phi\univ(2N'a,0)$, for
$a\in\ZZ/p$. Using the definition of $\chi\bad$, we obtain
\begin{equation} \label{monopropeq1}
  \sigmat\univ(a\cdot P\univ) = \chi\bad(\gamma)^{-1}\cdot a\cdot P\univ, 
  \qquad a\in\ZZ/p.
\end{equation}
In particular, $\sigmat\univ$ fixes the $0$-section of
$E\univ$. Similar to the proof of Lemma
\ref{stablem2}, one shows that $\sigmat$ is a $\gamma$-semilinear
automorphism of the generalized elliptic curve $E\univ$, i.e.\ is
compatible with the ``group law'' on $E\univ$. In particular,
$\sigma|_E:E\iso E$ is a $k$-linear automorphism of the generalized
elliptic curve $E$.
 
Recall that $E\univ=E_{R_0}\otimes_{R_0}R_Z$, where $E_{R_0}$ is the
universal deformation of $E$, defined over $R_0=W(k)[[t]]\subset
R_Z$. It follows that $\sigmat\univ$ is induced by a $\gamma_0$-linear
automorphism $\sigma_{R_0}:E_{R_0}\iso E_{R_0}$, where
$\gamma_0:=\gamma|_{R_0}$ is the monodromy action of $\sigma|_E$ on
$R_0$. If $E$ is smooth then $\sigma|_E=\Id_E$. So in this case we
have $\gamma(t)=t$, as claimed in Proposition \ref{monoprop}. 
If $E$ is not smooth, it is isomorphic to the standard $2n$-gon, for
some $n|N'$, see \cite{DelRap}, Section II.1.1. So we identify the group
of components $E\sm/E_0\sm$ with $\ZZ/2n$ and the individual
components $E_i$, $i\in\ZZ/2n$, with $\PP^1_k$. According to
\cite{DelRap}, Proposition II.1.10, the automorphism $\sigma|_E$ is given
by the formula
\[
     \sigma|_E(x,i) = (\zeta^ix,i), \qquad x\in\PP^1,\;i\in\ZZ/2n,
\]
for some $2n$th root of unity $\zeta$. But $\sigma|_{E_i}:E_i\iso E_i$
lifts to an element of the decomposition group $D(E_i)\subset G$,
which is dihedral of order $2N/n$, by Proposition \ref{redbpts} (d).
We conclude that $\zeta$ is a $d$th root of unity, where $d$ is the
greatest common divisor of $2n$ and $N'/n$. We set
$\chi\adm(\gamma):=\zeta$. It is obvious that this defines a group
homomorphism $\chi\adm:\Gamma\ab\to\mu_d$. If we identify $E\univ$
with the Tate elliptic curve $\mathbb{G}_m^t/q^{\ZZ}$, where
$q=t^{2n}$, then we find $\gamma(t)=\zeta\cdot t$, see
\cite{DelRap}, Chapitre VII.

It remains to prove the formula
$\gamma(\pi)\equiv\chi\bad(\gamma)\cdot\pi\pmod{\pi^2}$. Recall that
$\pi=T(P\univ)$, where $T$ is a formal parameter of $E\univ$ along the
$0$-section. Actually, $T$ was chosen as a formal parameter of
$E_{R_0}$, see the proof of Proposition \ref{G2Ndefprop1}.  It follows
from the preceding discussion that we may assume
$\sigmat\univ^*(T)=T$. Using
\zgl{monopropeq1} and the formal group law on $E\univ$, we get
\begin{equation} \label{monopropeq2}
    \pi':= T(\sigmat\univ(P\univ)) \;\equiv\; 
         \chi\bad(\gamma)^{-1}\cdot\pi  \pmod{\pi^2}.
\end{equation}
By definition, we have
\begin{equation} \label{monopropeq3}
   T - \gamma(\pi') \;=\; \sigmat\univ^*(T-\pi') 
    \;=\; u\,(T-\pi),\qquad u\in R_Z[[T]]^\times.
\end{equation}
Comparing coefficients, we find $u\equiv 1\pmod{\pi}$, and therefore
$\gamma(\pi')\equiv\pi\pmod{\pi^2}$. With \zgl{monopropeq2}, we
conclude that $\gamma(\pi)\equiv\chi\bad(\gamma)\cdot\pi\pmod{\pi^2}$.
\Endproof

\subsubsection{Proof of Proposition \ref{innermonoprop}}
\label{innermonoproof}

We assume that Condition \ref{cond2} holds. Recall that the curve $X$
consists of five components $X_0,\ldots,X_4$, where $X_i$ contains the
specialization of the branch point $x_i$, for $i=1,\ldots,4$. Fix
$i\in\{1,\ldots,4\}$ and let $W_i$ be a component of $f^{-1}(X_i)$.
The restriction of $f:Y\to X$ to $W_i$ is a $D(W_i)$-Galois cover
$W_i\to X_i$, branched at $2$ points. Over $x_i$, the cover $W_i\to
X_i$ is tamely ramified, with inertia type $C_i$. Over the point where
$X_i$ meets $X_0$, we have wild ramification, of order $2p$. It
follows from Condition \ref{cond2} (f) that $D(W_i)=G$, i.e.\ 
$W_i=f^{-1}(X_i)\to X_i$ is a $G$-cover. Let $w\in W_i$ be a point
where $W_i$ meets $E\subset Y$.  The inertia group $I(w)$ is dihedral
of order $2p$.

Since $C_G=1$, by Condition \ref{cond2} (e), every element $\gamma\in
\Gamma$ is induced by a unique element $\sigma\in\Aut_k(Y)$. It is
clear that $\sigma$ fixes the component $W_i$ and that $\sigma(w)$ is
again a singular point of $Y$. Moreover, there exists an element
$\tau\in G$ such that $\sigma':=\tau^{-1}\circ\sigma\in\Aut_k\ab(Y)$
fixes $w$.  The element $\tau$ is unique up to composition with an
element of $I(w)$ and normalizes $I(w)$. For $h\in I(w)$, we have
$\sigma' h(\sigma')^{-1}=\tau^{-1}h\tau$. By Condition
\ref{cond1} (d), we may assume that $\sigma'$ centralizes $I(w)$.
According to the Katz--Gabber Lemma, there exists a unique
$I(w)$-cover $Z_i\to X_i$ which is isomorphic to $W_i\to X_i$ in an
\'etale neighborhood of $w$ and tamely ramified above $x_i\in X_i$. In
fact, $Z_i$ is a component of $g^{-1}(X_i)$, where $g:Z\to X$ is the
auxiliary cover associated to $f$. The automorphism $\sigma'|_{W_i}$
induces an isomorphism $\sigma'_{Z_i}$ of $Z_i$ centralizing $I(w)$.
Recall that there exists an isomorphism $Z_i\cong\PP^1_k$ such that
the action of $I(w)$ on $Z_i$ is generated by the translation
$z\mapsto z+1$ and the reflection $z\mapsto -z$. Using this
identification it is easy to see that $\sigma'_{Z_i}=\Id_{Z_i}$, hence
$\sigma|_{W_i}=\tau|_{W_i}$. It follows that $\tau\in G$ centralizes
the action of $D(W)=G$ on $W$. But $G$ has trivial center by Condition
\ref{cond2} (e), so $\sigma|_{W_i}=\Id_{W_i}$.

We have shown that $\sigma$ is the identity on all components of $Y$
except possibly on those that lie above $X_0$. Therefore, $\sigma$
restricts to an automorphism $\sigma|_E:E\iso E$ of $E$ which fixes
all points where $E$ meets another component of $Y$. Recall that the
set of points where $E$ meets another component is the image of a
$\Gamma_2(N')$-structure $\phi':A'=\ZZ/2N'\times\ZZ/2\to E(k)$, with
$N'=N/p$. In particular, $\sigma$ fixes $0\in E$ and is thus an
isomorphism of $E$ as generalized elliptic curve (see the proof of
Proposition \ref{monoprop} above). 

We can now finish the proof of Proposition \ref{innermonoprop} by
applying the classification of automorphism groups of generalized
elliptic curves, see \cite{DelRap}. Assume that $N'>1$. Using that
$\phi'(A')$ meets every component of $E$ and contains points of order
$>2$ we conclude that $\sigma|_E=\Id_E$. Therefore,
$\chi\bad(\gamma)=\chi\adm(\gamma)=1$, so $\gamma=1$. 

Suppose $N'=1$. Then $\sigma|_E\in\gen{[-1]}$, where $[-1]:E\iso E$ is
the canonical involution of $E$. Therefore, $\chi\bad(\gamma)=\pm 1$
and $\chi\adm(\gamma)=1$. It remains to be shown that $[-1]:E\iso E$
lifts to an element $\sigma\in\Aut_k(Y)$ inducing $\gamma\in\Gamma$
with $\chi\bad(\gamma)=-1$. We can set $\sigma|_{W_i}:=\Id_{W_i}$ for
$W_i:=f^{-1}(X_i)$, $i=1,\ldots,4$, and define $\sigma|_{\tau(E)}$ as
the canonical involution on the generalized elliptic curve $\tau(E)$,
for all $\tau\in G$. This completes the proof.  \Endproof

%\bibliographystyle{abbrv} \bibliography{../hurwitz}

%\end{document}

%% file: sec5.tex
%\documentclass{article}
%\usepackage{def}
%\begin{document}

\section{Applications to good reduction}

In this section we apply the results obtained in the previous sections
to questions of good reduction of Galois covers. We extend Raynaud's
criterion for good reduction to our situation (Theorem \ref{ramthm}).
Since we are in a very special situation, we get a somewhat sharper bound.
For covers with $4$ branch points, this result is not useful to
produce covers with good reduction, in practice. The rigidity method,
which can be used to construct covers over fields with low
ramification, hardly ever works for $4$ branch points. 

However, our results on the reduction of the Hurwitz space allow us to
compute the number of covers with good reduction, for given type and
position of the branch points (Theorem \ref{goodredthm}). The rough
idea is this. In characteristic $0$, the structure of the Hurwitz
space $H$ as a cover of the $\lambda$-line is known, and has a nice
description in terms of the braid action. Our Reduction Theorem
describes the structure of $\hb\otimes\FF_p$. The {\em Cusp Principle}
links these two results.  It states that a cusp of $\hb$ has good
reduction if and only if it corresponds to an admissible cover with
prime-to-$p$ ramification over the singular point.

\subsection{Good and bad reduction}

\subsubsection{} \label{situation1}

Let $G$ be a finite group and $\Cl=(C_1,C_2,C_3,C_4)$ be a $4$-tuple
of conjugacy classes of $G$. Let $K_0$ be a field of characteristic
$0$ such that the individual conjugacy classes $C_i$ are rational over
$K_0$, i.e.\ $\QQ(\Cl)\subset K_0$. We fix an algebraic closure
$\Kb_0$ of $K_0$. We choose an element $\lambda\in K_0-\{0,1\}$ and
define $\Cov$ as the set of isomorphism classes of $G$-covers
$f:Y\to\PP^1$, defined over $\Kb_0$, of type
$(\Cl;0,1,\infty,\lambda)$. In other words, $\Cov=\pi^{-1}(\lambda)$,
where
\begin{equation} \label{pieq}
  \pi:H_4\inn(\Cl) \To \PP^1_\lambda-\{0,1,\infty\}
\end{equation}
is the natural map from the inner Hurwitz space of $G$-covers of type
$\Cl$ to the $\lambda$-line, and we see $\lambda$ as a $\Kb_0$-rational
point on $\PP^1_\lambda$. 

Since the domain of definition of $H_4\inn(\Cl)$ is contained in $K_0$
(by assumption), we obtain a natural action of $\Gal(\Kb_0/K_0)$ on
$\Cov$. In more concrete terms, this action is given as follows. For
$\sigma\in\Gal(\Kb_0/K_0)$ and $f\in\Cov$, we can form the twisted
cover $\op{\sigma}{f}:\op{\sigma}{Y}\to\PP^1$ by applying $\sigma$ to
the coefficients of the equations defining $f$ and the action of $G$
on $Y$. To each $f\in\Cov$ we can associate the {\em field of moduli}
of $f$ (relative to $K_0$), i.e.\ the fixed field of all
$\sigma\in\Gal(\Kb_0/K_0)$ such that $\op{\sigma}{f}\cong
f$. Equivalently, $K$ is the field of rationality of the point on
$H_4\inn(\Cl)$ corresponding to $f$. If the center of $G$ is trivial
then $f$ has a unique model $f_K:Y_K\to\PP^1_K$ over $K$. See e.g.\
\cite{FriedVoe91}.

\subsubsection{} \label{situation2}

In the situation of \ref{situation1}, we will now make the following
assumptions. The field $K_0$ is complete with respect to a discrete
valuation $v$. We denote by $R_0$ the corresponding valuation ring.
The residue field $k_0$ of $v$ is assumed to be algebraically closed
of characteristic $p$, where $p$ is an odd prime number. We assume
that Conditions \ref{cond1} and \ref{cond2} hold, with respect to the
class vector $\Cl$ and the prime $p$.  Finally, we assume that
\begin{equation} \label{lambdaeq}
  \lambda \not\equiv 0,1,\infty  \pmod{v}.
\end{equation}
For a $G$-cover $f:Y\to\PP^1$ in $\Cov$ we can ask whether it has good
or bad reduction at $v$, in the sense of Section \ref{reduction}
(under Condition \zgl{lambdaeq}, good reduction is equivalent to
admissible reduction). If a cover $f\in\Cov$ has good reduction, its
reduction $f_k:Y_k\to\PP^1_k$ is a $G$-cover over the field $k$ of
positive characteristic, of type $(\Cl;0,1,\infty,\bar{\lambda})$
($\bar{\lambda}\in k$ denotes the residue of $\lambda\in K_0$). By a
theorem of Grothendieck, all $G$-covers over $k$ of type
$(\Cl;0,1,\infty,\bar{\lambda})$ arise as the reduction of a unique
$G$-cover $f\in\Cov$ with good reduction. This motivates the
following question.

\begin{question} \label{goodredquest}
  How many covers $f\in\Cov$ have good reduction at $v$?
\end{question}

Theorem \ref{goodredthm} below answers this question in an explicit
way. Its proof relies on the Reduction Theorem \ref{redthm} and the
{\em Cusp Principle}, Proposition \ref{cuspprinciple}.

\subsubsection{}

Let $\hb=\hb_4\inn(\Cl)$ be the completion of the Hurwitz space
$H=H_4\inn(\Cl)$, see Section \ref{complete}. We understand that $\hb$
is defined over $\Lambda:=\O_{\QQ(\Cl),\p}$. Here $\QQ(\Cl)\subset
K_0$ is as in Section \ref{morevariants} and $\p$ is the prime ideal
corresponding to the restriction of $v$ to $\QQ(\Cl)$. By Theorem
\ref{redthm} (i), $\hb$ is a normal scheme of dimension $2$, proper
and flat over $\Lambda$. Let $f\in\Cov$ and be $K$ its field of moduli
(relative to $K_0$). Since $K/K_0$ is finite, $v$ extends uniquely to
$K$. We denote by $R$ the corresponding valuation ring. Since $\hb$ is
proper over $\Lambda$, the morphism $\Spec K\to H$ corresponding to
$f$ extends uniquely to a morphism $\phi:\Spec R\to\hb$, giving rise
to a $k$-rational point $s:\Spec k\to\hb$. By the definition of
$\hb\bad$, $f$ has good (resp.\ bad) reduction if and only if
$s\not\in\hb\bad$ (resp.\ $s\in\hb\bad$).

We can be more precise. Since we assume the center of $G$ to be
trivial, $f$ descents to a unique $G$-cover $f_K:Y_K\to\PP^1_K$
defined over $K$.  Let $K'/K$ be the minimal extension of $K$ over
which $f_K\otimes K'$ has a stable model
$f_{0,R'}:Y_{R'}\to\PP^1_{R'}$ over the valuation ring $R'\subset K'$
(see Section \ref{reduction}). We denote by $f_{0,k}:Y_k\to\PP^1_k$
the special fiber of $f_{0,R'}$. The morphism $f_{0,k}$ (together with
the induced marks on $Y_k$ and $\PP^1_k$ and the $G$-action on $Y_k$)
is an object of the Hurwitz stack $\Hb$ associated to $\hb$ and
corresponds to the $k$-point $s:\Spec k\to\hb$. The stable model
$f_{0,R'}$ is a deformation of $f_{0,k}$. Hence it corresponds to a
unique morphism
\begin{equation} \label{equiveq1}
     R_{Y_k} \To R' 
\end{equation}
of $W(k)$-algebras, where $R_{Y_k}$ is the universal deformation ring
of $f_{0,k}$, see Section \ref{precise1}. The extension $K'/K$ is
Galois, and its Galois group is a subgroup of $\Gamma$, the universal
monodromy group of $f_{0,k}$, see Section \ref{precise2}. The morphism
\zgl{equiveq1} is equivariant with respect to the injection
$\Gal(K'/K)\inj\Gamma$. The morphism
\begin{equation}\label{equiveq2}
   \Od_{\hb,s} = R_{Y_k}^{\Gamma} \To R 
\end{equation}
obtained from \zgl{equiveq1} by taking invariants corresponds to the
morphism $\phi:\Spec R \to\hb$. If $f$ has good reduction, then
$\Gamma=1$ and $K'=K$. On the other hand, if $f$ has bad reduction, it
has modular reduction of level $N$, where $N$ is an integer strictly
divisible by $p$, see Proposition \ref{modularprop}.

We say that $\lambda\in K_0$ is {\em ordinary} (resp.\ {\em
supersingular}) if the elliptic curve $y^2=x(x-1)(x-\lambda)$ has
ordinary (resp.\ supersingular) reduction modulo $v$. By
\cite{Hartshorne}, Corollary IV.4.22, $\lambda$ is supersingular if and
only if $h_p(\lambda)\equiv 0 \pmod{v}$, where $h_p(X)\in\ZZ[X]$ is an
explicit polynomial of degree $(p-1)/2$. By a theorem of Igusa,
$h_p(X)$ is separable modulo $p$, i.e.\ there are exactly $(p-1)/2$
supersingular values $\bar{\lambda}$.

\begin{thm} \label{ramthm}
  Assume that $f\in\Cov$ has bad reduction of level $N$. 
  Let $K$ be the field of moduli of $f$, relative to $K_0$,
  and denote by $e$ the ramification index of $p$ in $K$. Then
  \begin{equation}\label{badpropeq}
         e \;\geq\; \left\{\begin{array}{cll}
                       \;(p-1)/2 \quad & \text{if}\;\; & N=p \\
                         p-1     \quad & \text{if}\;\; & N>p. \\
                    \end{array}\right.  
  \end{equation}
  Moreover, if $\lambda$ is supersingular then the inequality
  \zgl{badpropeq} is strict.
\end{thm}

\proof As explained above, the cover $f\in\Cov$ gives rise to a local
ring homomorphism $R_{Y_k}\to R'$, where $R_{Y_k}$ is the universal
deformation ring of the reduction of $f$ and $R'/R$ is the minimal
extension over which $f$ has a stable model. We denote by $\pi\in R'$
the image of the element of $R_{Y_k}$ with the same name, given by
Theorem \ref{precisethm}. Clearly, $v(\pi)>0$. Since $\pi^{p-1}|p$
(Theorem \ref{precisethm} (iii) and (iv)), the ramification index of
$p$ in $K'$ is at least $p-1$. But $\Gal(K'/K)\inj\Gamma$, so the
inequality \zgl{badpropeq} follows from Proposition
\ref{innermonoprop}. If $\lambda$ is supersingular then
$p=\pi^{p-1}u$, where $v(u)>0$, showing that the inequality
\zgl{badpropeq} is strict.  \Endproof

\begin{rem} \label{ramrem}
  The inequality $e\geq(p-1)/2$ can also be deduced from
  \cite{Raynaud98}, Corollaire 4.2.5 and Th\'eor\`eme 5.1.1.
\end{rem}

\subsection{The Hurwitz space as cover of the $\lambda$-line}

\subsubsection{The Hurwitz classification and braid action}
\label{braid}

Let us for the moment consider the morphism $\pi$ in \zgl{pieq} as an
unramified cover of Riemann surfaces. We choose a point
$\lambda_0\in(-\infty,0)$ on the negative real line and a point
$x_0\in\{x\in\CC\mid \Im x>0\}$ on the upper half plane. Let
$\gamma_i$ be the unique element of
$\pi_1(\PP^1_{\CC}-\{0,1,\infty,\lambda_0\},x_0)$ represented by a
simple loop which crosses the real line exactly twice, turning around
the $i$th point in the list $(0,1,\infty,\lambda_0)$, in
counterclockwise orientation. The group
$\pi_1(\PP^1_{\CC}-\{0,1,\infty,\lambda_0\},x_0)$ is generated by the
$\gamma_i$, $i=1,\ldots,4$, with the only relation
$\prod_i\gamma_i=1$. With these choices made, there is a canonical
bijection
\begin{equation} \label{Nielseneq}
  \pi^{-1}(\lambda_0) \;\cong\; \Ni_4\inn(\Cl) \;:=\; 
        \{\;\g=(g_1,\ldots,g_4) \mid 
              G=\gen{g_i},\;g_i\in C_i,\;\prod_i g_i=1\;\}/G,
\end{equation} 
(here $G$ acts on the set of tuples $\g$ by diagonal conjugation).

The cover $\pi$ induces an action of
$\Pi:=\pi_1(\PP^1-\{0,1,\infty\},\lambda_0)$ on $\pi^{-1}(\lambda_0)$,
hence on $\Ni_4\inn(\Cl)$ via \zgl{Nielseneq}.  To describe this
action explicitly, we denote by $\HH_4$ the {\em Hurwitz braid group}
on $4$ strings, with generators $Q_1,Q_2,Q_3$ and relations
\cite{Fried95}, (3.1.a-c). We identify $\HH_4$ with the fundamental
group of $\U_4:=\{\x\subset\PP^1(\CC) \mid |\x|=4\}$, with base point
$\x_0=\{0,1,\infty,\lambda_0\}$. Thus we obtain an embedding
\begin{equation} \label{hurwemb}
  \Pi:=\pi_1(\PP^1-\{0,1,\infty\},\lambda_0) \;\inj\; \HH_4.
\end{equation}
The elements
\begin{equation} \label{aeq}
    a_0       := Q_3 Q_2 Q_1^2 Q_2^{-1} Q_3^{-1}, \quad
    a_1       := Q_3 Q_2^2 Q_3^{-1}, \quad
    a_\infty  := Q_3^2,
\end{equation}
lie in the image of \zgl{hurwemb} and define standard generators of
$\Pi$. In particular, $a_0a_1a_\infty=1$, and $a_w$ is represented by
a simple loop around $w$, for $w\in\{0,1,\infty\}$. One can check
using the formula \cite{Fried95}, (3.1.d), that the induced action of
$\Pi$ on $\Ni_4\inn(\Cl)$ is given by
\begin{equation} \label{hurwop} 
   \renewcommand{\arraystretch}{1.7}
   [g_1,g_2,g_3,g_4]\, a_w \;=\;
   \left\{\begin{array}{lll}
      \;\;[g_1^\gamma,g_2,g_3,g_4^\gamma],\quad & \gamma=g_4g_1,\qquad 
               &  \text{if}\;\;w=0, \\
      \;\;[g_1,g_2^\gamma,g_3^{[g_2^{-1},g_4^{-1}]},g_4^\gamma],\quad
               &  \gamma=g_2g_4,\qquad & \text{if}\;\;w=1, \\
      \;\;[g_1,g_2,g_3^\gamma,g_4^\gamma],\quad & \gamma=g_3g_4,\qquad
               &  \text{if}\;\;w=\infty
   \end{array}\right.
\end{equation}
(we use the notation $g^\gamma=\gamma^{-1}g\gamma$).

\subsubsection{The cusps}

Let $\pib:\hb\to\PP^1_\lambda$ denote the canonical map, which extends
$\pi$. We say that a $\bar{\QQ}$-point $c$ on $\hb$ is a {\em cusp}
if $\bar{\pi}(c)\in\{0,1,\infty\}$.  We write $\Cp(\Cl,w)$ for the set
of cusps above $w\in\{0,1,\infty\}$. The cusps are the ramification
points of the finite, tamely ramified morphism $\bar{\pi}_{\QQ}$.
By Section \ref{braid} we can identify cusps with certain braid orbits:
\begin{equation} \label{cuspeq}
   \Cp(\Cl,w) \;\cong\; \Ni_4\inn(\Cl)/\gen{a_w}, \qquad 
       w\in\{0,1,\infty\}.
\end{equation}

Let us fix a cusp $c$ above $w\in\{0,1,\infty\}$, represented by a
class $[\g]\in\Ni_4\inn(\Cl)$. The conjugacy class of the element
$\gamma\in G$ associated to $[\g]$ in \zgl{hurwop} does only depend on
the orbit of $[\g]$ under the action of $a_w$, and is thus canonically
associated to the cusp $c$. We call $n:=\mathop{\rm ord}(\gamma)$ the
{\em order} of the cusp $c$. Note that the length of the $a_w$-orbit
of $[\g]$ divides $n$.

As a point on $\hb$, the cusp $c$ corresponds to an admissible cover
\[
    f_K:Y_K \To X_K
\]
between stably marked curves over a number field $K$, together with an
action of $G$ on $Y_K$. The bottom curve $X_K$ is singular, consisting
of two components $X_{1,K}, X_{2,K}\cong\PP^1_K$, intersecting in one
point. It is shown e.g.\ in \cite{diss}, Section 4.3.3, that the order
$n$ is the ramification index of $f$ above the singular point of
$X_K$.

We choose once and for all a valuation $\bar{v}$ of $\bar{\QQ}$
extending the valuation of $\QQ(\Cl)$ corresponding to the prime ideal
$\p$. Since $\hb$ is proper over $\Lambda$, the $\bar{\QQ}$-valued
point $c$ reduces (with respect to the valuation $\bar{v}$) to an
$\FFbp$-valued point $\bar{c}$ on $\hb$. We say that the cusp $c$ has
{\em good} (resp.\ {\em bad}) {\em reduction} if
$\bar{c}\not\in\hb\bad$ (resp.\ $\bar{c}\in\hb\bad$). If $c$ has bad
reduction then the point $\bar{c}\in\hb\bad$ corresponds to a bad
cover $f:Y\to X$ of modular type of level $N$, for some integer $N$
strictly divisible by $p$. We say that $c$ has bad reduction of level
$N$.

\begin{prop}[The Cusp Principle] \label{cuspprinciple}
  A cusp $c$ of order $n$ has bad reduction if and only if $p|n$. In
  this case, $n$ divides the level $N$.
\end{prop}

\proof After a finite extension of $K$, the admissible cover $f_K$
corresponding to $c$ extends to an $R$- object $f_R:Y_R\to X_R$ of the
Hurwitz stack $\Hb$, where $R$ is the valuation ring of $K$
corresponding to $\bar{v}|_K$. The special fiber $f:Y\to X$ of $f_R$
is the object of $\Hb$ which induces the $\FFbp$-point $\bar{c}$ on
$\hb$. We are exactly in the situation of Section \ref{boundarypts}.
The Cusp Principle follows from Proposition \ref{redbpts}.  \Endproof

\subsubsection{The bad components}

According to Theorem \ref{redthm}, $\hb\bad\otimes\FFbp$ is a smooth
projective curve over $\FFbp$. We pick a connected component
$W\subset\hb\bad\otimes\FFbp$. We call $W$ a {\em bad component}. The
geometric points of $W$ correspond to bad covers of modular type. It
is clear that the level $N$ of these covers is constant on $W$.
Therefore, we may call $N$ the level of $W$. We let $\eta$ be the
generic point of $W$, and we denote by $m$ the multiplicity of $W$
inside $\hb\bad$. By definition, $m$ is the length of the local
ring of $\eta$ on $\hb\otimes\FFbp$. 

\begin{prop} \label{badprop1}
  If $N=p$ then $m=(p-1)/2$. Otherwise, $m=p-1$.
\end{prop}

\proof This follows directly from Corollary \ref{precisecor3} (iii) and
Proposition \ref{innermonoprop}.  \Endproof

By Corallary \ref{precisecor3} (i), the natural map $W\to
\PP^1_\lambda\otimes\FFbp$ is inseparable, with inseparability degree
$p$. Moreover, the induced map $W^{(p)}\to\PP^1_\lambda\otimes\FFbp$
is tamely ramified in $0$, $1$, $\infty$, and \'etale everywhere
else. We are interested in describing it in more detail. 
We write $N=pN'$, and let $X_2(N')$ be the coarse moduli space for
generalized elliptic curves with $\Gamma_2(N')$-structure, see Section
\ref{dihedral}. We fix a bijection $\alpha:\{1,2,3,4\}\cong\ZZ/2\times\ZZ/2$;
this determines a finite map $X_2(N')\to\PP^1_\lambda$, by Proposition
\ref{G2Nprop2} (i). Since $N'$ is
prime-to-$p$, $X_2(N')\otimes\FFbp$ is a smooth projective curve over
$\FFbp$ and the map $X_2(N')\otimes\FFbp\to\PP^1_\lambda\otimes\FFbp$
is finite, tamely ramified in $0$, $1$, $\infty$ and \'etale
everywhere else.

\begin{prop} \label{badprop}
  There exists a finite map $h:W\to X_2(N')\otimes\FFbp$, compatible
  with the maps to $\PP^1_\lambda\otimes\FFbp$. This map is the
  composition of a purely inseparable map of degree $p$ and an \'etale
  map.
\end{prop}

\proof Let $s:\Spec k\to W$ be a geometric point, corresponding to a
cover $f_0:Y\to X_0$ over $k$. Following Section \ref{reductionthm},
we associate to $f_0$ a generalized elliptic curve $E$ over $k$,
together with a $\Gamma_2(N)$-structure $\phi$. We may assume that the
ordering of the branch points $x_1,\ldots,x_4\in X_0$ is compatible
with the bijection $\alpha:\{1,2,3,4\}\cong\ZZ/2\times\ZZ/2$ we have
chosen above. As in the proof of Proposition \ref{G2Ndefprop2}, we let
$E'':=E/\phi(\ZZ/p)$ and $\phi''$ be the induced
$\Gamma_2(N')$-structure. The pair $(E'',\phi'')$ gives rise to a
geometric point $s':\Spec k\to X_2(N')\otimes\FFbp$. One easily checks
that $(E'',\phi'')$ is unique up to isomorphism. In other words,
$s\mapsto h(s):=s'$ is a well defined map on geometric points. We
claim that this map is induced by a finite morphism $h:W\to
X_2(N')\otimes\FFbp$, as in the statement of the proposition. In order
to prove this claim, it suffices to show the following.  Let $R_s$
(resp.\ $R_{s'}$) be the complete local ring of $s$ on $W$ (resp.\ of
$s'$ on $X_2(N')\otimes\FFbp$). Then there exists a finite morphism
$h_s^*:R_{s'}\to R_s$ of local $\FFbp$-algebras, purely inseparable of
degree $p$, with the following property. Let $\Kb$ be an algebraic
closure of the fraction field of $R_s$, $\tilde{s}:\Spec\Kb\to W$ the
tautological point and $\tilde{s}':\Spec\Kb\to X_2(N')\otimes\FFbp$
the point induced by $h_s^*$. Then $h(\tilde{s})=\tilde{s}'$.

Following Section \ref{precise3}, we identify $R_s=\Od_{W,s}$ with
$(R_Y\bad)^\Gamma=k[[t]]^\Gamma$. By Proposition \ref{innermonoprop},
$\Gamma$ acts trivially on $k[[t]]$, so $R_s=k[[t]]$. In particular,
the point $\tilde{s}$ corresponds to the generic fiber of
$Y\univ\otimes_{R_Y}k[[t]]$. In the same way, we can identify $R_{s'}$
with $R_{\phi''}^{\Gamma''}$, where $R_{\phi''}$ is the universal
deformation ring of $(E'',\phi'')$ and $\Gamma''$ the corresponding
monodromy group. By the proof of Proposition \ref{G2Ndefprop2}, we
have $R_{\phi''}=k[[t^p]]$.  Moreover, $\Gamma''$ is trivial, because
\[
      \Aut(E'',\phi'') \;=\; \left\{
       \begin{array}{cl}
          \ZZ/2 & \qquad\mbox{if $N'=1$}, \\
            1   & \qquad\mbox{otherwise.}
       \end{array}\right.
\]
If we define $h_s^*$ as the natural injection $k[[t^p]]\inj k[[t]]$,
the proposition follows.
\Endproof

\begin{rem}\label{badrem}
  One can show that the \'etale part $h^{(p)}:W^{(p)}\to
  X_2(N')\otimes\FFbp$ of $h$ is in fact an isomorphism. In this
  sense, bad components are ``modular''.
\end{rem}

\subsubsection{The number of covers with good reduction}

Let us denote by $\Cov\good$ the subset of $\Cov$ containing the
covers with good reduction. Define
\[
     d := |\,\Ni_4\inn(\Cl)\,|, \qquad
     d\bad := |\,\{\,[\g]\in\Ni_4\inn(\Cl) 
               \;\mid\; p|\mathop{\rm ord}(g_3g_4)\;\}\,|.
\]
We know from \zgl{Nielseneq} that $d=\deg\pib=|\Cov|$. Using the Cusp
Principle and the Reduction Theorem, we can show:

\begin{thm} \label{goodredthm}
  We have
  \[    \renewcommand{\arraystretch}{1.5}
     |\Cov\good| \;=\; \left\{\;\;
        \begin{array}{ll} 
           d- d\bad,\qquad & 
              \text{\rm if $\lambda$ is ordinary,} \\
           d-\frac{p+1}{p}d\bad,\qquad & 
              \text{\rm if $\lambda$ is supersingular.}
        \end{array}
     \right.
  \]
  In particular, if $\lambda$ is ordinary and
  $\Ni_4\inn(\Cl)\not=\emptyset$ then $\Cov\good\not=\emptyset$. 
\end{thm}

\proof According to Theorem \ref{redthm}, $\hb\good$ is a smooth curve
over $\FF_q$ and the natural map
$\pib\good:\hb\good\to\PP^1_\lambda\otimes\FF_q$ is finite. Let
$S:=\hb\good\cap\hb\bad$. By Theorem \ref{redthm} (iii), $S$ contains
exactly the points on $\hb\bad$ with a supersingular $\lambdab$-value.
By definition, we have $\hb\good-S=\hb\adm\otimes\FF_q$. It follows
that $\hb\good-S\to\PP_\lambda^1\otimes\FF_q$ is tamely ramified in
$0$, $1$ and $\infty$ and \'etale everywhere else. Comparing the
degrees of $\pib\good$ and of $\hb\adm\to\PP_\lambda^1$ above
$\infty$, we get
\begin{equation} \label{gooddegeq}
   \deg\pib\good \;=\; 
          \sum_c e_c 
               \;=\; d-d\bad.
\end{equation}
Here $c$ runs over the set of cusps above $\infty$ with good reduction,
and $e_c$ denotes the ramification index of $\pib$ in $c$ (which is
equal to the length of the $a_\infty$-orbit of $\Ni_4\inn(\Cl)$). The
second equality in \zgl{gooddegeq} is a consequence of the Cusp
Principle, Proposition \ref{cuspprinciple}.

In case $\lambda$ is ordinary, the statement of the theorem follows  
directly from \zgl{gooddegeq}. Assume that $\lambda$ is supersingular,
and let $S_\lambda$ be the set of points in $S$ above
$\bar{\lambda}$. For $s\in S_\lambda$, we denote by $m_s$ the
ramification index of $\pi\good$ in $s$. We obtain 
\begin{equation} \label{baddegeq1}
   |\,\Cov\good\,| = \deg\pib\good - \sum_{s\in S_\lambda} m_s.
\end{equation}
According to Corollary \ref{precisecor3} (ii) and (iii), $m_s$ equals
the multiplicity of the bad component $W_s$ meeting $\hb\good$ in $s$.
By Corollary \ref{precisecor3} (i), the natural map
$W_s\to\PP^1_\lambda\otimes\FFbp$ is the composition of a purely
inseparable map of degree $p$ and a map which is \'etale away from
$0$, $1$ and $\infty$.  Therefore,
\begin{equation} \label{baddegeq2}
   d \;=\; \deg\pib\good \;+\; p\cdot\sum_{s\in S_\lambda} m_s.
\end{equation}
The equations \zgl{gooddegeq}, \zgl{baddegeq1} and \zgl{baddegeq2}
together imply Theorem \ref{goodredthm} in the supersingular case.
\Endproof

%\bibliographystyle{abbrv}
%\bibliography{../hurwitz} 

%\end{document}

%% file: examples.tex
%\documentclass{article}
%\usepackage{def}
%\begin{document}

\subsection{Examples} \label{examples}

In this section, we explain how the results we obtained can be used to
compute the reduction of the Hurwitz space, in an explicit example. We
take $G=PSL_2(\ell)$, where $\ell\neq p$ is a prime such that $p$
exactly divides $|G|=(\ell^2-1)\ell/2$. Recall that $G$ has two
conjugacy classes of order $\ell$, which we denote by $\ell A$ and
$\ell B$. We take $\underline{C}=(\ell A, \ell A, \ell A, \ell A)$ or
$\underline{C}=(\ell A, \ell A, \ell B, \ell B)$. The normalizer of a
$p$-Sylow group of $G$ for $p\neq \ell$ is a dihedral group of order
$\ell+1$ or $\ell -1$, depending on whether $p|\ell +1$ or $p|\ell-1$,
\cite{Huppert}. Note that Condition \ref{cond2} is also satisfied. We
are interested in computing the reduction to characteristic $p$ of
$\bar{H}:= \bar{H}^{\rm in}_4(\underline{C})$.

We will explain the algorithm for computing the reduction of $\bar{H}$
in the special case $\ell=11$ and $\underline{C}=(\ell A, \ell A, \ell
B, \ell B)$. After that, we give a table with the reduction of
$\bar{H}$ for $\ell\leq 31$, to all odd primes $p\neq \ell$ exactly
dividing the order of $G$. 

Take $\ell=11$ and $\underline{C}=(\ell A, \ell A, \ell B, \ell B)$.
Since $|G|=660$, we know that $\bar{H}$ has good reduction to
characteristic $p\neq 2,3, 5, 11$. The reduction to characteristic 2
and 11 we cannot compute using our methods.  Let us first take $p=3$.
The normalizer of a 3-Sylow group is a dihedral group of order 12, so
the possible levels associated to the bad components are 3 and 6. To
decide which ones occur, we want to apply the Cusp Principle
\ref{cuspprinciple}. For this we need to know the order of the cusps
in characteristic zero.

Using the program {\em ho}, \cite{ho}, we compute the irreducible components
$Z$ of $\bar{H}\otimes \bar{\QQ}$, and the ramification indices of $Z\to
\PP^1_\lambda$, for each irreducible component.  Let us relate these
ramification indices to the order of the cusps. Let $s$ be a cusp of
$\bar{H}\otimes \bar{\QQ}$ defined over $K$, and let $f_K\!:Y_K\to X_K$ be
the corresponding admissible cover. Let $\tau$ be the unique singular
point of $X_K$ and $\rho$ a singular point of $Y_K$. Let $n$ be the
ramification index of $\rho$ and $G_1$ and $G_2$ the decomposition
groups of the two components of $Y_K$ passing through $\rho$.  Using
the description of the normalizers of elements in $G=PSL_2(\ell)$
given in \cite{Huppert}, Abschnitt II.8, one easily checks that the
ramification index $e$ of $s$ in $\bar{H}\otimes \CC\to \PP^1_\lambda$
is equal to $n$, unless $n=\ell$ and $G_i\simeq \ZZ/\ell$ for some
$i$, see.\ Section \ref{braid}. In particular, $p|n$ iff $p|e$.
\[
\begin{array}{|l|r|r|r|}
\hline
\multicolumn{4}{|c|}{\ell=11,\quad  \underline{C}=(\ell A, \ell A, \ell
B, \ell B)}\\
\hline
\mbox{ramification}& \mbox{deg} & g &\mbox{num}\\
\hline
2^1;-;1^2&2&0&1\\
\hline
2^2 6^2;-;1^4 3^4&16&1&1\\
\hline
2^3 1^5 6^2 5^2;-;5^211^13^4& 33& 2& 1\\
\hline
\end{array}
\]

The notation is as follows. Each row corresponds to an irreducible
component; the last entry of each row gives the number of isomorphic
components. The first entry gives the ramification of $\bar{H}\otimes
\CC\to \PP^1_\lambda$ over $0,1,\infty$. Here $2^1$ means one
ramification point of order 2, and $1^2$ means two ramification points
of order 1. A ``$-$'' indicated that over this point, the ramification
indices are the same as over the previous point. The next entries give
the genus of the components and its degree over the $\lambda$-line.

The Cusp Principle implies that the first component $W_1$ has good
reduction to characteristic 3, since 3 does not divide the order of
any of the cusps. The second and third component, which we will denote
by $W_2$ and $W_3$, have bad reduction to characteristic 3. Since both
these components have a cusp of order $n=6$ and the order of the cusp
divides the level of the bad component it reduces to, we see that in
both cases there will be a bad component of level $N=6$. A bad
component of level 6 is a cover of $X_2(N')\otimes\bar{\FF}_3$, purely
inseparable of degree $3$, with $N'=N/p=2$, see Remark \ref{badrem}.
The curve $X_2(2)\otimes\bar{\FF}_3$ is a cover of degree 2 of the
$\lambda$-line, branched at 0 and 1 and unbranched at $\infty$.  In
Proposition \ref{badprop1} we computed the multiplicity of the bad
components. A bad component of level 6 has multiplicity $p-1=2$. We
conclude that the reduction of the irreducible components $W_2$ and
$W_3$ each have one bad component, and it is of level 6.

Now let us have a look at the good components. The good and the bad
components intersect over the supersingular $\lambda$'s. In
characteristic 3, there is only one supersingular $\lambda$, namely
$\lambda=-1 \pmod{3}$. This means that the good and bad components
meet in two points. Since the multiplicity of the bad component is
two, the good components will be ramified of order two in these
intersection points. From this we can compute the number of covers
with good reduction, for each value of $\lambda$.
\[
  |{\rm Cov}(W_2,\lambda)\good| \;=\;
\left\{\begin{array}{ll} 
4& \mbox{ if } \lambda\not\equiv -1 \pmod{3},\\
0&  \mbox{ if } \lambda\equiv -1 \pmod{3}.
\end{array}
\right.
\]
\[
  |{\rm Cov}(W_3,\lambda)\good| \;=\;
\left\{\begin{array}{ll} 
21& \mbox{ if } \lambda\not\equiv -1 \pmod{3},\\
17&  \mbox{ if } \lambda\equiv -1 \pmod{3}.
\end{array}
\right.
\]

In general we are not able to calculate the number of good components.
However, since the degree of $W_2$ is sufficiently small, we can
describe what its good part looks like. As remarked before, the degree
of $W_2^{\rm good}$ over $\PP^1_\lambda$ is 4 and it is ramified over
$-1$ at two points of order two. Outside the supersingular
$\lambda$'s, the ramification is as in characteristic zero. So over
$0$ and 1, there are two ramification points of order two, and over
$\infty$ it is unramified. From this it follows that $W_2^{\rm good}$
is connected.

\bigskip\noindent Now let us have a look at $p=5$. In this case, the
components $W_1$ and $W_2$ both have good reduction, since 5 does not
divide the order of any of the cusps. The component $W_3$ has bad
reduction. Note that the normalizer of a 5-Sylow group of $G$ is a
dihedral group of order 10, so the only possibility for the level is
5. A bad component of level 5 is a cover of
$X_2(1)\otimes\bar{\FF}_5$, purely inseparable of degree $5$. The
curve $X_2(1)\otimes\bar{\FF}_5$ is isomorphic to the
$\lambda$-line. It has multiplicity $(p-1)/2=2$. In characteristic 5,
there are two supersingular $\lambda$-values: the primitive sixth
roots of unity. In these points the good part will have an ``extra''
ramification of order two. So as before we are able to compute the
number of covers with good reduction.

\[
   |{\rm Cov}(W_3,\lambda)\good| \;=\;
\left\{\begin{array}{ll} 
23& \mbox{ if } \lambda\not\equiv \zeta_6, \zeta_6^5 \pmod{5},\\
21&  \mbox{ if } \lambda\equiv  \zeta_6, \zeta_6^5\pmod{5}.
\end{array}
\right.
\]

The results are summarized in the following lemma.

\begin{lem}
  Let $W_1, W_2, W_3$ be the three irreducible components of $H^{\rm
    in}_4(PSL_2(11))\otimes \bar{\QQ}$, as described above.
\begin{itemize}
\item[(a)] Then $W_1$ has good reduction to characteristic $p\neq 2,
  11$.
\item[(b)] The component $W_2$ has good reduction to characteristic
  $p\neq 2,3,11$. In characteristic 3, it as two irreducible
  components: a bad component of level 6 and a good component. The
  degree of the good component over the $\lambda$-line is 4.
\item[(c)] The component $W_3$ has good reduction to characteristic
  $p\neq 2,3,5,11.$ In characteristic 3, there is one bad component,
  of level 6, and the degree of the good part over the $\lambda$-line
  is 21. In characteristic 5, there is one bad component, of level 5,
  and the degree of the good part over the $\lambda$-line is 23.
\end{itemize} 
\end{lem}

The Hurwitz space might have irreducible components of large
degree having good reduction at many primes. For example, take
$\ell=31$ and $\underline{C}=(\ell A, \ell A, \ell B, \ell B)$. The
Hurwitz space in characteristic zero has an irreducible component of
genus 37, whose degree over the $\lambda$-line is 128. The cusps of
this component all have ramification index a power of 2. We conclude
from the Cups Principle that this component has good reduction to all
primes $p\neq 2, 31$.  For details, see the table below.

\begin{exa}[Raynaud]\label{Raynaudexa}  
  Let ${H}$ be the inner Hurwitz space parameterizing Galois covers of
  $\PP^1$ with Galois group $A_5$ which are branched at four points of
  order 3. Let $\bar{H}$ be its completion, over $\ZZ_{(5)}$. (Raynaud
  considered the absolute Hurwitz space; it is easy to make the
  adaption to that case.)  In characteristic zero, $\bar{H}\otimes
  \bar{\QQ}$ is connected and has degree 18 over the $\lambda$-line.
  Over $0,1,\infty$ this cover has ramification $3^2 5^2 1^2$. We
  conclude as above that the reduction of the Hurwitz space to
  characteristic 5 has one bad component of level 5. So the good
  degree is 8. If $\lambda$ reduces to a supersingular value, there
  are 6 covers with good reduction.

  The above example was presented by Raynaud in his talk in
  Oberwolfach, June 1997. In the problem session of the same
  conference he proposed the exercise of computing the number of
  covers with good reduction to characteristic 3 for $G=A_5$ and
  ramification of order 5. The answer to this exercise appears in the
  first rows of the table below.
\end{exa}

\bigskip\noindent The following table describes the reduction of the
Hurwitz space $\bar{H}:=\bar{H}^{\rm in}_4(\underline{C})$, where
$G=PSL_2(\ell)$ and $\underline{C}$ is either $(\ell A, \ell A, \ell
B, \ell B)$ or $(\ell A, \ell A, \ell A, \ell A)$. Every row
corresponds to an isomorphism class of irreducible components in
$\bar{H}$; the entry ``num'' gives the number of isomorphic
components. The first column gives the $\ell$, the second column gives
the class vector. The third column gives the ramification of
$\bar{H}\to \PP^1_\lambda$ in characteristic zero over $0,1,\infty$.
Here $a^b$ means: $b$ ramification points of order $a$ and $-$ means:
the same as the previous point. The entries ``deg'' and $g$ give the
degree over the $\lambda$-line and the genus of the component (in
characteristic zero). The last three entries describe the reduction
for odd primes different from $\ell$ which exactly divide the order of
$G$. No statement is made for other primes. A dash means: the
component has good reduction to all such primes. Each prime $p$ such
that the component has bad reduction to characteristic $p$, is listed
on a separate row. Under ``bad components'', for each prime, all the
bad components are listed. The last entry gives the degree of the good
part over the $\lambda$-line. This is the number of covers with good
reduction, for $\lambda$ ordinary. The number of covers with good
reduction for supersingular $\lambda$ can be computed by Theorem
\ref{goodredthm}. A component has good reduction to characteristic $p$
for odd primes $p\neq \ell$ which are not listed and which strictly
divide the order of $G$. The prime $\ell=13$ is missing from the table
because our computer refused to run {\em ho} for $PSL_2(13)$. The prime
$\ell=17$ is missing because there are no primes $p\neq 17$ which
exactly divide the order of $PSL_2(17)$.

\renewcommand{\arraystretch}{1.065}
\[\begin{array}{|r|l|l|r|r|r|l|l|r|}
\hline
\ell&\mbox{ Ni}& \mbox{ramification} & \mbox{deg}& g&\mbox{num}&p &\mbox{ bad comp}&
\mbox{gdeg}\\
\hline
\hline
5&AABB&2^1;-;1^2&2&0&1&-&-&-\\
\hline
5&AABB&3^1 1^2;3^1 2^1;-& 5& 0& 1& 3& 1\times N=3& 2\\
\hline
5&AAAA& 5^13^1 1^2;-;-& 10&0& 1& 3& 1\times N=3& 7\\ 
\hline
7&AABB&4^2;-;1^42^2& 8& 0& 1&  -& -&-\\
\hline
7&AABB& 1^3 4^2 3^1;-;7^1 3^1 2^2& 14& 0& 1& 3& 1\times N=3& 11\\
\hline 
7&AAAA&1^2; 2^2;-& 2&0&3& -&-&- \\
\hline
7&AAAA&2^23^1;-;-& 7&0&1&3&  1\times N=3&4 \\
\hline
11&AABB& 2^1;-;1^2&2&0&1&-& -&-\\
\hline
11&AABB& 2^26^2;-;1^43^4&16&1&1&3& 1\times N=6 &4 \\
\hline
11&AABB&2^3 1^5 6^2 5^2;-; 5^211^1 3^4& 33& 2& 1& 3& 1\times N=6 &
21\\
&&&&&&5& 1\times N=5& 23\\

\hline 
11& AAAA & 1^13^1;-;-& 4& 0& 4& 3& 1\times N=3 & 1\\
\hline
11& AAAA & 3^4 5^2; -; -& 22& 3& 1& 3& 4\times N=3&10 \\
&&&&&& 5& 1\times N=5 &12\\ 
\hline
19&AABB&  2^1;-;1^2&2&0&1&-& -&-\\
\hline
19&AABB&2^4 10^4; -; 1^8 5^8& 48& 9& 1& 5& 1\times N=10&16\\
\hline
19& AABB& 2^5 1^9 3^3 9^3 10^4 ;-; 5^8 9^3 3^3 19^1& 95& 17& 1& 5&
1\times N=10&55\\
\hline
19& AAAA& 1^2 5^2;-;-& 12& 1& 4& 5& 1\times N=5&2\\
\hline
19& AAAA& 3^3 5^8 9^3; -; -& 76& 18& 1& 5& 4\times N=5&36\\
\hline
23& AABB& 4^2; -; 1^42^2& 8& 0&1& -&-&-\\
\hline
23& AABB& 4^4 12^4;-;1^86^43^82^4& 64& 13& 1& 3& 1\times N=12&16\\
\hline
23 &AABB& 1^{11} 4^6 11^5 12^4; -; & 138& 32& 1 &
3& 1\times N=12&90\\
&& 6^4 11^5 3^8 2^6 23^1 &&&& 11&1\times N=11&83\\
\hline
23&AAAA & 1^2;2^1;-& 2&0&3&-&-&-\\
\hline 
23&AAAA& 3^1 1^1;-;-& 4& 0&4&3&1\times N=3&1\\
\hline
23&AAAA& 2^2 6^2; -; 3^41^4& 16& 1& 2& 3& 1\times N=6&4\\
\hline
23&AAAA&2^6 3^8 11^5 6^4;-;-& 115& 24&1&3& 3\times N=6, 2\times N=3&67\\
&&&&&& 11&1\times N=11& 60\\
\hline
29& AABB&  2^1;-;1^2&2&0&1&-& -&-\\
\hline
29& AABB& 2^6 14^6 ;-; 1^{12} 7^{12}& 96& 25& 1& 7& 1\times N=14&12\\
\hline
 29& AABB& 5^6 3^5 2^7 14^6 15^4;-;& 203& 54& 1&
3& 1\times N=15, 1\times N=3&128\\
&& 1^4 15^4 7^{12} 5^6 3^5&&&& 5& 1\times N=15, 1\times N=5&113\\
&&&&&& 7& 1\times N=14&119\\
\hline
29&AAAA& 1^37^3;-;-& 24& 4& 4& 7& 1\times N=7&3\\
\hline
29& AAAA& 1^{14} 3^5 7^{12} 5^{6} 15^4 29^1; -; -& 232& 54& 1& 3& 1\times N=15, 1\times N=3&157\\
&&&&&& 5& 1\times N=15, 1\times N=5&142\\
&&&&&& 7& 2\times N=7&148\\ 
\hline
31&AABB& 16^8;-;1^{16} 2^8 4^8 8^8& 128 & 37 & 1 &-&-&-\\
\hline
31& AABB& 1^{15} 3^5 16^8 5^6 15^4 ;-; & 248 & 67 & 1& 3& 1\times N=15, 1\times N=3&173\\
&& 5^6 15^4 31^{11} 2^8 4^8 3^5
8^8 &&&& 5& 1\times N=15, 1\times N=5&158\\
\hline
31&AAAA&  2^1;-;1^2&2&0&3&-&-&-\\
\hline
31&AAAA& 1^4 2^2; 4^2;-& 8 & 0 & 3& -&-&-\\
\hline
31& AAAA& 1^8 2^4 4^4; 8^4; -& 32& 0 & 3& -&-&-\\
\hline
31& AAAA& 3^5 2^8 15^4 5^6 4^8 8^8;-;-& 217 &51& 1& 3&  1\times N=15, 1\times N=3&142\\
&&&&&& 5& 1\times N=15, 1\times N=5&127\\
\hline
\end{array}
\]

%\bibliographystyle{abbrv}
%\bibliography{hurwitz}

%\end{document}

%%% Local Variables: 
%%% mode: latex
%%% TeX-master: t
%%% End: 